\input amstex
\NoBlackBoxes

\define\phom{ \, \underset{p}\to\sim \, }

\documentstyle{amsppt}
\pagewidth{6.5in}
\topmatter
\title 
Refined Transfer
\endtitle
\author
J. C. Becker
\endauthor
\endtopmatter
\document

\subhead 1. Introduction \endsubhead The index theorem of Dwyer, Weiss, and Williams \cite{4} asserts that two tranfer like maps which arise
in algebraic $K$-theory are homotopic. Let $E\downarrow B$ be a smooth bundle and let $A$ denote Walhausen's algebraic $K$-thory of retractive spaces
\cite{8}. The usual transfer defines a map $T_A:|B|@>>>A(E)$, where $|B|$ denotes the realization of the singular complex of $B$. On the other hand there is the characteristic $\chi_A:|B|@>>>A(E)$. The index theorem asserts a preferred homotopy
$$T_A \ \sim \ \chi_A \tag 1.1$$
There results a definition of torsion as a secondary index \cite{1}. Let $R$ be a ring and $V$ a local system of $R$-modules over $E$.
There is the map $ch:A(E)@>>>K(R)$ which assigns to a retractive space its chain complex with coefficients in $V$. Consider
$$|B|@>T_A \ \sim \ \chi_A>>A(E)@>ch>>K(R) \tag 1.2$$
One may have a trivialization $ch\chi_A\ \sim \ *$. For example, if the fiber of $E$ is acyclic rel. $V$. On the other hand a non zero 
vertical vector field on $E$ defines a trivialization $T_A \ \sim \ *$, so a second trivialization of $ch\chi_A$ is given by
$ch\chi_A \ \sim \ ch T_A \ \sim \ *$. Gluing these together defines a torsion map
$$\Cal T:|B|@>>>\Omega K(R) \tag 1.3$$

Let Vec$_B(E)$ denote the set of homotopy classes of non zero vertical vector fields on $E$. Since there is a preferred path $T_A \ \sim \ \chi_A$,
one should have (for fixed $ch\chi_A \ \sim \ *$) a refined torsion map
$$\Cal T:\text{Vec}_B(E)\times |B|@>>>\Omega K(R)  \tag1.4$$ 
which restricts to the refined torsion of Turaev \cite{6, 7} when $B$ is a point. Toward this end we will consider a refined notion of transfer. 
Recall that transfer is characterized by certain axioms \cite{3}. The refinement specifies preferred paths  by which the axioms are satisfied.

\bigskip
Let $x_0, x_1$ be points in $X$. By a {\it path class} from $x_0$
to $x_1$ we mean a component of the space of paths from $x_0$ to $x_1$. If $h$ is a path from $x_0$ to $x_1$, we will
write $h:x_0\phom x_1$ to denote the path class containing $h$

\bigskip
We will want to consider relative bundles in both the horizontal and vertical direction. Given $E_1\downarrow B$,
let $E_0<E$ denote that $E_0$ is a codimension $0$ subbundle whose frontier is neatly embedded in $E$.

\bigskip
Let $L(B,A)=|B|/|A|$. There is the natural pairing $\omega_L:L\wedge L@>>>L$.

\bigskip
Let $F$ be a homology theory from pairs of spaces to pointed spaces. We assume $F$ is a loop functor, and is equipped with a pairing
$F\wedge L@>>>F$, and augmentation $\varepsilon:L@>>>F$.

\bigskip
A refined transfer $T$ assigns to each \ $E_{01}=(E_0<E_1)\downarrow(B,A)$ \ a map
$$T(E_{01}):L(B,A)\longrightarrow F(E_1,E_1\downarrow A)$$
The axioms imply the following. 

\bigskip
\item{\bf I.} Let $E$ be diffeomorphic to an orthogonal disk bundle, and let $\Delta:B@>>>E-\dot E$ be a section. A path class 
$$\bold d(E,\Delta):F(\Delta)\varepsilon(B,A) \ \phom \ T(E,E\downarrow A)$$ 

\bigskip
\item{\bf II.} For each excision \ $(\widetilde f, f):(E_0<E_1)@>>>(\overline E_0< \overline E_1)$ \ a path class
$$T(\widetilde f,f):F(\tilde f)T(E_{01}) \ \phom \ T(\overline E_{01})L(f)$$

\bigskip
\item{\bf III.} For each
$$E_{012}=(E_0<E_1<E_2)\downarrow (B,A)$$
an additivity path class
$$T(E_{012}): T(E_{02}) \ \phom \ F(i_{12}) T(E_{01})+T(E_{12})$$
Here $+$ is the loop product in $F$, and $i_{12}:E_1@>>>E_2$ is the inclusion map.

\bigskip
\item{\bf IV.} For a pair of spaces $(X,Y)$ consider
$$
\CD
F(E_1,E_1\downarrow A)\wedge L(X,Y)@> \omega_F>> F((E_1,E_1\downarrow A)\times (X,Y))\\
@AA T(E_{01})\wedge 1 A       @AA T(E_{01}\times (X,Y)) A\\
L(B,A)\wedge L(X,Y)@> \omega_L >> L((B,A)\times (X,Y))\\
\endCD
$$
A  multiplicativity path
$$T(E_{01},(X,Y)):\omega_F (T(E_{01})\wedge 1) \ \phom \ T(E_{01}\times (X,Y))\omega_L$$

\bigskip
\proclaim{(1.5) Theorem} A homology theory $F$ has a refined transfer. If \ $T,\overline{T}:L@>>> F$ \ are refined transfers, there is
a unique function $\Lambda$ which assigns to each bundle $E_{01}$ 
a path class 
$$\Lambda:T(E_{01}) \ \phom \ \overline{T}(E_{01})$$
such that $\Lambda$ commutes with the path classes {\bf I}-{\bf IV}.
\endproclaim

Let us briefly indicate how theorem (1.5) leads to refined torsion map as in (1.4). The theorem cannot be applied directly since $A$ is not a homology theory. We have
$$P@>a>>A@>b>>A^s$$
where $A^s$ is the stabilization of $A$, and $P$ is Waldhausen's partition  space \cite{9}. Fundamental results of Waldhausen \cite{8, 9,10}state 
firstly, that $A^s$ is a homology theory, and secondly, that $ba$ is a homotopy equivalence. A non zero vector field defines a trivialization 
$T_A^s \ \sim \  *$. We then have a trivialization of $\chi_A^s$  by  
$$\chi_A^s \overset {\Lambda}\to {\sim } \ T_A^s \ \sim  \ * $$
The uniqueness part of theorem (1.5) is used to show that this path class does not depend on the choices involved in constructing $T_A$.
Since $\chi_A$ factors through $P$, this trivialization can be pulled back to $P$, then pushed forward to $A$, obtaining a trivialization of
$\chi_A$.

\bigskip
\subhead 2. Homology \endsubhead  Let $\Cal T^2$ denote the category of pairs of compact spaces, and $\Cal T_0$ the
category of compact pointed spaces. Suppose that $F:\Cal T_0@>>>\Cal T_0$ is a functor such that for every one point space $*$, \ $F(*)$ is a one point space, and 
$$\omega:F(X)\wedge |Y| \longrightarrow F(X\wedge Y)$$
is a natural transformation such that:
\bigskip
\itemitem{(2.1)}  $F(X)=G(X)\wedge |S^0|@>\omega >>
F(X\wedge S^0)=F(X)$ is the identity. 
\medskip
\itemitem{(2.2)} Let $j:|Y|\wedge |Z|\rightarrow |Y\wedge Z|$ denote the canonical map. The following commutes.
$$
\CD
F(X)\wedge |Y|\wedge |Z|@> \omega \wedge 1>> F(X\wedge Y)\wedge |Z|\\
@VV 1\wedge j V  @VV\omega V\\
F(X)\wedge |Y\wedge Z| @>\omega >> F(X\wedge Y\wedge Z).\\
\endCD
$$

\bigskip
\noindent
We will refer to $(F,\omega)$ as a {\it multiplicative} functor. 

Note that $(F,\omega)$ defines a functor from homotopies to homotopies: Given $h_0, h_1:X\rightarrow Y$, and a homotopy 
$H:X\wedge I^+\rightarrow Y$ from $h_0$ to $h_1$, the composite
$$F(X)\wedge I^+ @>1\wedge\lambda>> F(X)\wedge |I^+| @>\omega>>
F(X\wedge I^+)@> F(H)>> F(Y),$$ 
where $\lambda$ is the canonical map, is a homotopy from $F(h_0)$ to $F(h_1)$. We will write \newline
$\widehat F(H):F(X)\wedge I^+\longrightarrow F(Y)$ \ 
to denote this homotopy.  
Note also that if $(F,\omega)$ is a multiplicative functor we have a multiplicative functor 
$\Omega(F,\omega)=(\Omega F,\omega')$, where
$\omega'$ is given by
$$\Omega (F(X))\wedge |Y|\longrightarrow \Omega(F(X)\wedge |Y|)\overset{\Omega(\omega)}\to\longrightarrow 
\Omega(F(X\wedge Y).$$

A multiplicative functor $(F,\omega_F)$ is a {\it homology theory} if 
for each cofibration \newline $A@>i>> X@>q>> X/A$, \ 
the natural map
$$F(A)\longrightarrow \text{Fiber}\left(F(X)@>F(q)>> F(X/A)\right)$$
is a weak equivalence. 

\bigskip
The {\it stabilization} of $(F,\omega)$ is the multiplicative functor $(F^s,\omega^s)$, where
$$F^s(X)=\lim_{k}  \ \text{Map}\, (|S^1|^k;F(X\wedge S^k), \tag 2.3$$
the limit taken with respect to
$$\text{Map}\, (|S^1|^k;F(X\wedge S^k)) @> \text{\underbar{\hskip.1in}}\wedge 1>>
\text{Map}\, (|S^1|^{k+1};F(X\wedge S^k)\wedge |S^1|)@>\omega_*>>
\text{Map}\, (|S^1|^{k+1};F(X\wedge S^{k+1})).$$
The pairing 
$$\omega^s:F^s(X)\wedge |Y| \longrightarrow F^s(X\wedge Y) \tag 2.4$$ 
is the limit of 
$$\align
\text{Map}\, (|S^1|^k;&F(X\wedge S^k))\wedge |Y|@>>>\text{Map}\, (|S^1|^k;F(X\wedge S^k)\wedge |Y|)@>\omega_*>>
\\&\text{Map}\, (|S^1|^k;F(X\wedge S^k\wedge Y)@>F(\alpha)_*>>
\text{Map}\, (|S^1|^k;F(X\wedge Y\wedge S^k),\\
\endalign
$$
where $\alpha:X\wedge S^k\wedge Y\rightarrow X\wedge Y\wedge S^k$ is the interchange map. Then $(F,\omega)$ is 
{\it stable} if the natural transformation $(F,\omega)\rightarrow (F^s,\omega^s)$ is a weak equivalence. Stability of
$(F,\omega)$ is equivalent the property:
\bigskip
\itemitem{(2.5)} \ $F(X)\rightarrow\text{Map}(|S^1|;F(X\wedge S^1))$ \ is a weak equivalence.

\bigskip
\noindent
The property is clearly sufficient. Since $(F^s,\omega^s)$ has the property for any 
$(F,\omega)$, it is also necessary. 

A proof of the following is straightforward.

\proclaim{(2.6) Theorem}  A homology theory is stable. 
\endproclaim

It will be convenient to rewrite (2.5) in more familiar form. Let $\alpha:(I,\dot I)\rightarrow (S^1,\infty)$ be a 
relative homeomorphism such that $\alpha:(0,1)\rightarrow \Bbb R$ is increasing. We have 
$\beta:(I,\dot I)\rightarrow (|S^1|,*)$ by $(I,\dot I)@>>>(|I|,|\dot I|)@>|\alpha|>>(|S^1|,*)$. Define
$$\eta:F(X)\longrightarrow \Omega F(X\wedge S^1) \tag 2.7 $$
by
$$F(X)@>>>\text{Map}(|S^1|,F(X\wedge S^1)@>\beta^*>>\Omega F(X\wedge S^1)$$

An {\it augmentation} of $(F,\omega)$ is a natural transformation
$$\varepsilon:|X|\longrightarrow F(X)$$
such that
$$
\CD
F(X)\wedge |Y|@> \omega \wedge 1>> F(X\wedge Y)\\
@AA \varepsilon\wedge 1 A  @AA\varepsilon  A\\
|X|\wedge |Y| @>>> |X\wedge Y|.\\
\endCD 
$$
commutes.

\medskip
In what follows,  a homology theory $(F,\omega_F)$,  will be assumed firstly to be a loop functor, \newline
$(F,\omega_F)=\Omega (\widetilde F,\omega_{\widetilde F})$ for specified  $(\widetilde F,\omega_{\widetilde F})$, and secondly to be augmented.

\medskip
We now consider functors $\Cal T^2@>>>\Cal T_0$. Define $(L,\omega_L)$
by 
$$L(X,A)=|X|/|A|, \tag 2.8$$
$$\omega_L:L(X,A)\wedge L(Y,B)\rightarrow L((X,A)\times (Y,B)), \tag 2.9$$
the natural map $|X|/|A|\wedge |Y|/|B|\rightarrow |X\times Y|/|X\times B\cup A\times Y|$. 
A {\it multiplicative} functor is a pair $(F,\omega)$, where $F:\Cal T^2@>>>\Cal T_0$  is a functor and
$$\omega:F(X,A)\wedge L(Y,B)\longrightarrow F((X,A)\times (Y,B))$$
is a natural transformation such that

\medskip
\itemitem{(2.10) \, }  $F(X,A)=F(X,A)\wedge L(*)@>\omega >>F((X,A)\times *)=F(X,A)$ is the identity. 

\medskip
\itemitem{(2.11) \, } The following commutes.
$$
\CD
F(X,A)\wedge L(Y,B)\wedge L(Z,D)@> \omega \wedge 1>> F((X,A)\times (Y,B))\wedge L(Z,D)\\
@VV 1\wedge \omega_L V  @VV\omega_G V\\
F(X,A)\wedge L((Y,B)\times (Z,D)) @>\omega>> F((X,A)\times (Y,B)\times (Z,D)).\\
\endCD
$$

An {\it augmentation} of $(F,\omega)$ is a natural transformation
$$\varepsilon:L(X,A)\longrightarrow F(X,A)$$
such that
$$
\CD
F(X,A)\wedge L(Y,B)@> \omega \wedge 1>> ((X,A)\times (Y,B))\\
@AA \varepsilon\wedge 1 A  @AA\varepsilon  A\\
L(X,A)\wedge L(Y,B) @>\omega_L>> L((X,A)\times (Y,B)).\\
\endCD 
$$
commutes.

If $(F, \omega_F):\Cal T_0@>>>\Cal T_0$ is a multiplicative functor, we obtain a multiplicative functor 
$(G,\omega_G):\Cal T^2@>>>\Cal T_0$, by  $G(X,A)=F(X/A)$, and $\omega_G$  by
$$F(X/A)\wedge |Y|/|B|@>>>F(X/A)\wedge |Y/B|@>>>F(X/A\wedge Y/B)$$ 
In this case we say that $(G,\omega_G)$ is a homology theory
if $(F,\omega_F)$ is a homology theory. 

\bigskip
\subhead 3. Simplicial Transformations \endsubhead Let us recall some basic terminology. For \newline $p=0,1, \dots$, \ let 
$$[p]=\{0,1, \dots,p\}$$
Let $\Delta$ denote the category with objects $[p]$ and morphisms the order preseving maps. We may also regard $[p]$ as a category with one morphism $(i,j):\{i\} @>>> \{j\}$ whenever $i\le j$. The functors $[q] @>>> [p]$ are then precisely the order preserving maps. 

Recall that a simplicial set (category, etc.) is a contravariant functor from $\Delta$ to the category of sets 
(categories, etc.). If $\Cal C$ is a category, the nerve $\Cal N(\Cal C)$ of $\Cal C$ is the simplicial set
$$
\align
[p] & \ \longrightarrow \ \text{Functor}([p],\Cal C)\\
(\alpha:[q]@>>>[p]) & \ \longrightarrow \ (\alpha^*:\text{Functor}([p],\Cal C)) @>>> \text{Functor}([q],\Cal C))\\
\endalign
$$
\medskip
Let $\Cal C$ be a category and $F,\,G$ functors from $\Cal C$ to the category of pointed spaces. Let $\Cal S$ be a 
subcomplex of $\Cal N(\Cal C)$. By a {\it simplicial transformation} \ $T:F\rightarrow G$ \ on \ $\Cal S\subset\Cal N(\Cal C)$, \
we mean a function which assigns to each $p$-simplex $C$ in $\Cal S$,
$$C=\left(C_0@>  \ c_0^1 \ >> \ \dots\dots  \ @>c_{(p-1)}^p>>C_p\right)  \tag 3.1$$
a map
$$T(C):F(C_0)\wedge \Delta_p^+\longrightarrow G(C_p), \tag 3.2$$
such that for $\alpha:[q]@>>>[p]$, the following commutes
$$
\CD
F(C_0)\wedge \Delta_q^+ @>1\wedge\alpha>> F(C_0)\wedge\Delta_p^+ @>T(C)>> G(C_p)\\
@VV F(c_0^{\alpha(0)})\wedge 1 V    @.      @VV 1 V\\
F(C_{\alpha(0)})\wedge \Delta_q^+ @> T(\alpha^*(C))>> G(C_{\alpha(q)}) @>G(c_{\alpha(q)}^p)>> G(C_p)\\
\endCD \tag 3.3
$$
Equivalently, the face and degeneracy maps satisfy: 
$$
\align
&T(C)(1\wedge d_0)=T(d_0(C))(F(c_0^1)\wedge 1)\\
&T(C)(1\wedge d_j)=T(d_j(C)),\quad 1\le j\le p-1  \tag 3.4\\
&T(C)(1\wedge d_p)=G(c_{(p-1)}^p)T(d_p(C))\\
&T(C)(1\wedge s_j)=T(s_j(C)),\quad 0\le j\le p\\
\endalign 
$$

For $p=0$, we have
$$T(C):F(C)\longrightarrow G(C)$$
These are the maps of interest. The remaining maps specify naturality up to homotopy. Looking explicitly at the cases $p=1,2$, we see that for $p=1$,  
$$T(C):F(C_0)\wedge \Delta_1^+\longrightarrow G(C_1)$$ 
is a path validating homotopy commutativity of the diagram
$$
\CD
G(C_0)@>G(c_0^1)>>G(C_1)\\
@AA T(C_0)A     @AA  T(C_1) A\\
F(C_0) @> F(c_0^1)>> F(C_1)\\
\endCD
$$

For $p=2$, the diagram 
$$
\CD
G(C_0)@>G(c_0^1)>>G(C_1)@>G(c_1^2)>>G(C_2)\\
@AA T(C_0)A  @AA T(C_1)A  @AA T(C_2) A\\
F(C_0) @> F(c_0^1)>> F(C_1) @> F(c_1^2)>> F(C_2)\\
\endCD
$$
is homotopy commutative in two ways. Then $T(C)$ provides a contractible choice of paths between the two, namely any path of the form
$$F(C_0)\wedge \Delta_1^+\wedge I^+@>1\times h>> F(C_0)\wedge \Delta_2^+ @>T(\sigma)>> G(C_2),$$
where $h:\Delta_1\times I\rightarrow\Delta_2$ is a path from $d_2\cdot d_0$ to $d_1$. The linear path would be a canonical choice. 

Note that a natural tranformation $T:F\rightarrow G$ on $\Cal C$ defines a simplicial transformation  by defining $T(C)$ to be
$$F(C_0)\wedge \Delta_p^+ @>{\text{proj.}}>> F(C_0)@>F(c_0^p)>>F(C_p)@>
T(C_p)>>G(C_p)$$

\bigskip
Let $T:F\rightarrow G$ be a simplicial transformation on $\Cal S@>\lambda>>\Cal N(\Cal C)$. If $S:G\rightarrow G'$ is a natural transformation on $\Cal C$, we have a simplicial transformation 
$$ST:F\rightarrow G' \qquad\text{ on }\qquad \Cal S@>\lambda>>\Cal N(\Cal C)$$ 
defined by \ $ST(\sigma)=S(C_p)T(\sigma)$.

\bigskip
Similarly, if $S:F'\rightarrow F$ is a natural transformation on $\Cal C$, we have a simplicial transformation 
$$TS:F'\rightarrow G \qquad\text{ on }\qquad \Cal S@>\lambda>>\Cal N(\Cal S)$$ 
defined by \ $TS(\sigma)=T(\sigma)(S(C_0)\wedge 1)$. 

\bigskip
Given a pointed space $X$ and functor $F$, let $F\wedge X$ denote the functor $C\rightarrow F(C)\wedge X$. Let $T:F\wedge X@>>>G$
be a simplicial transformation. If $A$ is a pointed subspace of $X$, let $T|F\wedge A:F\wedge A@>>>G$ denote $T$ composed with the
natural transformation given by inclusion $F\wedge A@>>>F\wedge X$. Also, for $x\in X$, let $T_x:F@>>>G$ denote $T$ composed with
$\theta_x:F(C)@>>>F(C)\wedge X$, \  $\theta_x( - )=( - ,x)$.

Let $S,T:F@>>>G$ be simplicial transformations. A {\it simplicial homotopy} from $S$ to $T$ is a simplicial transformation $H:F\wedge I^+\rightarrow G$ such that $H_0=T$, $H_1=S$. 

Let us record the following elementary fact.
\proclaim{(3.5) Lemma} Given \ $S^{p-1}@>f>>Y_0@>\eta>>Y_1$ \ where $\eta$ is a homotopy equivalence, suppose $\eta f$ has an
extension $\overline{ \, \eta f}:D^p@>>>Y_1$. Then $f$ has an extension $\overline f:D^p@>>>Y_0$ such that 
$\eta\overline f \ \sim \ \overline{ \, \eta f}$ rel. $S^{p-1}$.
\endproclaim

\bigskip
Suppose given 
$$
\CD
@.                  G_0\\
@.               @V \eta VV\\
F\wedge (I^k)^+  @>\tau>> G_1\\
\endCD
$$
where $\tau$ is a simplicial transformation on $\Cal S\subset \Cal N(\Cal C)$, and $\eta$ is a natural transformation such that 
$\eta(C)$ is a homotopy equivalence for each object $C\in\Cal C$. Let $S:F\wedge \partial(I^k)^+@>>>G_0$ be a simplicial transformation 
and $K:F\wedge \partial(I^k)^+\wedge I^+@>>>G_1$ a simplicial homotopy such that $K_0=\eta S$, \  $K_1=\tau|F\wedge \partial(I^k)^+$.

\proclaim{(3.6) Theorem} There is a simplicial tranformation $T:F\wedge (I^k)^+@>>>G_0$ such that $T|F\wedge \partial(I^k)^+=S$, 
and simplicial homotopy $H:F\wedge (I^k)^+\wedge I^+@>>>G_1$  such that $H_0=\eta T$, $H_1=\tau$, and $H|F\wedge \partial(I^k)^+\wedge I^+=K$.
\endproclaim

\demo{Proof} Let $\Cal S^{(p)}$ denote the $p$-skeleton of $\Cal S$. Suppose on $\Cal S^{(p-1)}$ we have
$$T^{p-1}:F\wedge (I^k)^+@>>>G_0 \qquad\qquad H^{p-1}:F\wedge (I^k)^+\wedge I^+@>>>G_1$$
satisfying the conditions. Let 
$$C=(C_0@>c_{01}>> \dots @>c_{(p-1)p}>>C_p)$$
be non degenerate. Let $A=I^k\times\Delta_p$. We have
$$\widetilde T^p(C):F(C_0)\wedge \dot A^+\longrightarrow G_0(C_p)$$
where \ $\widetilde T^p(C)|F(C_0)\wedge \partial(I^k)^+\wedge \Delta_p^+=S(C)$, \ and $\widetilde T^p(C)|F(C_0)\wedge (I^k)^+\wedge \dot\Delta_p^+$ \ 
are defined from $T^{(p-1)}(C)$ by the conditions ( \ \ ). Define
$$P(C):F(C_0)\wedge (\dot A \times I\cup A\times \{1\})^+\longrightarrow G_1(C_p)$$
so that that its restriction to:
$$
\align
&F(C_0)\wedge (\partial(I^k)\times \Delta_p\times I)^+ \text{ \ is \ } K(C)\\
&F(C_0)\wedge (I^k\times \dot\Delta_p\times I)^+ \text{ \ is defined by (3.4) from \ } H^{p-1}(C)\\
&F(C_0)\wedge (I^k\times \Delta_p\times \{1\})^+ \text{ \ is \ } \tau(C)\\
\endalign
$$
Since $(\dot A \times I\cup A\times \{1\},\dot A\times\{0\})$ is homeomorphic to $(D^{k+p},S^{k+p-1})$, there is, by the lemma,
$$Q(C):F(C_0)\wedge (\dot A \times I\cup A\times \{1\})^+\longrightarrow G_0(C_p)$$
such that $Q_0(C)=\widetilde T^p(C)$, \ and 
$$L(C):F(C_0)\wedge (\dot A \times I\cup A\times \{1\}\times I)^+\longrightarrow G_1(C_p)$$
such that \ $L_0(C)=P(C)$, \  $L_1(C)=\eta Q(C)$, \ and $L$ is rel. $F(C_0)\wedge (\dot A\times \{0\})^+$. Let

$$W=(\dot A \times I\cup A\times \{1\}\times I)/(\sim),\qquad   (\dot a,0,s) \ \sim \ (\dot a,0,t), \quad  0\le s,t\le 1$$
Then $L(C)$ defines
$$\overline L(C):F(C_0)\wedge W^+\longrightarrow G_1(C_p)$$
Let \ $\psi:A@>>>\dot A \times I\cup A\times \{1\}$ \ be such that $\psi(\dot a)=(\dot a,0)$, \ and let $\phi:A\times I@>>>W$ \ 
be such that $\phi(y)=(y,0)$, \ $y\in \dot A \times I\cup A\times \{1\}$, \ and \ $\phi(a,0)=(\psi(a),1)$. Define
$$T^p(C)=Q(C)(1\wedge \psi^+):F(C_0)\wedge A^+\longrightarrow G_0(C_p)$$
$$H^P(C)=\overline L(C)(1\wedge\phi^+):F(C_0)\wedge (A\times I)^+\longrightarrow G_1(C_p)$$
This defines $T^p(C)$, $H^p(C)$, if $C$ is a non degenerate $p$-simplex. If $C$ is a non degenerate $q$-simplex, $q<p$, let
$T^p(C)=T^q(C)$, $H^p(C)=H^q(C)$. If $C$ is a degenerate simplex in $\Cal S^{(p)}$, define $T^p(C)$, $H^p(C)$ by condition (3.3).
We now have simplicial maps $T^p$, $H^p$ on $\Cal S^{(p)}$ which satisfy the conditions. Finally, define $T$, $H$ on $\Cal S$ 
by \ $T(C)=T^p(C)$, and $H(C)=H^p(C)$, \ if $C\in \Cal S^{(p)}$.
\enddemo

\bigskip
Let $T:F\longrightarrow G$ be a simplicial transformation on $\Cal C$. Let Ar$(\Cal C)$ denote the arrow category of $\Cal C$,
and extend $F$, $G$ to functors on Ar$(\Cal C)$ by
$$F(C_0@>c_{01}>>C_{1})=F(C_0)\qquad\qquad G(C_0@>c_{01}>>C_{1})=G(C_1)$$
Then $T$ defines a simplicial transformation
$$\widehat T:F\wedge I^+\longrightarrow G\qquad \text{on}\qquad \text{Ar}(\Cal C) \tag 3.7$$
as follows. Let $\bold C\in \Cal N_p(\text{Ar}(\Cal C))$,
$$\bold C=
\CD
C_{01}@>c_{01}^{11}>> \dots @>c_{(p-1)1}^{p1}>> C_{p1}\\
@AA c_{00}^{01} A @. @AA c_{p0}^{p1} A\\
C_{00}@>c_{00}^{10}>> \dots @>c_{(p-1)0}^{p0}>> C_{p0}\\
\endCD
$$
Let $I\times \Delta_p$ have the standard simplicial decomposition. The $(p+1)$ simplexes 
are
$$\Delta^j=((0,e_0),\dots,(0,e_j),(1,e_j),\dots,(1,e_p))$$
There is the corresponding element of $\Cal N_{p+1}(\Cal C)$,
$$\bold C^j=(C_{00}@>>>\dots @>>>C_{j0}@>>> C_{j1}@>>>\dots @>>> C_{p1})$$
Define \ $\widehat T(\bold C)(\_,s,t)=T(\bold C^j)(\_,s,t)$, \ if \ $(s,t)\in \Delta^j$. \ 
Let \ $\bold C_j=(C_{0j}@>>>\dots @>>>C_{pj})$, \ $j=0,1$. We see that 
$$\widehat T_0(\bold C)=G(c_{p0}^{p1})T(\bold C_0),\qquad T_1(\bold C)=T(\bold C_1)F(c_{00}^{01}) \tag 3.8$$

\bigskip
\subhead 4. Multiplicative structures \endsubhead Let $\Cal Pt.\subset \Cal T^2$ denote the subcategory of 
pairs $(*,\Phi)$ where $*$ is a one point space. Let $\Cal C$ be a category. A {\it multiplicative structure} for $\Cal C$ 
is a  functor 
$$\times:\Cal C\times \Cal T^2\longrightarrow \Cal C, \hskip.7in (C,(X,Y))\longrightarrow C\times (X,Y) \tag 4.1$$
together with natural equivalences

\medskip
\noindent
(i) From   
$$\text{proj.}:\Cal C\times\Cal Pt.\longrightarrow  \Cal C$$
to 
$$\times:\Cal C\times\Cal Pt.\longrightarrow  \Cal C$$

\medskip
\noindent
(ii) From 
$$\Cal C\times \Cal T^2\times \Cal T^2\longrightarrow \Cal C ,\qquad  (C,(X,Y),(X',Y'))\longrightarrow C\times ((X,Y)) \times (X',Y'))$$ 
to
$$(C,(X,Y),(X',Y'))\longrightarrow (C\times (X,Y))\times (X',Y')$$
It is assumed that the following commutes.
$$
\CD
C @>>> C\times *\\
@VVV   @VV V\\
C \times (*\times *') @>>> (C \times *)\times *'\\
\endCD
$$
 
\medskip
Let $\Cal C$ be a multiplicative category. A  multiplicative functor \ $(F,\omega_F):\Cal C@>>>\Cal T_0$
consists of a functor $F:\Cal C@>>>\Cal T_0$, together with a natural transformation
$$\omega_F:F(C)\wedge L((X,Y)\longrightarrow F(C\times(X,Y))$$
such that the following commute
\medskip
$$
\CD 
F(C)@> 1 >> F(C)\\
@VVV   @VVV\\
F(C)\wedge L(*)@>\omega_F>>F(C\times *)\\
\endCD \tag i
$$
\medskip
$$
\CD
F(C)\wedge L(X,Y)\wedge L(X',Y') @> \omega_F \wedge 1>> F(C \times (X,Y))\wedge L(X',Y')\\
@VV 1\wedge\omega_L V    @VV \omega_F V\\
F(C)\wedge L((X,Y)\times(X',Y')) @.   F((C\times (X,Y))\times(X',Y')\\
@VV \omega_F V   @VVV\\
F(C\times ((X,Y)\times(X',Y'))) @> = >>   F(C\times ((X,Y))\times(X',Y')))\\ 
\endCD \tag ii
$$

\bigskip
Let $(F,\omega_F):\Cal C@>>>\Cal T_0$ be a multiplicative functor. We have the associated functors
$$
\align
F^{\times}:\Cal C\times\Cal T^2@>>>\Cal T_0,\qquad \qquad &F^{\times}(C,(X,Y))=F(C\times (X,Y)) \tag 4.2\\
F\wedge L:\Cal C\times\Cal T^2@>>>\Cal T_0,\qquad \qquad &F\wedge L(C,(X,Y))=F(C)\wedge L(X,Y))\\
\endalign
$$
and natural transformation
$$\omega_F:F\wedge L@>>>F^{\times} \tag 4.3$$
Suppose now that $(F,\omega_F)$ and $(G,\omega_G)$ are multiplicative functors $\Cal C@>>>\Cal T_0$, and
$$T:F\longrightarrow G$$
is a simplicial transformation on $\Cal C$. There is the associated simplicial transformation
$$T^{\times}:F^{\times}\longrightarrow G^{\times} \tag 4.4$$
on $\Cal C\times\Cal T^2$ defined as follows: Under the identification 
$\Cal N(\Cal C\times \Cal T^2)=\Cal N(\Cal C)\times \Cal N(\Cal T^2)$, a $p$-simplex of $\Cal N(\Cal C\times \Cal T^2)$
is a pair $(\bold C,(\bold X,\bold Y))$, where $\bold C\in \Cal N_p(\Cal C)$, \ $(\bold X,\bold Y)\in \Cal N_p(\Cal T^2)$.
Writing
$$\bold C=C_0@>>>\dots@>>>C_p, \qquad (\bold X,\bold Y)=(X_0,Y_0)@>>>\dots@>>>(X_p,Y_p)$$
let
$$\bold C\times (\bold X,\bold Y)=C_0\times(X_0,Y_0)@>>>\dots@>>>C_p\times(X_p,Y_p)$$
Define \ $T^{\times}(\bold C,(\bold X,\bold Y))=T(\bold C\times (\bold X,\bold Y))$.

\bigskip
Consider the diagram of simplicial transformations on \, $\Cal C\times\Cal T^2$
$$
\CD
F\wedge L@>\omega_F>> F^{\times} \\
@VV T\wedge 1 V         @VV T^{\times} V\\
G\wedge L @>\omega_G>> G^{\times}\\
\endCD \tag 4.5
$$
A {\it multiplicative structure} for $T$ is a simplicial homotopy on \, $\Cal C\times\Cal T^2$
$$\mu_T:F\wedge L\wedge I^+\longrightarrow G^{\times}  \tag 4.6$$
such that \ $(\mu_T)_0=\omega_g(T\wedge 1)$, \ $(\mu_T)_1=T^{\times } \omega_F$. In addition, $\mu_T$ is assumed to have the normality properties described below.

\medskip\noindent
On \, $\Cal C\times \Cal Pt.$, \, the diagram (4.5) reduces to
$$
\CD
F@> 1 >> F \\
@VV T V   @VV T V\\
G @> 1 >> G\\
\endCD 
$$
We assume
\medskip\noindent 
{\bf I.} $\mu_T|\Cal C\times \Cal Pt.$ is the constant simplicial homotopy.

\medskip\noindent
We have functors 
$$L^{\times}((X,Y),(X',Y'))=L((X,Y)\times(X',Y'))$$
on $\Cal T^2\times\Cal T^2$,   
$$F^{\times\times}(C,(X,Y),(X',Y'))=F(C\times (X,Y)\times (X',Y'))$$
on $\Cal C\times\Cal T^2\times\Cal T^2$, and simplicial tranformation \ 
$T^{\times\times}:F^{\times\times}@>>>G^{\times\times}$,
$$T^{\times\times}(\bold C,(\bold X,\bold Y),(\bold X',\bold Y'))=
T(\bold C\times (\bold X,\bold Y)\times(\bold X',\bold Y'))$$
We have
$$\nu'_T:F\wedge L\wedge L\wedge (I^2)^+\longrightarrow G^{\times\times} \tag 4.7$$
with edge path
$$
\CD
F\wedge L\wedge L @>\omega_G\wedge 1 > \sim > F^{\times}\wedge L @>\omega_F >\sim > F^{\times\times}\\
@V\sim V T\wedge 1\wedge 1 V @V\sim V T^{\times}\wedge 1 V    @V\sim V T^{\times\times} V\\
G\wedge L\wedge L @>\omega_G\wedge 1 > \sim >G^{\times}\wedge L @>\omega_G >\sim >G^{\times\times}\\
\endCD
$$
The left square is filled in with $\mu_T\wedge 1$, and the right by $\mu'_T$, where
$$\mu'_T(\bold C,(\bold X,\bold Y),(\bold X',\bold Y'))=\mu_T(\bold C\times (\bold X,\bold Y),(\bold X',\bold Y'))$$
Secondly, we have
$$\nu{''}_T:F\wedge L\wedge L\wedge (I^2)^+\longrightarrow G^{\times\times} \tag 4.8$$
with edge path
$$
\CD
F\wedge L\wedge L @>1\wedge \omega_L  > \sim > F\wedge L^{\times} @>\omega_F >\sim > F^{\times\times}\\
@V\sim V T\wedge 1\wedge 1 V @V\sim V T\wedge 1 V    @V\sim V T^{\times\times} V\\
G\wedge L\wedge L @>1\wedge \omega_L > \sim >G\wedge L^{\times} @>\omega_G >\sim >G^{\times\times}\\
\endCD
$$
The left square commutes and the right is filled in by $\mu{''}_T$, where
$$\mu{''}_T(\bold C,(\bold X,\bold Y),(\bold X',\bold Y'))=\mu_T(\bold C, (\bold X,\bold Y)\times(\bold X',\bold Y'))$$
We assume associativity
\medskip\noindent  
{\bf II.} $\nu'_T$ and $\nu{''}_T$ are simplicially homotopic rel. $\partial(I^2)$.

\bigskip
\subhead 5. Permutation paths \endsubhead We suppose given a category $\Cal C$ and functors 
$L,F:\Cal C\rightarrow \Cal T_0$. Let us suppose further that $F$ has the form
$$\Cal C @>\eta>> \Cal T_0 @>F_0>>  \Cal T_0$$
where $F_0$ is a homology theory. Recall the assumption that $F_0$ is a loop functor, $F_0=\Omega G_0$. Setting $G=G_0\eta$, we have 
$F=\Omega G$. Let
$$T_1,\dots,  T_n:L\longrightarrow F $$
be simplicial tranformations on $\Cal S\subset \Cal N(\Cal C)$, and $\rho$  a permutation of $\{1,\dots,n\}$.  We note here that the 
given data determines a simplicial path
$$\rho:T_1+\cdots + T_n \ \sim \ T_{\rho(1)}+\cdots +T_{\rho(n)} \tag 5.1$$

\bigskip
Given $\sigma\in \Omega (X)$, let $\sigma^{(k)}\in \Omega (X)$, \ $1\le k\le n$, \ denote the loop \ 
$\sigma^{(k)}(t)=\sigma(nt-k+1)$, \ $(k-1)/n\le t\le k/n$, \ and \ $\sigma^{(k)}(t)=*$, \ otherwise.  For $f:Y@>>>X$, define 
$f^{(k)}:Y@>>>X$ by $f^{(k)}(y)=f(y)^{(k)}$.

\bigskip
For $X\in \Cal T_0$,  we have natural transformations on $\Cal T_0$
$$\bigvee_1^n G_0(X)@>\hskip.1in\phi_0\hskip.1in>>  G_0(\bigvee_1^n X)@>\hskip.1in\psi_0\hskip.1in>> \prod_1^n G_0(X) \tag 5.2$$
and $\psi_0\phi_0$ is the inclusion. Let
$$\Gamma_n F(C)= F_0(\bigvee_1^n \eta(C)) \tag 5.3$$
From (5.2), we have natural transformations on $\Cal C$
$$\Omega\bigvee_1^n G(C)@>\hskip.1in\phi\hskip.1in>>  \Gamma_nF(C)@>\hskip.1in\psi\hskip.1in>> \Omega\prod_1^n G(C)= 
\prod_1^n F(C)\tag 5.4$$
and $\psi\phi$ is the inclusion. Consider
$$
\CD
@.          \Gamma_n F\\
@.           @VV \psi V\\
L @>>>      \dsize{\prod}_1^n F\\
\endCD \tag 5.5
$$
Note that \ $(T_1^{(\rho(1))},\dots,T_n^{(\rho(n))})$  \ has image in $\dsize{\Omega\bigvee_1^n G(C)}$. Define
$$\langle T_1^{(\rho(1))},\dots,T_n^{(\rho(n))}\rangle=\phi (\psi\phi)^{-1}(T_1^{(\rho(1))},\dots,T_n^{(\rho(n))}) \tag 5.6$$
Then 
$$\psi\langle T_1^{(\rho(1))},\dots,T_n^{(\rho(n))}\rangle =(T_1^{(\rho(1))},\dots,T_n^{(\rho(n))})$$ 
There is the canonical homotopy
$$\widetilde H:(T_1^{(1)},\dots,T_n^{(n)}) \ \sim \ (T_1^{(\rho^{-1}(1))},\dots,T_n^{(\rho^{-1}(n))})$$
By the lifting propety (3.6), there is a simplicial homotopy
$$H:G\longrightarrow \Gamma_n F$$
such that 
$$H_0=\langle T_1^{(1)},\dots,T_n^{(n)}\rangle \qquad\qquad H_1=\langle T_1^{(\rho^{-1}(1))},\dots,T_n^{(\rho^{-1}(n))}\rangle$$
and \ $\psi H \ \sim \ \widetilde H$ \ rel. $\dot I$. Let 
$$\lambda:\bigvee_1^n \eta(C)@>>>\eta(C)$$
denote the folding map. We have 
$$F_0(\lambda)H:T_1+\dots +T_n \ \sim \ T_{\rho(1)}+\dots +T_{\rho(n)} \tag 5.7$$ 
Any two such are simplicially homotopic rel. $\dot I$.

\proclaim{5.8 Definition} A permutation path
$$\rho:T_1+\dots +T_n \ \sim \ T_{\rho(1)}+\dots +T_{\rho(n)}$$
is a simplicial path which is homotopic rel. $\dot I$ to one of the form (5.7). 
\endproclaim

\bigskip
\subhead 6. Transfer\endsubhead Suppose given for $\sigma\in\Omega X$, products
$$\prod_1^m \sigma^{\epsilon_j},\qquad \prod_1^n \sigma^{\delta_j}$$
where \ $\epsilon_j, \delta_j\in\{0,-1,1\}$, \ and \ $\sum\epsilon_j=\sum\delta_j$. Let $\gamma:I@>>>S^1$, \ $\gamma(t)=e^{2\pi i t}$.
Then $\sigma$ factors as $I@>\gamma>>S^1@>\widetilde\sigma>>$. By a {\it standard} path
$$\frak s:\prod_1^m \sigma^{\epsilon_j} \ \sim \ \prod_1^n \sigma^{\delta_j}$$
we mean one of the form $\widetilde\sigma h$, where \ $h:\prod_1^m \gamma^{\epsilon_j} \ \sim \ \prod_1^n \gamma^{\delta_j}$. The space of
standard paths is contractible.

\bigskip
Let
$$H=\{ \ x\in R^n\mid x_1\ge 0 \ \},\quad  H_{+}=\{ \ x\in H\mid x_2\ge 0 \ \},\quad H_{-}=\{ \ x\in H\mid x_2\le 0 \ \}$$
Let $M$ be a compact smooth manifold and $M_0$ a compact codimension zero submanifold. Suppose that for each point $x\in M_0$ there is an open neighborhood $U$ of $x$ and diffeomorphism $\phi:U\rightarrow \phi(U)\subset H$ such that
$$\phi(U\cap M_0)=\phi(U)\cap H_{-},\qquad \phi(U\cap \overline{M-M_0})=\phi(U)\cap H_{+}$$
It follows that the frontier fr$(M_0 \subset M)=M_0\cap\overline {M-M_0}$ is a neatly embedded smooth submanifold of $M$ \cite{5}, and $M_0$ inherits a smooth structure by straightening the angle at points of the boundary of fr$(M_0\subset M)$. In this case we shall write $M_0<M$.
\bigskip

If $M_0<M_1<M_2$, it need not be the case that $M_0<M_2$. Let us say that a diagram

$$
\matrix
M_{0q} & < & \dots & < & M_{pq}\\
\vee &&&& \vee\\
\vdots &&&& \vdots\\
\vee &&&& \vee\\
M_{00} & < & \dots & < & M_{p0}\\
\endmatrix \tag 6.1
$$
is {\it regular} if $M_{ij}<M_{i'j'}$ for all $(i,j)\le(i',j')$. Then every subdiagram is also regular. Hereafter we assume that diagrams are regular.
A diagram 
$$
\matrix
M_{01} & < &  M_{11}\\
\vee && \vee\\
M_{00} & <  & M_{10}\\
\endmatrix
$$
is {\it excisive} if
$$M_{00}=M_{10}\cap M_{01} \qquad\qquad M_{11}=M_{10}\cup M_{01}$$
Generally, a diagram as in (6.1) is excisive if every $2\times 2$ subdiagram is excisive.

\bigskip
Let $E$ be a fiberwise smooth bundle over $B$, and $E_0$ a subbundle. We write $E_0<E$ 
if for $b\in B$, we have $(E_0)_b<E_b$. Let $\Cal B^p$ denote the multiplicative category with objects \, $(E_0<\dots<E_p)\downarrow (B,A)$. The arrows are fiberwise smooth maps. Let
$$(\widetilde f,f):(E_0<E_1)\downarrow (B,A)\longrightarrow (\overline E_0<\overline E_1)\downarrow (\overline B,\overline A)$$
be a fiberwise smooth embedding. It is an {\it excision} if for $b\in B$, we have an excisive diagram

$$
\matrix
(\overline E_{01})_{f(b)} & < &  (\overline E_{11})_{f(b)}\\
\vee && \vee\\
\widetilde f((E_{00})_b) & <  & \widetilde f((E_{10})_b)\\
\endmatrix
$$

Let us say that a $q$-simplex in the nerve of $\Cal B^p$,
$$E_0@> \ \ (\widetilde f_{01},f_{01}) \ \ >> \dots @>(\widetilde f_{(p-1)p},f_{(p-1)p})>> E_q,\qquad\qquad 
E_j=(E_{0j}<\dots<E_{pj})\downarrow (B_j,A_j)$$
is excisive if for all $(i_0,j_0)\le (i_1,j_1)$, 
$$(E_{i_0j_0}<E_{i_1j_0})@> \ \ (\widetilde f_{j_0j_1},f_{j_0j_1}) >>(E_{i_0j_1}<E_{i_1j_1})$$
is excisive. These form a subcomplex $\Cal S^p$ of $\Cal N(\Cal B^p)$.

\bigskip
Let $F:\Cal T^2@>>>\Cal T_0$ be a homology theory. Recall that $L(X,A)=|X|/|A|$. We extend $L$ and 
$F$ to multiplicative functors on $\Cal B^p$ by
$$L((E_0<\dots<E_p)\downarrow (B,A))=L(B,A),\qquad\qquad F((E_0<\dots<E_p)\downarrow (B,A))=F(E_p,E_p\downarrow A)$$

\bigskip
Let $\Cal D$ denote the category whose objects are pairs $(E\downarrow (B,A),\Delta)$, where $E$ is diffeomorphic to an
orthogonal disk bundle, and $\Delta:B@>>>E-\dot E$ is a section. The arrows are fiberwise diffeomorphisms which preserve 
the section. Let
$$d(E,\Delta)=F(\Delta)\varepsilon(B,A):L(B,A)\longrightarrow F(E\downarrow(B,A) \tag 6.2$$

\bigskip
We will use the notation \ $E_{0\dots p}=E_0<\dots<E_p$.

\bigskip\noindent
{\bf (6.3) Definition.} A {\it transfer} for $F$ consists of multiplicative simplicial transformations
$$
\align
& T:L\wedge(I^{p-1})^+\longrightarrow F \quad \text{on}\quad \Cal S^p\subset\Cal N(\Cal B^p),\quad 1\le p\le 4\\
&\bold d:L\wedge I^+\longrightarrow F \quad \text{on}\quad \Cal N(\Cal D),\qquad \bold d_0=T,\quad \bold d_1=d\\
\endalign
$$
such that

\bigskip\noindent
(6.3.1) $T(E_0<E_1<E_2):T(E_0<E_2) \ \sim \ F(i_{12})T(E_0<E_1)+T(E_1<E_2)$, \  \ $i_{12}$ inclusion.

\bigskip\noindent
(6.3.2) $T(E_{0123})$ \  has boundary values 
$$
\CD
T(E_{03})@>T(E_{013}) > \sim > F(i_{13})T(E_{01})+T(E_{13})\\
@V \sim V T(E_{023}) V   @V \sim V 1+T(E_{123}) V\\
F(i_{23})T(E_{02})+T(E_{23})@>F(i_{23})T(E_{012})+1 > \sim > F(i_{13})T(E_{01})+
F(i_{23})T(E_{12})+T(E_{23})\\
\endCD 
$$

\bigskip\noindent
(6.3.3) $T(E_{01234})$ \ has boundary values
$$
T(E_{01234})(\_,t_1,t_2,t_3)=
\cases
T(E_{0234})(\_,t_2,t_3),\qquad &t_1=0\\
T(E_{0134})(\_,t_1,t_3),\qquad &t_2=0\\
T(E_{0124})(\_,t_1,t_2),\qquad &t_3=0\\
F(i_{14}) T(E_{01})+T(E_{1234})(\_,t_2,t_3),\qquad &t_1=1\\
F(i_{24}) T(E_{012})(\_,t_1)+T(E_{234})(\_,t_3),\qquad &t_2=1\\
F(i_{34}) T(E_{0123})(\_,t_1,t_2)+T(E_{34}),\qquad &t_3=1\\
\endcases
$$

\bigskip\noindent
(6.3.4) $T|\Cal N(\Cal T^2)=\varepsilon$, \ and \ $\bold d|\Cal N(\Cal T^2):\varepsilon \ \sim \ \varepsilon$
is the constant path.

\bigskip\noindent
(6.3.5)  $T(E_0<E_0)=*$, and
$$
\align
&T(E_{001}):T(E_{01})\ \sim \ *+T(E_{01})\\
&T(E_{011}):T(E_{01})\ \sim \ T(E_{01})+*\\
\endalign
$$
are standard.

\bigskip\noindent
(6.3.6) $T(E_{0023})$ reduces by a standard homotopy to
a map with edge path
$$
\CD
T(E_{03})@>T(E_{023}) > \sim > F(i_{23})T(E_{02})+T(E_{23})\\
@V \sim V 1 V   @V \sim V 1 V\\
T(E_{03})@>T(E_{023}) > \sim > F(i_{23})T(E_{02})+T(E_{23})\\
\endCD
$$
We assume $T(E_{0023})$ is homotopic to $\left((\_,t_1,t_2)\longrightarrow T(E_{023})(\_,t_1)\right)$,  by a homotopy which is standard on  
$\partial(I^2)$. A similar assumption is made for \ $T(E_{0113})$, and for $T(E_{0122})$. 

\bigskip
An excisive diagram
$$
\Cal E=
\matrix
E_1&< &E\\
\vee && \vee\\
E_0&< &E_2\\
\endmatrix,\qquad\qquad \Cal E=\Cal E\downarrow (B,A)
$$
corresponds to an excisive map \ $(i,1):(E_0<E_2)@>>>(E_1<E)$, \ where $i$ is inclusion. In this case we will write $T(\Cal E)$ in place of $T(i,1)$. 
We have
$$T(\Cal E):F(i_2)T(E_0<E_2) \ \sim T(E_1<E)$$
Define
$$\bold a(\Cal E):T(E_0<E) \ \sim \ F(i_{01})T(E_0<E_1)+F(i_{02})T(E_0<E_2)$$
by
$$T(E_0<E) \ \overset{(1)}\to {\sim} \ F(i_{01})T(E_0<E_1)+)T(E_1<E) \ \overset{(2)}\to {\sim} \ F(i_{01})T(E_0<E_1)+F(i_{02})T(E_0<E_2)$$
\medskip\noindent
where $(1)=T(E_0<E_1<E)$, \  $(2)=1+T(\Cal E)^{-1}$. The dual diagram
$$
\Cal E^*=
\matrix
E_2&< &E\\
\vee && \vee\\
E_0&< &E_1\\
\endmatrix
$$
is excisive, so we have
$$\bold a(\Cal E^*):T(E_0<E) \ \sim \ F(i_{02})T(E_0<E_2)+F(i_{01})T(E_0<E_1)$$
Let $\bold b(\Cal E)=\bold a(\Cal E^*)\cdot \bold a(\Cal E)^{-1}$.
$$\bold b(\Cal E):F(i_{01})T(E_0<E_1)+F(i_{02})T(E_0<E_2)\ \sim \ F(i_{02})T(E_0<E_2)+F(i_{01})T(E_0<E_1)$$
We assume

\bigskip\noindent
(6.3.7) $\bold b(\Cal E)$ \ is a permutation path.

\bigskip
The final property  concerns expansions. Suppose given \ $E_0<E_1$ \ and a pair $(\phi,\lambda)$ where $\phi:E_0@>>>E_1$ is a fiberwise 
diffeomorphism and 
$\lambda:E_0\times I@>>>E_1$ is a map such that $\lambda_0=i_{01}$ and $\lambda_1=\phi$. We will refer to such a pair as an {\it expansion}.
Define
$$\bold c(\phi,\lambda):F(i_{01})T(E_0) \ \sim \ T(E_1) $$
by
$$F(i_{01})T(E_0) \ \overset{(1)}\to\sim \ F(\phi)T(E_0) \ 
\overset{T(\phi)}\to\sim \ T(E_1)$$
where $(1)=\widehat F(\lambda)(T(E_0)\wedge 1)$. There is the associated path class
$$\overline {\bold c}(\phi,\lambda):T(E_0<E_1) \ \phom \ * $$
defined by 
$$T(E_0<E_1) \ \overset{(1)}\to\sim \  -F(i_{01})T(E_0)+T(E_1)  \ \overset{(2)}\to\sim \  -T(E_1)+T(E_1) \  \overset{\frak s}\to\sim \ *$$
where $(1)$ is a rearrangement of $T(\Phi<E_0<E_1)$, and $(2)=-\bold c(\phi,\lambda)+1$.
\medskip
Suppose that we have an excision
$$
\Cal E=
\matrix
E_0 & < & E_1\\
\vee & & \vee\\
\widetilde E_0 & < & \widetilde E_1\\
\endmatrix
$$
and expansion $(\phi,\lambda)$ for $E_0<E_1$, such that $\phi$ restricts to a diffeomorphism $\widetilde\phi:\widetilde E_0@>>>\widetilde E_1$,
and $\lambda$ restricts to $\widetilde \lambda:\widetilde E_0\times I@>>>\widetilde E_1$. Then $(\widetilde \phi,\widetilde \lambda)$ is an expansion
for $\widetilde E_0<\widetilde E_1$. We assume 

\bigskip\noindent
(6.3.8) The following commutes.
$$
\CD
T(E_0<E_1)@>\overline{\bold c}(\phi,\lambda)> \phom > \ \ *\\
@A \phom A T(\Cal E) A         @A \phom A 1 A\\
F(i)T(\widetilde E_0<\widetilde E_1)@>F(i)\overline {\bold c}(\widetilde \phi,\widetilde \lambda)> \phom > \ \ *\\
\endCD 
$$

\bigskip
With the exception of (6.3.8), which is only needed at the path class level, the other properties should be understood as statemants
about multiplicative simplicial transformations. For (6.3.7), $\bold b(\Cal E)$ is defined  using (3.7), as a multiplicative simplicial transformation
on a subcomplex of $\Cal N(\text{Ar}(\Cal B^1))$. 

\bigskip
\subhead 7. Homotopy property \endsubhead  Recall that for a homotopy $H:X\times I@>>>Y$,
and multiplicative functor $(F,\omega_F)$, the induced homotopy \newline
$\widehat H:F(X)\wedge I^+@>>>F(Y)$, is given by
$$F(X)\wedge I^+@>1\wedge\iota>>F(X)\wedge L(I)@>\omega_F>>F(X\times I)@>F(H)>>F(Y)$$
where $\iota:I^+@>>>L(I)$ is the canonical map. In particular, for the identity map $1_{X\times I}:X\times I@>>>X\times I$,
we have $\widehat F(1_{X\times I})=\omega_F(1\wedge\iota)$.

Suppose now that we are given $H:X\times I\times I@>>>Y$. For $(a,b)\in I\times I$, let 
$$
\align
&H^a:X\times I@>>>Y,\qquad H^a(x,t)=H(x,a,t)\\
&H_b:X\times I@>>>Y,\qquad H_b(x,s)=H(x,s,b)\\
&H_{(a,b)}:X@>>>Y,\qquad H_{(a,b)}(x)=H(x,a,b)\\
\endalign
$$ 

We have the  edge-path diagram.
$$
\CD
H_{(0,1)}@>H_1>\phom>H_{(1,1)}\\
@A H^0 A\phom A  @A H^1 A\phom A\\
H_{(0,0)}@>H_0>\phom>H_{(1,0)}\\
\endCD
$$

Let \ $(\tilde i_j,i_j):(E\times \{j\},B\times \{j\})@>>>(E\times I,B\times I)$, \ $j=0,1$, \ denote inclusions. 

\proclaim{(7.1) Lemma} Let $E=E\downarrow (B,A)$. The following commutes.
$$
\CD
T(E\times I)L(i_o) @> T(E\times I)\widehat L(1_{B\times I})>\phom> T(E\times I)L(i_1)\\
@A\phom A T(\widetilde i_0,i_0) A         @A\phom A T(\widetilde i_1,i_1) A\\
F(\widetilde i_0)T(E) @>\widehat F(1_{E\times I}))(T(E)\wedge 1)>\phom > F(\widetilde i_1) T(E)\\
\endCD
$$
\endproclaim
\demo{Proof}We have \ $T(E,I):L(B)\wedge L(I)\wedge I^+\longrightarrow F(E\times I)$. Let 
$$\overline T(E,I):L(B)\wedge I^+\wedge I^+\longrightarrow F(E\times I)$$
denote 
$$L(B)\wedge I^+\wedge I^+@>1\wedge\iota\wedge 1>>L(B)\wedge L(I)\wedge I^+@>T(E,I)>> F(E\times I)$$ 
Its edge path diagram has the form
$$
\CD
T(E\times I)L(i_o) @>T(E\times I)\widehat L(1_{B\times I})>\phom> T(E\times I)L(i_1)\\
@A\phom A \overline T(E,I)^0 A         @A\phom A \overline T(E,I)^1 A\\
F(\widetilde i_0)T(E) @>\widehat F(1_{E\times I})(T(E)\wedge 1)>\phom > F(\widetilde i_1)T(E).\\
\endCD
$$

To show that \ $\overline T(E,I)^0 \, \sim \, T(\tilde i_0,i_0)$ \ rel. $\dot I$, \ let 
$$\sigma=((E,\{0\})@>(1_E,l)>>(E,I)) \tag 7.2$$
where $l$ is inclusion. We have
$$T(\sigma):L(B)\wedge L(\{0\})\wedge I^+\wedge\Delta_1^+\longrightarrow F(E\times I)$$
Let $\overline T(\sigma)$ denote $T(\sigma)$ composed with \
$$L(B)\wedge I^+\wedge\Delta_1^+@>>>L(B)\wedge L(\{0\})\wedge I^+\wedge\Delta_1^+$$

We have the $2$-simplex
$$ \tau=(E@>(\widetilde a_,a_)>>E\times\{0\}@>(\widetilde b,b)>>E\times I)\tag 7.3$$
where the maps are  inclusions. Noting that $(\widetilde b\widetilde a,ba)=(\tilde i_0,i_0)$, the edge-path diagram of 
$\overline T(\sigma)$ is 
$$
\CD
F(\tilde i_o)T(E)@>\overline T(E,I)^0>\sim> T(E\times I)L(i_0)\\
@A \sim A 1 A      @A \sim A T(\widetilde b,b)L(a) A\\
F(\tilde i_o)T(E)@>F(\widetilde b)\overline T(E,\{0\})>\sim> F(\widetilde b)T(E\times \{0\})L(a)\\
\endCD \tag 7.4
$$
where $\overline T(E,\{0\})$ is $T(E,\{0\})$ composed with 
$$L(B)\wedge I^+@>>>L(B)\wedge L(\{0\})\wedge I^+$$ 
From the normality property $\bold I$, we have
$\overline T(E,\{0\})=T(\widetilde a,a)$. So
$$\overline T(E,I)^0 \ \sim \ F(\widetilde b)T(\widetilde a,a)\cdot T(\widetilde b,b)L(a) \qquad \text{rel. } \dot I$$ 
On the other hand, we have
$$T(\tau):T(\tilde i_0,i_0) \ \sim \ F(\widetilde b)T(\widetilde a,a) \cdot T(\widetilde b,b)L(a)
\qquad \text{rel. } \dot I$$
Therefore, \ $\overline T(E,I)^0 \sim T(\tilde i_0,i_0)$ \ rel. $\dot I$. The proof that \
$\overline T(E,I)^1 \sim T(\tilde i_1,i_1)$ \ rel. $\dot I$ \ is similar.
\enddemo

\bigskip
\proclaim{(7.5) Theorem} For $(\tilde H,H):(E\times I,B\times I,p\times 1)\rightarrow(E',B',p')$,
let $(\tilde h_j,h_j)=(\tilde H \tilde i_j, Hi_j)$, $j=0,1$. The following commutes.
$$
\CD
T(E')L(h_o) @>T(E')\widehat L(H)>\phom> T(E')L(h_1)\\
@A\phom A T(\tilde h_0,h_0) A         @A\phom A T(\tilde h_1,h_1) A\\
F(\tilde h_0) T(E) @>\widehat F(\tilde H)(T(E)\wedge 1)>\phom > F(\tilde h_1)T(E)\\
\endCD
$$
\endproclaim
\demo{Proof} In the diagram
$$
\CD
T(E')L(h_o) @>T(E')\widehat L(H)>\phom> T(E')L(h_1)\\
@A\phom A T(\tilde H,H)(L(i_0)\wedge 1) A         @A\phom A T(\tilde H,H)(L(i_1)\wedge 1) A\\
F(\tilde H)T(E\times I)L(i_o) @>F(\tilde H)T(E\times I)\widehat L(1_{B\times I})>\phom> 
F(\tilde H)T(E\times I)L(i_1)\\
@A\phom A F(\tilde H)T(\tilde i_0,i_0) A         @A\phom A F(\tilde H)T(\tilde i_1,i_1) A\\
F(\tilde h_0)T(E) @>\widehat F(\tilde H)(T(E)\wedge 1)>\phom > F(\tilde h_1)T(E),\\
\endCD
$$
the bottom square commutes by the lemma. The top square commutes as it is the edge path diagram of the map
$$L(B)\wedge I^+\wedge\Delta_1^+@>1\wedge\iota\wedge 1>>L(B)\wedge L(I)\wedge \Delta_1^+
@>\omega_L\wedge 1>> L(B\times I)\wedge \Delta_1^+@>T(\tilde H,H)>> F(E').$$
Finally,  \ $\sigma_j=(E@>(\tilde i_j,i_j)>>E\times I@>(\tilde H,H)>>E')$, \ defines
$$T(\tilde h_j,h_j) \, \sim \, F(\tilde H)T(\tilde i_j,i_j)\cdot T(\tilde H,H)(L(i_j)\wedge 1)
\qquad\text{rel. } \dot I$$
\enddemo

\bigskip
Suppose that $E=E\downarrow (B,A)$ is an orthogonal disk bundle, and $\Delta_0,\Delta_1:B@>>>E-\dot E$ are sections.
Let $\phi:B\times I@>>>E$ denote the linear path \ $b@>>>(1-t)\Delta_0(b)+t\Delta_1(B)$, \ and let $\Psi:E\times I@>>>E$
be a fiberwise diffeomorphism over the projection $\pi:B\times I@>>>B$, such that $\Psi_0=1_E$ and $\Psi\Delta_0=\phi$
Theorem (7.5) applied to $(\Psi,\pi)$ together with the naturality of $\bold d$, gives the following.

\bigskip
\proclaim{(7.6) Theorem} The following commutes.
\endproclaim
$$
\CD
T(E) @> 1>\phom> T(E)\\
@V\phom V \bold d(E,\Delta_0) V         @V\phom V \bold d(E,\Delta_1) V   \\
F(\Delta_0)\varepsilon(B,A) @>\widehat F(\phi)(\varepsilon(B,A)\wedge 1)>\phom > F(\Delta_1)\varepsilon(B,A)\\
\endCD
$$

\bigskip
\subhead 8. Vertical reduction \endsubhead Given a map $f:Y\rightarrow X$, where $X$ is pointed, define the {\it support} of 
$f$ to be 
$$\text{supp}(f)=\overline{Y-f^{-1}(*)}$$
Suppose now that we are given a bundle $E_2\downarrow (B,A)$, a map 
$$\gamma:E_2\downarrow (B,A)@>>>(X,*)$$ 
and $E_0<E_1<E_2$ such that 
$$\text{supp}(\gamma)\subset E_1-\text{fr}(E_1<E_2)$$ 
We may then define a {\it reduction path}
$$\frak r(E_{02}/E_{01}):F(\gamma)T((E_0<E_2)\downarrow (B,A)) \ \sim \  F(\gamma)T((E_0<E_1)\downarrow (B,A)) \tag 8.1$$
as follows. Choose $U_1<E_1$ such that \ $\text{supp}(\gamma)\subset E_1-U_1$, and $U_1$ is a neighborhood in $E_1$ of 
$\text{fr}(E_1<E_2)$. Let \ $U_2=U_1\cup (E-E_1)$. We have an excision
$$
\matrix
E_1&<&E_2\\
\vee & &\vee\\
U_1&<&U_2\\
\endmatrix 
$$
and $\gamma(U_2)=*$. Define $\frak r(E_{02}/E_{01})$ by applying $F(\gamma)$ to
$$
\align
T((E_0<E_2) &\sim \ F(i)T((E_0<E_1)+T((E_1<E_2)\\
                     &\ \sim \ F(i)T((E_0<E_1)+F(i)T((U_1<U_2)
\endalign                     
$$
The path class does not depend on the choice of $U_1$. Let $C_1\subset E_1$ be a collar about $\text{fr}(E_1<E_2)$
such that $C_1\subset U_1$, and let $C_2=C_1\cup (E-E_1)$. The assertion follows by consideration of
$$
\matrix
E_1&<&E_2\\
\vee & &\vee\\
U_1&<&U_2\\
\vee & &\vee\\
C_1&<&C_2\\
\endmatrix 
$$
and the fact that
$$F(\gamma)T\left(
\matrix
U_1&<&U_2\\
\vee & &\vee\\
C_1&<&C_2\\
\endmatrix \right)=*
$$

\bigskip
{\it Relation with additivity.} Suppose that we have \ $E_0<E_1<E_2<E_3$ \ where
$$\text{supp}(\gamma)\subset E_2-\text{fr}(E_2<E_3)$$
A map
$$T(E_{013}/E_{012}):L(B,A)\wedge (I^2)^+\longrightarrow F(X,*) \tag 8.2$$
with edge path
$$
\CD
F(\gamma)T(E_0<E_3) @> \frak r(E_{03}/E_{02}) > \sim> F(\gamma)T(E_0<E_2)\\
@V \sim V F(\gamma)T(E_0<E_1<E_3) V   @V\sim V F(\gamma)T(E_0<E_1<E_2) V \\
F(\gamma)T(E_0<E_1)+F(\gamma)T(E_1<E_3)  @>1+\frak r(E_{13}/E_{12}) > \sim > F(\gamma)T(E_0<E_1)+F(\gamma)T(E_1<E_2)\\
\endCD
$$
is defined as follows. In the diagram
{\eightpoint
$$
\CD
T(E_{03}) @> T(E_{023}) > \sim > F(i)T(E_{02})+T(E_{23}) @>  > \sim > F(i)T(E_{02})+F(i)T(U_{23})\\
@V \sim V T(E_{013}) V   @V\sim V T(E_{012})+1 V  @V\sim V T(E_{012})+1 V\\
F(i)T(E_{01})+T(E_{13}) @> 1+T(E_{123}) > \sim > F(i)T(E_{01})+F(i)T(E_{12})+T(E_{23})
@>  > \sim > F(i)T(E_{01})+F(i)T(E_{12})+T(U_{23})\\
\endCD 
$$}
the left square is filled in by the associativity path \ $T(E_{0123})$. Define 
(9.2) by applying $F(\gamma)$ to the diagram.

\bigskip
{\it Relation with excision.} Suppose we have an excisive diagram
$$\matrix
E_0 & < & E_1 & < & E_2\\
\vee & & \vee & & \vee\\
\widetilde E_0 & < & \widetilde E_1 & < & \widetilde E_2\\
\endmatrix
$$
Let
$$\Cal E_{ij}=\matrix
E_i & < & E_j\\
\vee & & \vee\\
\widetilde E_i & < & \widetilde E_j\\
\endmatrix
$$
We will define
$$\frak r(\Cal E_{02}/\Cal E_{01}):L(B,A)\wedge (I^2)^+\longrightarrow F(X) \tag 8.3$$
with edge path
$$
\CD
F(\gamma)T(\widetilde E_0<\widetilde E_2) @>\frak r(\widetilde E_{02}/\widetilde E_{01}) > \sim > 
F(\gamma)T(\widetilde E_0<\widetilde E_1)\\
@V \sim V F(\gamma)T(\Cal E_{02}) V   @V \sim V F(\gamma)T(\Cal E_{01}) V \\
F(\gamma)T(E_0<E_2) @>\frak r(E_{02}/E_{01}) > \sim > F(\gamma)T(E_0<E_1)\\
\endCD
$$
Since $\text{fr}(E_1<E_2)=\text{fr}(\widetilde E_1<\widetilde E_2)$, we may choose $U_1\subset \widetilde E_1$.
Then $U_2\subset \widetilde E_2$, and we have an excisive diagram
$$\matrix
E_1 & < & E_2\\
\vee & & \vee\\
\widetilde E_1 & < & \widetilde E_2\\
\vee & & \vee\\
U_1 & < & U_2\\
\endmatrix
$$
We define $\frak r(\Cal E_{02}/\Cal E_{01})$ by applying $F(\gamma)$ to the map with edge path
$$
\CD
F(i)T(\widetilde E_0<\widetilde E_2) @> T(\Cal E_{02}) > \sim >  T(E_0<E_2)\\
@V \sim V F(i)T(\widetilde E_{012}) V   @V \sim V T(E_{012}) V \\
F(i)T(\widetilde E_0<\widetilde E_1)+F(i)T(\widetilde E_1<\widetilde E_2) @> F(i)T(\Cal E_{01})+ T(\Cal E_{12})> 
\sim > F(i)T(E_0<E_1)+T(E_1<E_2)\\
@V \sim V 1+F(i)T(\widetilde {\Cal D}_{12})^{-1} V    @V \sim V 1+T(\Cal D_{12})^{-1} V \\
F(i)T(\widetilde E_0<\widetilde E_1)+F(i)T(U_1<U_2) @> T(\Cal E_{01})+ 1> \sim > F(i)T(E_0<E_1)+F(i)T(U_1<U_2)\\
\endCD
$$
\bigskip
\noindent
where
$$
\widetilde{\Cal D}_{12}=\matrix
\widetilde E_1 & < & \widetilde E_2\\
\vee & & \vee\\
U_1 & < & U_2\\
\endmatrix
\qquad\quad
\Cal D_{12}=\matrix
E_1 & < & E_2\\
\vee & & \vee\\
U_1 & < & U_2\\
\endmatrix
$$

\bigskip
Suppose given a fiberwise diffeomorphism
$$(\overline f,f):\widetilde E_0<\widetilde E_1<\widetilde E_2\longrightarrow E_0<E_1<E_2$$
Let $\widetilde \gamma=\gamma \overline f$,  and let $\overline f_{ij}:\widetilde E_i<\widetilde E_j\longrightarrow E_i<E_j$ denote the restriction of $(\overline f,f)$. Similar to the definition of (9.3), we define
$$\frak r(\overline f_{02}/ \overline f_{01}):L(B,A)\wedge (I^2)^+\longrightarrow F(X,*) \tag 8.4$$
with edge path
$$
\CD
F(\widetilde \gamma)T(\widetilde E_0<\widetilde E_2) @>\frak r(\widetilde E_{02}/\widetilde E_{01})> \sim > 
F(\widetilde \gamma)T(\widetilde E_0<\widetilde E_1)\\
@V \sim V F(\gamma)T(\overline f_{02},f) V   @V \sim V F(\gamma)T(\overline f_{01},f) V \\
F(\gamma)T(E_0<E_2)L(f) @>\frak r(E_{02}/E_{01})L(f)> \sim > F(\gamma)T(E_0<E_1)L(f)\\
\endCD
$$

\bigskip
{\it Iteration.} Given $E_0<E_1<E_2<E_3$ such that
$$\text{supp}(\gamma)\subset E_1-\text{fr}(E_1<E_2)$$
we shall define 
$$\frak r(E_{03}/E_{02}/E_{01}):L(B,A)\wedge(I^2)^+\longrightarrow F(X,*) \tag 8.5$$
with edge path
$$
\CD
F(\gamma)T(E_0<E_3)@>\frak r(E_{03}/E_{02})>\sim> F(\gamma)T(E_0<E_2)\\
@V\sim V \frak r(E_{03}/E_{01}) V     @V\sim V \frak r(E_{02}/E_{01}) V\\
F(\gamma)T(E_0<E_1)@> 1 >\sim> F(\gamma)T(E_0<E_1)\\
\endCD
$$

Consider
\bigskip
{\eightpoint
$$
\CD
T(E_0<E_3) @>>\sim> T(E_0<E_2)+T(E_2<E)@>>\sim> T(E_0<E_2)+T(U_2<U_3)\\
@V\sim VV         @V\sim VV @V\sim VV\\
T(E_0<E_1)+T(E_1<E_3) @>>\sim> T(E_0<E_1)+T(E_1<E_2)+T(E_2<E_3)@>>\sim>T(E_0<E_1)+T(E_1<E_2)+T(U_2<U_3)\\
@V\sim VV         @V\sim VV  @V\sim VV\\
T(E_0<E_1)+T(U_1<U_3) @>>\sim> T(E_0<E_1)+T(U_1<U_2)+T(U_2<U_3)@>>\sim>T(E_0<E_1)+T(U_1<U_2)+T(U_2<U_3)\\
\endCD
$$}
\bigskip
\noindent
The upper left square is filled in by $T(E_0<E_1<E_2<E_3)$, the lower left by
$$1+T\left(\matrix
E_1&<&E_2&<&E_3\\
\vee & &\vee & &\vee\\
U_1&<&U_2&<&U_3\\
\endmatrix \right)
$$
The right hand squares are filled in by the obvious maps. Now $\frak r(E_{03}/E_{02}/E_{01})$ is defined by applying 
$F(\gamma)$ to the filled in diagram.

\bigskip
The map $\frak r(E_{03}/E_{02}/E_{01})$ is related to additivity and excision in a similar way as $\frak r(E_{02}/E_{01})$
in (9.1). For additivity, suppose given $E_0<E_1<E_2<E_3<E_4$ with \ $\text{supp}(\gamma)\subset E_2-\text{fr}(E_2<E_3)$. \
We have the excision
$$
\Cal D_{234}=
\matrix
E_2&<&E_3&<&E_4\\
\vee & &\vee & &\vee\\
U_2&<&U_3&<&U_4\\
\endmatrix 
$$
Let
$$\Cal D_{ij}=\matrix
E_i&<&E_j\\
\vee & &\vee\\
U_i&<&U_j\\
\endmatrix 
$$
We have
$$L(B,A)\wedge(I^3)^+\longrightarrow F(E_4\downarrow (B,A)))$$
by
$$
(\_,t_1,t_2,t_3)\longrightarrow
\cases
T(E_{01234})(\_,2t_1,2t_2,t_3) \quad & 0\le t_1,t_2\le 1/2\\
T(E_{0123})(\_,2t_2,t_3)+T(\Cal D_{34})(\_,2-2t_1) \quad & 1/2\le t_1\le 1, \ 0\le t_2\le 1/2\\
T(E_{012})(\_,t_3)+T(\Cal D_{234})(\_,2t_1,2-2t_2) \quad & 0\le t_1\le 1/2, \ 1/2\le t_2 \le 1\\
T(E_{012}(\_,t_3)+T(\Cal D_{23})(\_,2-2t_2)+ T(\Cal D_{34})(\_,(2-2t_1)(2-2t_2))\quad & 1/2\le t_1\le 1, \ 1/2\le t_2 \le 1\\
\endcases
$$
Applying $F(\gamma)$ gives a map
$$T(E_{0124}/E_{0123}):L(B,A)\wedge(I^3)^+\longrightarrow F(X,*)  \tag 8.6$$
such that
$$
T(E_{0124}/E_{0123})(\_,t_1,t_2,t_3)=
\cases
\frak r(E_{04}/E_{03}/E_{02})(\_,t_1,t_2) \qquad & t_3=0\\
F(\gamma)T(E_{01})(\_) + \frak r(E_{14}/E_{13}/E_{12})(\_,t_1,t_2) \qquad & t_3=1\\
T(E_{014}/E_{012})(\_,t_2,t_3) \qquad & t_1=0\\
T(E_{013}/E_{012})(\_,t_2,t_3) \qquad & t_1=1\\
T(E_{014}/E_{013})(\_,t_1,t_3) \qquad & t_2=0\\
F(\gamma)T(E_{012})(\_,t_3) \qquad & t_2=1\\
\endcases
$$

\bigskip
For excision, suppose we have \ $E_0<E_1<E_2<E_3$ with
$$\text{supp}(\gamma)\subset E_1-\text{fr}(E_1<E_2)$$
and an excisive diagram
$$\Cal E=
\matrix
E_0 & < & E_1 & < & E_2 & < & E_3\\
\vee & & \vee & & \vee & & \vee\\
\widetilde E_0 & < & \widetilde E_1 & < & \widetilde E_2  & < & \widetilde E_3\\
\endmatrix
$$
Let
$$\Cal E_{ij}=\matrix
E_i & < & E_j\\
\vee & & \vee\\
\widetilde E_i & < & \widetilde E_j\\
\endmatrix
\qquad\qquad
\Cal E_{ijk}=\matrix
E_i & < & E_j & < & E_k\\
\vee & & \vee& & \vee\\
\widetilde E_i & < & \widetilde E_j& < & \widetilde E_k\\
\endmatrix
$$

$$
\Cal U=
\matrix
E_1 & < & E_2 & < & E_3\\
\vee & & \vee & & \vee\\
\widetilde E_1 & < & \widetilde E_2  & < & \widetilde E_3\\
\vee & & \vee & & \vee\\
U_1 & < & U_2 & < & U_3\\
\endmatrix
\qquad\qquad
\Cal U_{ij}=\matrix
E_i & < & E_j\\
\vee & & \vee\\
\widetilde E_i & < & \widetilde E_j\\
\vee & & \vee\\
U_i & < & U_j\\
\endmatrix 
$$

We have 
$$\align
T(\Cal U)&:L(B,A)\wedge I\wedge \Delta_2^+\longrightarrow F(E_3\downarrow (B,A))\\
T(\Cal U_{ij})&:L(B,A)\wedge \Delta_2^+\longrightarrow F(E_j\downarrow (B,A))\\
\endalign
$$

Relative to the standard simplicial decomposition of $I^2$, let 
$$
\align
\phi:I^2@>>> \Delta_2,\qquad  &\phi(0,0)=\phi(0,1)=e_0, \quad \phi(1,0)=e_1, \quad \phi(1,1)=e_2\\
\psi:I^3@>>> \Delta_2,\qquad  &\psi(t_1,t_2,t_3)=(1-t_2)\phi(t_1,t_3)+t_2e_0\\
\endalign
$$

Define
$$L(B,A)\wedge(I^3)^+\longrightarrow F(E_3\downarrow (B,A)))$$
by 
$(\_,t_1,t_2,t_3)\longrightarrow$
$$
\cases
T(\Cal E)(\_,2t_1,2t_2,t_3) \quad & 0\le t_1,t_2\le 1/2\\
T(\Cal E_{012})(\_,2t_2,t_3)+T(\Cal U_{23}))(\_,\phi(2-2t_1,t_3)) \quad & 1/2\le t_1\le 1, \ 0\le t_2\le 1/2\\
T(\Cal E_{01})(\_,t_3)+T(\Cal U)(\_,2t_1,\phi(2-2t_2,t_3)) \quad & 0\le t_1\le 1/2, \ 1/2\le t_2 \le 1\\
T(\Cal E_{01})(\_,t_3)+T(\Cal U_{12})(\_,\phi(2t_2-1,t_3))+ T(\Cal U_{23})(\_,\psi(2-2t_1,2t_2-1,t_3))
\quad & 1/2\le t_1\le 1, \ 1/2\le t_2 \le 1\\
\endcases
$$
Applying $F(\gamma)$ defines
$$\frak r(\Cal E_{03}/\Cal E_{02}/\Cal E_{01}):L(B,A)\wedge(I^3)^+\longrightarrow F(X,*) \tag 8.7$$
such that
$$
\frak r(\Cal E_{03}/\Cal E_{02}/\Cal E_{01})(\_,t_1,t_2,t_3)=
\cases
\frak r(\widetilde E_{03}/\widetilde E_{02}/\widetilde E_{01})(\_,t_1,t_2) \qquad & t_3=0\\
\frak r(E_{03}/E_{02}/E_{01})(\_,t_1,t_2) \qquad & t_3=1\\
\frak r(\Cal E_{03}/\Cal E_{01})(\_,t_2,t_3) \qquad & t_1=0\\
\frak r(\Cal E_{02}/\Cal E_{01})(\_,t_2,t_3) \qquad & t_1=1\\
\frak r(\Cal E_{03}/\Cal E_{02})(\_,t_1,t_3)\qquad & t_2=0\\
F(\gamma)T(\Cal E_{01})(\_,t_3) \qquad & t_2=1\\
\endcases
$$

\bigskip
Similarly, given a fiberwise diffeomorphism
$$(\overline f,f):\widetilde E_0<\widetilde E_1<\widetilde E_2<\widetilde E_3 \longrightarrow E_0<E_1<E_2<E_3$$
we define
$$\frak r(\overline f_{03}/\overline f_{02}/\overline f_{01}):L(B,A)\wedge(I^3)^+\longrightarrow F(X,*) \tag 8.8$$
such that
$$
\frak r(\overline f_{03}/\overline f_{02}/\overline f_{01})(\_,t_1,t_2,t_3)=
\cases
\frak r(\widetilde E_{03}/\widetilde E_{02}/\widetilde E_{01})(\_,t_1,t_2) \qquad & t_3=0\\
\frak r(E_{03}/E_{02}/E_{01})L(f)(\_,t_1,t_2) \qquad & t_3=1\\
\frak r(\overline f_{03})/\overline f_{01})(\_,t_2,t_3) \qquad & t_1=0\\
\frak r(\overline f_{02}/\overline f_{01})(\_,t_2,t_3) \qquad & t_1=1\\
\frak r(\overline f_{03}/\overline f_{02})(\_,t_1,t_3)\qquad & t_2=0\\
F(\gamma)T((\overline f_{01},f))(\_,t_3) \qquad & t_2=1\\
\endcases
$$

\bigskip
We shall need one more level of iteration. Suppose given $E_0<E_1<E_2<E_3<E_4$, such that
$$\text{supp}(\gamma)\subset E_1-\text{fr}(E_1<E_2)$$
The following is defined from (5.5), as (9.4) is defined from (5.4).
$$\frak r(E_{04}/E_{03}/E_{02}/E_{01}):L(B,A)\wedge(I^3)^+\longrightarrow F(X,*) \tag 8.9$$
such that 

$$
\frak r(E_{04}/E_{03}/E_{02}/E_{01})(\_,t_1,t_2,t_3)=
\cases
\frak r(E_{04}/E_{03}/E_{02})(\_,t_1,t_2) \qquad & t_3=0\\
F(\gamma)T(E_{01})(\_) \qquad & t_3=1\\
\frak r(E_{04}/E_{02}/E_{01})(\_,t_2,t_3) \qquad & t_1=0\\
\frak r(E_{03}/E_{02}/E_{01})(\_,t_2,t_3) \qquad & t_1=1\\
\frak r(E_{04}/E_{03}/E_{01})(\_,t_1,t_3) \qquad & t_2=0\\
\frak r(E_{02}/E_{01})(\_,t_3) \qquad & t_2=1\\
\endcases
$$

\bigskip
\subhead 9. Reduction of $\bold a(\Cal E)$\endsubhead For an excisive diagam of the form
$$\Cal E= \matrix
E_1 & < & E\\
\vee & & \vee\\
E_0 & < & E_2\\
\endmatrix 
$$
there is the path
$$\bold a(\Cal E):T(E_0<E) \ \sim \ F(i)T(E_0<E_1)+F(i)T(E_0<E_2) \tag 9.1$$
We shall discuss various reductions of $F(\gamma)\bold a(\Cal E)$. Let
$$\Cal E_1= E_1\cap\Cal E:=
\matrix
E_1 & < & E_1\\
\vee & & \vee\\
E_0 & < & E_0\\
\endmatrix \qquad
\Cal E_2= E_2\cap\Cal E:=
\matrix
E_0 & < & E_2\\
\vee & & \vee\\
E_0 & < & E_2\\
\endmatrix \qquad
\Cal E_0= E_0\cap\Cal E:=
\matrix
E_0 & < & E_0\\
\vee & & \vee\\
E_0 & < & E_0\\
\endmatrix 
$$

\bigskip\noindent
(i) Suppose that $\text{supp}(\gamma)\subset E_1-\text{fr}(E_1<E)$. Combining (8.2) and (8.3) defines
$$\bold a(\Cal E/\Cal E_1):L(B,A)\wedge (I^2)^+\longrightarrow F(X,*) \tag 9.2$$
with edge path
$$
\CD
F(\gamma)T(E_0<E) @>F(\gamma)\bold a(\Cal E) > \sim > F(\gamma)T(E_0<E_1)+F(\gamma)T(E_0<E_2)\\
@V \sim V \frak r(E/E_1) V   @V\sim V 1+\frak r(E_2/E_0) V \\
F(\gamma)T(E_0<E_1) @>F(\gamma)\bold a(\Cal E_1) > \sim > F(\gamma)T(E_0<E_1)+F(\gamma)T(E_0<E_0)\\
\endCD
$$

Note that $T(E_0<E_0)=*$, and $\bold a(\Cal E_1)$ is a standard path, so the edge path of $\bold a(\Cal E/\Cal E_1)$
simplifies to
$$
\CD
F(\gamma)T(E_0<E) @>F(\gamma)\bold a(\Cal E) > \sim > F(\gamma)T(E_0<E_1)+F(\gamma)T(E_0<E_2)\\
@V \sim V \frak r(E/E_1) V   @V\sim V 1+\frak r(E_2/E_0) V \\
F(\gamma)T(E_0<E_1) @> \frak s > \sim > F(\gamma)T(E_0<E_1)+*\\
\endCD
$$

\bigskip\noindent
(ii) Suppose that $\text{supp}(\gamma)\subset E_2-\text{fr}(E_2<E)$. We shall define
$$\bold a(\Cal E/\Cal E_2):L(B,A)\wedge (I^2)^+\longrightarrow F(X,*) \tag 9.3$$
with edge path
$$
\CD
F(\gamma)T(E_0<E) @>F(\gamma)\bold a(\Cal E) > \sim > F(\gamma)T(E_0<E_1)+F(\gamma)T(E_0<E_2)\\
@V \sim V \frak r(E/E_2) V   @V\sim V \frak r(E_1/E_0)+1 V \\
F(\gamma)T(E_0<E_2) @> F(\gamma)\bold a(\Cal E_2)=\frak s > \sim > *+F(\gamma)T(E_0<E_2)\\
\endCD
$$

There is the dual diagram
$$\Cal E^*=\matrix
E_2&<&E\\
\vee & &\vee\\
E_0&<&E_1\\
\endmatrix 
$$
and $(\Cal E_2)^*=(\Cal E^*)_1$. From (i), we have $\bold a(\Cal E^*/(\Cal E^*)_1)$ with edge path
$$
\CD
F(\gamma)T(E_0<E) @>F(\gamma)\bold a(\Cal E^*) > \sim > F(\gamma)T(E_0<E_2)+F(\gamma)T(E_0<E_1)\\
@V \sim V \frak r(E/E_2) V   @V\sim V 1+\frak r(E_1/E_0) V \\
F(\gamma)T(E_0<E_2) @> \frak s > \sim > F(\gamma)T(E_0<E_2)+*\\
\endCD
$$

By normality, we have that \ $\bold a(\Cal E^*)^{-1}\cdot\bold a(\Cal E)$ \ is a permutation path. Consequently, we have
$$\bold c(\Cal E):L(B,A)\wedge (I^2)^+\longrightarrow F(X,*) \tag 9.4$$ 
with edge path
$$
\CD
F(\gamma)T(E_0<E) @>F(\gamma)\bold a(\Cal E^*) > \sim > F(\gamma)T(E_0<E_2)+F(\gamma)T(E_0<E_1)\\
@V \sim V 1 V   @V\sim V \rho V \\
F(\gamma)T(E_0<E) @> F(\gamma)\bold a(\Cal E) > \sim > F(\gamma)T(E_0<E_1)+F(\gamma)T(E_0<E_2)\\
\endCD
$$
where $\rho$ denotes a $(2,1)$ permutation path. Let 
$$Q(\Cal E/\Cal E_2):L(B,A)\wedge (I^3)^+\longrightarrow F(X,*)$$
be such that
$$
Q(\Cal E/\Cal E_2)(\_,t_1,t_2,t_3)=
\cases
\frak r(E/E_2)(\_,t_2)\qquad &t_1=0\\
\rho(F(\gamma)T(E_0<E_2),\frak r(E_1/E_0))(\_,t_2,t_3)\qquad &t_1=1\\
\bold c(\Cal E)(\_,t_1,t_3)\qquad &t_2=0\\
\bold c(\Cal E_2)(\_,t_1,t_3)\qquad &t_2=1\\
\bold a(\Cal E^*/(\Cal E^*)_1)(\_,t_1,t_2)\qquad &t_3=0\\
\endcases
$$
Now define $\bold a(\Cal E/\Cal E_2)(\_,t_1,t_2)=Q(\Cal E/\Cal E_2)(\_,t_1,t_2,1)$. 

\bigskip\noindent
(iii) Suppose that $\text{supp}(\gamma)\subset E_0-\text{fr}(E_0<E)$. Let
$$\Cal V=\matrix V_1&<&V\\
                 \vee&&\vee\\
                 V_0&<&V_2\\
\endmatrix \qquad\qquad
\Cal D=\matrix \Cal E\\
                   \vee\\
                   \Cal V\\
\endmatrix
$$
Define
$$\bold a(\Cal E/\Cal E_0):L(B,A)\wedge (I^2)^+\longrightarrow F(X,*),\qquad
\bold a(\Cal E/\Cal E_0)(\_,t_1,t_2)=F(\gamma)\bold a (\Cal D)(\_,t_1,1-t_2) \tag 9.5$$
It has edge path
$$
\CD
F(\gamma)T(E_0<E) @>F(\gamma)\bold a(\Cal E) > \sim > F(\gamma)T(E_0<E_1)+F(\gamma)T(E_0<E_2)\\
@V \sim V \frak r(E/E_0) V   @V\sim V \frak r(E_1/E_0)+\frak r(E_2/E_0) V \\
* @> \frak s > \sim > *+ *\\
\endCD
$$
                                                                                                              
Let 
$$\bold a(\Cal E/\Cal E_1/\Cal E_0):L(B,A)\wedge(I^3)^+\longrightarrow F(X,*) \tag 9.6$$
by
$$\bold a(\Cal E/\Cal E_1/\Cal E_0)(\_,t_1,t_2,t_3)=F(\gamma)\bold a(\Cal D/\Cal D_1)(\_,t_1,t_2,1-t_3)$$
where $\bold a(\Cal D/\Cal D_1)$ is defined in a similar manner as $\bold a(\Cal E/\Cal E_1)$. We have
$$
\bold a(\Cal E/\Cal E_1/\Cal E_0)(\_,t_1,t_2,t_3)=
\cases
\bold a(\Cal E/\Cal E_1)(\_,t_1,t_2) \qquad & t_3=0\\
* \qquad & t_3=1\\
\frak r(E/E_1/E_0)(\_,t_2,t_3) \qquad & t_1=0\\
\frak r(E_1/E_1/E_0)(\_,t_2,t_3)+ \frak r(E_2/E_0/E_0)(\_,t_2,t_3) \qquad & t_1=1\\
\bold a(\Cal E/\Cal E_0)(\_,t_1,t_3) \qquad & t_2=0\\
\bold a(\Cal E_1/\Cal E_0)(\_,t_1,t_3) \qquad & t_2=1\\
\endcases
$$

\bigskip
We shall define
$$\bold a(\Cal E/\Cal E_2/\Cal E_0):L(B,A)\wedge(I^3)^+\longrightarrow F(X,*) \tag 9.7$$
such that
$$
\bold a(\Cal E/\Cal E_2/\Cal E_0)(\_,t_1,t_2,t_3)=
\cases
\bold a(\Cal E/\Cal E_2)(\_,t_1,t_2) \qquad & t_3=0\\
* \qquad & t_3=1\\
\frak r(E/E_2/E_0)(\_,t_2,t_3) \qquad & t_1=0\\
\frak r(E_1/E_0/E_0)(\_,t_2,t_3)+ \frak r(E_2/E_2/E_0)(\_,t_2,t_3) \qquad & t_1=1\\
\bold a(\Cal E/\Cal E_0)(\_,t_1,t_3) \qquad & t_2=0\\
\bold a(\Cal E_2/\Cal E_0)(\_,t_3) \qquad & t_2=1\\
\endcases
$$
By normality, we have
$$\bold c(\Cal D):L(B,A)\wedge (I^3)^+\longrightarrow F(X,*)$$ 
with edge path
{\eightpoint
$$
\CD
F(\gamma)T\left(\matrix E_0&<&E\\
                        \vee&&\vee\\
                        V_0&<&V\\
                 \endmatrix\right)
@>F(\gamma)\bold a(\Cal D^*) > \sim > 
F(\gamma)T\left(\matrix E_0&<&E_2\\
                        \vee&&\vee\\
                        V_0&<&V_2\\
                 \endmatrix\right)
+F(\gamma)T\left(\matrix E_0&<&E_1\\
                        \vee&&\vee\\
                        V_0&<&V_1\\
                 \endmatrix\right)\\
@V \sim V 1 V   @V\sim V \rho V \\
F(\gamma)T\left(\matrix E_0&<&E\\
                        \vee&&\vee\\
                        V_0&<&V\\
                 \endmatrix\right)
@>F(\gamma)\bold a(\Cal D) > \sim > 
F(\gamma)T\left(\matrix E_0&<&E_1\\
                        \vee&&\vee\\
                        V_0&<&V_1\\
                 \endmatrix\right)
+
F(\gamma)T\left(\matrix E_0&<&E_2\\
                        \vee&&\vee\\
                        V_0&<&V_2\\
                 \endmatrix\right)\\
\endCD
$$}
where $\rho$ is a $(2,1)$ permutation path. Let
$$Q(\Cal E/\Cal E_2/\Cal E_0):L(B,A)\wedge (I^4)^+\longrightarrow F(X,*)$$
be such that
$$
Q(\Cal E/\Cal E_2/\Cal E_0)(\_,t_1,t_2,t_3,t_4)=
\cases
\frak r(E/E_2/E_0)(\_,t_2,t_3)\qquad &t_1=0\\
\rho(\frak r(E_2/E_2/E_0),\frak r(E_1/E_0/E_0))(\_,t_2,t_3,t_4)\qquad &t_1=1\\
\bold c(\Cal D)(\_,t_1,1-t_3,t_4)\qquad &t_2=0\\
\bold c(\Cal D_2)(\_,t_1,1-t_3,t_4)\qquad &t_2=1\\
Q(\Cal E/\Cal E_2)(\_,t_1,t_2,t_4)\qquad &t_3=0\\
*\qquad &t_3=1\\
\bold a(\Cal E^*/(\Cal E^*)_1/(\Cal E^*)_0)(\_,t_1,t_2,t_3)\qquad &t_4=0\\
\endcases
$$
Define 
$$\bold a(\Cal E/\Cal E_2/\Cal E_0)(\_,t_1,t_2,t_3)=Q(\Cal E/\Cal E_2/\Cal E_0)(\_,t_1,t_2,t_3,1))$$

\bigskip
\subhead 10. Sections \endsubhead The proof of the main theorem is a matter of refining the axiomatic characterization of transfer given in [Becker-Schultz]. The main construction is discussed here. 
Given bundles $E_0$, $E$ over $B$ and a fiberwise smooth embedding $f:E_0@>>>E$ over the identity, suppose that
$U$ is an orthogonal vector bundle over $E_0$ and $\phi:U@>>>E$ is a codimension $0$ fiberwise smooth embedding such that $\phi(0_e)=f(e)$. 
We shall refer to $\phi$ as a {\it regular extension} of $f$.

Let $p:E\rightarrow B$ be a bundle. Let
$\dot p:\dot E\rightarrow B$ denote the fiberwise boundary of $E$ and let
$\theta:\dot E\times [-1,1]\rightarrow E$  be a collar with 
$\theta(\dot e,1)=\dot e$. Set
$$ E(s,t)=\theta(\dot E\times [s,t]), \qquad
E(s)=E-\theta(\dot E\times(s,1])\qquad 0\le s\le t\le 1$$

\medskip
Let $A$ be a closed subspace of $B$. We suppose given a pointed space $X$, a positive number $\epsilon^*$, 
and a family of maps
$$\gamma_{\epsilon}:E\downarrow (B,A)\longrightarrow (X,*),\qquad 0<\epsilon\le \epsilon^*$$ 
such that for $\epsilon<\epsilon'$, \ $\text{supp}(\gamma_{\epsilon})\subset \text{supp}(\gamma_{\epsilon'})$.

\medskip
Let \ $\Delta:B\longrightarrow E$ \ be a section such that $\Delta(B)\subset E-\dot E$. 
We suppose given a decomposition $B=B_1\cup B_2$, of $B$ into closed subspaces such that $\{B_1,B_2\}$ 
has the following properties.
\medskip
\noindent
{\bf I}. For every neighborhood $N$ of $\Delta(B_1)$ in $E\downarrow B_1$ there is $\epsilon$ such that
$$\text{supp}(\gamma_{\epsilon}|E\downarrow B_1)\subset N$$

\medskip
\noindent
{\bf II}. $\Delta(B_2)$ lies in the fiberwise interior of $E(0,1)\downarrow B_2$. Defining
$$ \delta:B_2\longrightarrow \dot E\downarrow B_2 \qquad \overline\delta:B_2\longrightarrow I$$
by $\Delta(b)=\theta(\delta(b),\overline\delta(b)$, for every neighborhood $N$ of $\delta(B_2)$ in $\dot E\downarrow B_2$ 
there is $\epsilon$ such that
$$\text{supp}(\gamma_{\epsilon}|E\downarrow B_2)\subset \theta(N\times (0,1])$$
\medskip
\noindent
Let $B_0=B_1\cap B_2$. We will also assume that $(B,B_1\cup A)$ and $(B_2,B_0\cup (B_2\cap A))$ have the homotopy type of a CW-pair.
Under these assumptions we will define for sufficiently small $\epsilon$, a path
$$\frak h(E\downarrow (B,A)):F(\gamma)T(E\downarrow (B,A))\sim F(\gamma\Delta)\varepsilon(B,A)\qquad 
\gamma=\gamma_{\epsilon} \tag 10.1$$

To simplify notation we will consider the absolute case. The relative case requires only minor changes.
Let $V@>>>B$ be an orthogonal vector bundle and $\psi:V\longrightarrow E$ a
regular extension of $\Delta$. Let
$$P=\psi(D(V)),\qquad  Q=E-\psi(D^{1/2}(V)),\qquad  R=P\cap Q$$
There is the excisive diagram
$$\Cal A=
\matrix
P&<&E\\
\vee & &\vee\\
R&<&Q\\
\endmatrix 
$$
We will define  
$$\frak i:F(\gamma)T(E))\sim F(\gamma)T(P) \tag 10.2$$
Then $\frak h$ is defined by

$$F(\gamma)T(E) \ \overset{\frak i}\to\sim \ F(\gamma)T(P)\overset{(1)}\to\sim \ F(\gamma\Delta)\varepsilon,\qquad\qquad 
(1)=F(\gamma)\bold d((P,\Delta))^{-1}$$ 
 
By additivity and excision, $\frak i$ is determined by a path 
$$\overline{\frak j}:F(\gamma)T(R<Q) \ \sim \ * \tag 10.3$$

For the definition of $\frak j$, by condition {\bf I}, we have $\gamma(Q\downarrow B_1)=*$, for sufficiently small $\epsilon$. Consider the diagram
$$\CD
(X,*)@> 1 >>   (X,*)@> 1 >>  (X,*)\\
@AA \gamma A       @AA \gamma A @AA \gamma A\\
Q@>\tilde a>> (Q,Q\downarrow B_1) @<\tilde b<< (Q\downarrow B_2,Q\downarrow B_0)\\
@VVV               @VVV                @VVV\\
B @>a>> (B,B_1) @<b<< (B_2,B_0),\\
\endCD$$ 
where $a$, $b$ are inclusions. A path 
$$\frak l:F(\gamma)T(R<Q)\downarrow(B_2,B_0))\sim * \tag 10.4$$
defines $\frak j$ as follows: Let
$$\frak l_0:F(\gamma)T((R<Q)\downarrow (B,B_1))L(b)\ \sim \ * \tag 10.5$$
denote
$$F(\gamma)T((R<Q)\downarrow (B,B_1))L(b) \ \overset{\text{(1)}}\to\sim \ F(\gamma)T((R<Q)\downarrow (B_2,B_0)) \ 
\overset{\frak l}\to\sim \ *$$
where (1)=$F(\gamma)\bold T(\tilde b,b)^{-1}$. Since $(B,B_1)$ and $(B_2,B_0)$ have the homotopy type
of a CW-pair, $L(b)$ is a homotopy equivalence. Thus, $\frak l_0$ determines a path 
$$\frak l_1:F(\gamma)T((R<Q)\downarrow (B,B_1)) \, \sim \, * \tag 10.6$$
Then $\frak j$ is defined to be
$$F(\gamma)T(R<Q) \ \overset{(1)}\to\sim \ F(\gamma)T((R<Q)\downarrow (B,B_1))L(a)
 \ \overset{(2)}\to\sim \ *$$
where $(1)=F(\gamma)T(\tilde a,a)$, and $(2)=\frak l_1L(a)$.

\medskip
We turn now to the definition of $\frak l$ in (10.4). From the excision $\Cal A$, \ $\frak l$ is determined by
$$\frak g:F(\gamma)T(P<E)\downarrow(B_2,B_0))\sim * \tag 10.7$$
In what follows, all bundles are over $(B_2,B_0)$. Let $U@>>>B_2$ be an orthogonal vector bundle and \newline
$\widehat\psi:U\longrightarrow \dot E\downarrow B_2$ a regular extension of \ $\delta$. Let 
$$\widehat\theta=\theta(\widehat\psi\times 1):D(U)\times [-1,1]\longrightarrow E\downarrow B_2$$
By {\bf II}, we may take $P$ and $\gamma$ so that $P$ and  $\text{supp}(\gamma)$ are in the fiberwise interior of
$\widehat\theta(D(U)\times I)$. 

\medskip
Let $H(r)$ denote the half disk
$$H(r)=\{(u,t)\in U\times R \ | \ |(u,t-1)|\le r, \ t\le 1 \ \}$$
and let $\widehat P$ denote the image of  $H(2)$. We have
$$P<\widehat P<E,\qquad\qquad \text{supp}(\gamma)\subset \widehat P-\text{fr}(\widehat P<E)$$

There is the reduction path
$$\frak r(E/\widehat P):F(\gamma)T(E) \ \sim \ F(\gamma)T(\widehat P) \tag 10.8$$

For $P$ and $\widehat P$ we have
$$\bold d(P,\Delta):F(\Delta)\varepsilon(B_2,B_0) \ \sim \ T(P), \qquad  
\bold d(\widehat P,\Delta):F(\Delta)\varepsilon(B_2,B_0) \ \sim \ T(\widehat P)$$
Let 
$$\frak q:F(\gamma)T(\widehat P) \ \sim \ F(\gamma)T(P),\qquad \qquad \frak q=F(\gamma)\bold d(P,\Delta)\cdot  
F(\gamma)\bold d(\widehat P,\Delta)^{-1}\tag 10.9$$
Let
$$\overline{\frak g}:F(\gamma)T(E) \ \sim \ F(\gamma)T(P) \tag 10.10$$
be (7.8) followed by (7.9). Finally, let $\frak g$ in (10.7) be the path associated to $\overline{\frak g}$.

\bigskip
The path class of $\frak h$ is natural with respect to fiberwise diffeomorphisms in the following sense. Given
\ $(\widetilde f,f):E'\downarrow (B',A')\longrightarrow E\downarrow (B,A)$, let $\gamma'=\gamma\widetilde f$.
With respect to the data pulled back from $E$, we have $\frak h(E'\downarrow (B',A'))$. It is straigtforward to show 
that the following commutes.
$$
\CD
F(\gamma')T(E'\downarrow (B',A'))@>F(\gamma)T(\widetilde f,f)>\phom> F(\gamma)T(E\downarrow (B,A))L(f)\\ 
@V \phom V \frak h(E'\downarrow (B',A')) V     @V \phom V \frak h(E\downarrow (B,A))L(f) V\\
F(\gamma'\Delta')\varepsilon(B',A') @> 1> \phom > F(\gamma'\Delta')\varepsilon(B',A')\\
\endCD \tag 10.11
$$

We note that the path class of $\frak h$ does not depend on the choice of regular extension of $\Delta$. Suppose that
$\psi:V@>>>E$, \  $\psi':V'@>>>E$ are two such regular extensions. There is $\overline\psi:V\times I@>>>E\times I$,
a regular extension of $\Delta\times 1$, such that $\overline\psi_0=\psi$, and $\overline\psi_1=\psi'\eta$, where
$\eta:V@>>>V'$ is orthogonal. Let 
$$\overline\gamma:(E\times I)\downarrow(B\times I,A\times I)\times I@>>>(X,*)$$
denote \ $\text{proj.}:E\times I@>>>E$ followed by $\gamma$. Relative to the decomposition $\{B_1\times I,B_2\times I\}$,
and regular extensions $\overline\psi$ of $\Delta\times 1$, and $\phi\times 1$ of $\delta\times 1$, we have
$\frak h_{\overline\psi}((E\times I)\downarrow ((B\times I,A\times I))$. We have a commutative diagram
$$
\CD
F(\gamma)T(E)@>>\phom> F(\overline\gamma)T(E\times I)L(i_0)@>(1)>\phom> F(\overline\gamma)T(E\times I)L(i_1)@>>\phom>F(\gamma)T(E)\\ 
@V \phom V \frak h_{\psi}(E) V     @V \phom V \frak h_{\overline\psi}(E\times I)L(i_0)V 
@V \phom V \frak h_{\overline\psi}(E\times I)L(i_0)V @V \phom V \frak h_{\psi'\eta}(E)V\\
F(\gamma\Delta)\varepsilon(B,A) @> 1> \phom > F(\gamma\Delta)\varepsilon(B,A) @> 1> \phom >F(\gamma\Delta)\varepsilon(B,A) 
@>1>\phom> F(\gamma\Delta)\varepsilon(B,A)\\
\endCD 
$$
Where $(1)=F(\gamma)T(E\times I)\widehat L(1_B\times 1)$. By (7.5), together with the observation that
$\frak h_{\psi'\eta}(E)=\frak h_{\psi'}(E)$, the diagram compresses to
$$
\CD
F(\gamma)T(E)@>1>\phom>F(\gamma)T(E)\\ 
@V \phom V \frak h_{\psi}(E) V     @V \phom V \frak h_{\psi'}(E) V\\
F(\gamma\Delta)\varepsilon(B,A) @>1>\phom > F(\gamma\Delta)\varepsilon(B,A)\\
\endCD 
$$
By a similar argument, the path class of $\frak h$ does not depend on the choice of regular extension of $\delta$.

\bigskip
\subhead 11. Expansions\endsubhead Given \ $E_0<E_1<E_2$, and expansion maps 
$$(\phi_{01},\lambda_{01}):E_0@>>>E_1, \qquad (\phi_{12},\lambda_{12}):E_1@>>>E_2$$
we have an expansion map 
$$(\phi_{02},\lambda_{02}):E_0@>>>E_2$$ 
$\phi_{02}=\phi_{12}\phi_{01}$, \  $(\lambda_{02})_t=(\lambda_{12})_t(\lambda_{01})_t$, \ $0\le t\le 1$.

\proclaim{(11.1) Theorem} \ $\bold c(\phi_{02},\lambda_{02})  \ \sim \ \bold c(\phi_{12},\lambda_{12})\cdot F(i_{12})\bold c(\phi_{01},\lambda_{01})$
rel. $\dot I$.
\endproclaim

\demo{Proof} In the diagram
$$\CD
F(i_{02})T(E_0) @> F(i_{12})\bold c(\phi_{01},\lambda_{01})>\sim> F(i_{12})T(E_1)\\
@V \sim V \widetilde F(\lambda_{02}) V     @V \sim V (1) V\\
F(\phi_{02})T(E_0)  @> F(\phi_{12})T(\phi_{01})>\sim>  F(\phi_{12})T(E_1)\\
@V \sim V T(\phi_{02}) V     @V \sim V T(\phi_{12}) V\\
T(E_2)  @> 1 >\sim> T(E_2)\\
\endCD
$$
where $(1)=\widehat F(\lambda)(T(E_1)\wedge 1)$, the top square commutes. The map
$$T(E_0@>\phi_{01}>>E_1@>\phi_{12}>>E_2):L(B,A)\wedge\Delta_2^+\longrightarrow F(E_2,E_2\downarrow A)$$
provides a homotopy validating homotopy commutativity of the bottom square.
\enddemo

\bigskip
Suppose now that for $E_0<E_1$, bothe $E_0$ and $E_1$ are fiberwise diffeomorphic to an orthogonal disk bundle. Let $\Delta_j:B@>>>E_j-\dot E_j$
be sections such that $i_{01}$ is section preserving. For an expansion $(\phi,\lambda)$, we will write
$$(\phi,\lambda):(E_0,\Delta_0)@>>>(E_1,\Delta_1)$$
if $\phi$ and $\lambda_t$, $0\le t\le 1$, are section preserving. For $(E_j,\Delta_j)$ we have
$$\bold d(E_j,\Delta_j):T(E_j) \ \sim \ F(\Delta_j)\varepsilon(B,A)$$

\proclaim{(11.2) Theorem} \ $\bold c(\phi,\lambda)  \ \sim \ \bold d((E_1,\Delta_1)\cdot F(i_{01})\bold d(E_0,\Delta_0)^{-1}$
rel. $\dot I$.
\endproclaim

\demo{Proof} In the diagram
$$\CD
F(i_{01})T(E_0) @> (1)>\sim> F(i_{12})T(E_1)  @>T(\phi)>\sim> T(E_1) \\
@A \sim A F(i_{01})\bold d(E_0,\Delta_0) A     @A \sim A  F(\phi)\bold d(E_0,\Delta_0) A   @A \sim A \bold d(E_1,\Delta_1) A\\
F(\Delta_1)\varepsilon(B,A)  @> 1 >\sim>  F(\Delta_1)\varepsilon(B,A)  @> 1 >\sim> F(\Delta_1)\varepsilon(B,A)\\
\endCD
$$
where $(1)=\widehat F(\lambda)(T(E_0)\wedge 1)$, the left square is filled in by $\widehat F(\lambda)(\bold d(E_0,\Delta_0)\wedge 1)$. 
The right square is filled in by
$$\bold d((E_0,\Delta_0)@>\phi>>(E_1,\Delta_1)):L(B,A)\wedge I^+\wedge \Delta_1^+\longrightarrow F(E_1,E_1\downarrow A)$$
\enddemo

\bigskip
\subhead 12. Disk bundles\endsubhead Suppose that $E$ in section 7 is diffeomorphic to an orthogonal disk bundle.

\proclaim{(12.1) Theorem}  
$$\frak h(E\downarrow (B,A))=F(\gamma)\bold d(E,\Delta)^{-1}:F(\gamma)T(E\downarrow (B,A) \ \phom \ F(\gamma\Delta)\varepsilon(B,A)$$
\endproclaim
\demo{Proof} Consider bundles over $(B_2,B_0)$. We have $\widehat P=\lambda(H(2))$. Let
$$\widehat R=\widehat P-\text {Int}(\lambda(H(\sqrt 2)) \qquad\qquad \widehat Q=E-\text {Int}(\lambda(H(\sqrt 2))$$
We have an excision
$$\widehat{\Cal A}=
\matrix
\widehat P&<&E\\
\vee & &\vee\\
\widehat R&<& \widehat Q\\
\endmatrix 
$$
We can construct an extension \ $(\phi^0, \lambda^0):(\widehat P,\widehat R,\Delta)@>>>(E,\widehat Q,\Delta)$.
The path \ $\overline{\frak r}(E/\widehat P)$ \ corresponding to \ $\frak r(E/\widehat P)$ \ in (7.8) is 
$$\overline{\frak r}(E/\widehat P)=F(\gamma)T(\Cal A)^{-1} \tag 12.2$$
On the other hand, by property ( / / )-(vii),
$$
\CD
T(\widehat P<E)  @> \ \overline{\bold c}(\phi^0,\lambda^0) \ >\phom>  *\\
@A\phom A T(\Cal {A}) A  @A\phom A 1 A\\
T(\widehat R< \widehat Q)  @> \ \overline{\bold c}(\phi^0,\lambda^0) \ >\phom>  *\\
\endCD
$$
commutes. Since $\gamma(\widehat Q)=*$, we have
$$F(\gamma)\overline{\bold c}(\phi^0,\lambda^0)=F(\gamma)T(\Cal A)^{-1} \tag 12.3$$
Therefore
$$\frak r(E/\widehat P)=F(\gamma)\bold c(\phi^0,\lambda^0)^{-1} \tag 12.4$$

\bigskip\noindent
We can construct an expansion $(\phi^1,\lambda^1):(P,\Delta)@>>>(\widehat P,\Delta)$. From (11.2),
$$\bold d(\widehat P,\Delta)\cdot F(i)\bold d(P,\Delta)^{-1}=\bold c(\phi^1,\lambda^1) \tag 12.5$$
Let \  $(\phi^2,\lambda^2):(P,\Delta)@>>>(E,\Delta)$ \ denote $(\phi^1,\lambda^1)$ followed by $(\phi^0,\lambda^0)$.
From (11.1) together with (12.4) and (12.5), we have that $\frak g$ in (10.7) is
$$\frak g=F(\gamma )\overline{\bold c}(\phi^2,\lambda^2):F(\gamma)T(P<E)\downarrow(B_2,B_0))\sim *  \tag 12.6$$

\bigskip
Let $(\phi,\lambda):(P,R,\Delta)@>>>(E,Q,\Delta)$ be an expansion over $B$. By (11.2),
$$\overline{\bold c}(\phi,\lambda)=\overline{\bold c}(\phi^2,\lambda^2):T(P<E)\downarrow(B_2,B_0))\sim *  \tag 12.7$$
So we have
$$\frak g=F(\gamma )\overline{\bold c}(\phi,\lambda):F(\gamma)T(P<E)\downarrow(B_2,B_0))\sim *  \tag 12.8$$
It follows that $\frak i$ in (10.2) is
$$\frak i=F(\gamma)\bold c(\phi,\lambda)^{-1}=F(\gamma)\bold d(P,\Delta)\cdot F(\gamma)\bold d(E,\Delta)^{-1} \tag 12.9$$
Therefore \ $\frak h=F(\gamma)\bold d(E,\Delta)^{-1}$. This completes the proof.
\enddemo

\bigskip
\subhead 13. Reduction of $\frak h$\endsubhead Suppose we have $E_1<E$ and
$$\text{supp}(\gamma)\subset E_1-\text{fr}(E_1<E),\qquad \Delta(B)\subset E_1-\text{fr}(E_1<E)$$
Let $U_1< E_1$ be a neighborhood of fr$(E_1<E)$ such that supp$(\gamma)\subset E_1-U_1$, and
$\Delta(B)\subset E_1-U_1$. Let $\theta_1:\dot E_1\times I@>>>E_1$ be a collar
such that $\theta_1=\theta$ on $(\dot E_1\cap (E_1-U_1)\times I$. Then the conditions 
are satisfied for $E_1$, so we have
$$\frak h(E_1):F(\gamma)T(E_1\downarrow (B,A)) \ \sim \ F(\gamma\Delta)\varepsilon(B,A)$$
We will define
$$\frak h(E/E_1):L(B,A)\wedge (I^2)^+\longrightarrow F(X,*) \tag 13.1$$
with edge path
$$
\CD
F(\gamma)T(E\downarrow (B,A)) @>\frak r(E/E_1)> \sim > F(\gamma)T(E_1\downarrow (B,A))\\
@V \sim V \frak h(E) V @V \sim V \frak h(E_1) V\\
F(\gamma\Delta)\varepsilon(B,A) @> 1 > \sim > F(\gamma\Delta)\varepsilon(B,A)\\
\endCD
$$
For $E_1$ we can take $P_1=P$ and $\widehat P_1=\widehat P$. Then $R_1=R$, $\widehat R_1=\widehat R$, and  
$\frak d(E_1)=\frak d(E)$. We have $\widehat R<Q_1<Q$, and $\frak r(Q/Q_1/\widehat R)$ with edge path
$$
\CD
F(\gamma)T((R<Q)\downarrow (B_2,B_0)) @>\frak r(Q/Q_1)> \sim > F(\gamma)T((R<Q_1)\downarrow (B_2,B_0))\\
@V \sim V \frak r(Q/\widehat R) V @V \sim V \frak r(Q_1/\widehat R)V\\
F(\gamma)T((R<\widehat R)\downarrow (B_2,B_0)) @> 1 > \sim > F(\gamma)T((R<\widehat R)\downarrow (B_2,B_0))\\
\endCD
$$

Define 
$\frak c(E/E_1)$ by
$$
\frak c(E/E_1)(\_,t_1,t_2)=
\cases
\frak r(Q/Q_1/\widehat R)(\_,t_1,2t_2) \qquad &0\le t_2\le 1/2\\
\frak d(E)(\_,2t_2-1) \qquad &1/2\le t_2\le 1\\
\endcases
$$
In the diagram
$$
\CD
F(\gamma)T((R<Q)\downarrow (B,B_1))L(b) @>\frak r(Q/Q_1)L(b)> \sim > F(\gamma)T((R<Q_1)\downarrow (B,B_1))L(b)\\
@V \sim V F(\gamma)T(\widetilde b,b)^{-1} V @V \sim V F(\gamma)T(\widetilde b,b)^{-1} V\\
F(\gamma)T((R<Q)\downarrow (B_2,B_0)) @> \frak r(Q/Q_1) > \sim > F(\gamma)T((R<Q)\downarrow (B_2,B_0))\\
@V \sim V \frak c(E) V @V \sim V \frak c(E_1) V\\
* @> 1 > \sim > *\\
\endCD
$$
the top square is filled in by $\frak r(\widetilde b/\widetilde b_1)$, \ the bottom by $\frak c(E/E_1)$. This defines
$$\widehat{\frak c}(E/E_1):L(B_2,B_0)\wedge(I^2)^+\longrightarrow F(X,*)$$
with edge path
$$
\CD
F(\gamma)T((R<Q)\downarrow (B,B_1))L(b) @>\frak r(Q/Q_1)L(b)> \sim > F(\gamma)T((R<Q_1)\downarrow (B,B_1))L(b)\\
@V \sim V \widehat{\frak c}(E) \ \sim \ \widehat{\widehat{\frak c}}(E)L(b) V @V \sim V \widehat{\frak c}(E_1) \ \sim \ \widehat{\widehat{\frak c}}(E_1)L(b) V\\
* @> 1 > \sim > *\\
\endCD
$$
Since $L(b)$ is a homotopy equivalence, $\widehat{\frak c}(E/E_1)$ defines $\widehat{\widehat{\frak c}}(E/E_1)$
with edge path
$$
\CD
F(\gamma)T((R<Q)\downarrow (B,B_1)) @>\frak r(Q/Q_1)> \sim > F(\gamma)T((R<Q_1)\downarrow (B,B_1))\\
@V \sim V \widehat{\widehat{\frak c}}(E) V @V \sim V \widehat{\widehat{\frak c}}(E_1) V\\
* @> 1 > \sim > *\\
\endCD
$$
Then 
$\frak b(E/E_1)$ 
with edge path
$$
\CD
F(\gamma)T(R<Q) @>\frak r(Q/Q_1)> \sim > F(\gamma)T(R<Q_1)\\
@V \sim V \frak b(E) V @V \sim V \frak b(E_1) V\\
* @> 1 > \sim > *\\
\endCD
$$
is defined by 
$$
\CD
F(\gamma)T(R<Q) @>\frak r(Q/Q_1)> \sim > F(\gamma)T(R<Q_1)\\
@V \sim V F(\gamma)T(\widetilde a,a) V @V \sim V F(\gamma)T(\widetilde a,a) V\\
F(\gamma)T((R<Q)\downarrow (B,B_1))L(a) @> \frak r(Q/Q_1)L(a) > \sim > F(\gamma)T((R<Q)\downarrow (B,B_1))L(a)\\
@V \sim V \widehat{\widehat{\frak c}}(E)L(a) V @V \sim V \widehat{\widehat{\frak c}}(E_1)L(a) V\\
* @> 1 > \sim > *\\
\endCD
$$
where $\frak r(\widetilde a,a)$ fills in the top square, and $\widehat{\widehat{\frak c}}(E/E_1)L(a)$
the bottom square.
 
Finally, $\frak h(E/E_1)$ is defined by
$$
\CD
F(\gamma)T(E) @>\frak r(E/E_1)> \sim > F(\gamma)T(E_1)\\
@V \sim V V @V \sim V V\\
F(\gamma)T(P)+F(\gamma)T(P<E) @>1+\frak r(E/E_1)> \sim > F(\gamma)T(P)+F(\gamma)T(P<E_1)\\
@V \sim V V @V \sim V V\\
F(\gamma)T(P)+F(\gamma)T(R<Q) @>1+\frak r(Q/Q_1)> \sim > F(\gamma)T(P)+F(\gamma)T(R<Q_1)\\
@V \sim V 1+\frak b(E) V @V \sim V 1+\frak b(E_1) V\\
F(\gamma)T(P)+* @> 1 > \sim > F(\gamma)T(P)+*\\
\endCD
$$
where (9.2) and (9.3) are used to fill in the top two squares.

\bigskip
If we are given  $E_1<E_2<E$, and
$$\text{supp}(\gamma)\subset E_1-\text{fr}(E_1<E),\qquad \Delta(B)\subset E_1-\text{fr}(E_1<E)$$
we have
$$\frak h(E/E_2/E_1):L(B,A)\wedge (I^3)^+\longrightarrow F(X,*) \tag 13.2$$
such that
$$
\frak h(E/E_2/E_1)(\_,t_1,t_2,t_3)=
\cases
\frak r(E/E_2/E_1)(\_,t_1,t_2) \qquad & t_3=0\\
F(\gamma\Delta)\varepsilon(\_) \qquad & t_3=1\\
\frak h(E/E_1)(\_,t_2,t_3)  \qquad & t_1=0\\
\frak h(E_2/E_1)(\_,t_2,t_3) \qquad & t_1=1\\
\frak h(E/E_2)(\_,t_1,t_3) \qquad & t_2=0\\
\frak h(E_1)(\_,t_3) \qquad & t_2=1\\
\endcases
$$
It is defined from $\frak r(Q/Q_2/Q_1/\widehat R)$ in the same manner as $\frak h(E/E_1)$ from 
$\frak r(Q/Q_1/\widehat R)$.

\bigskip
\subhead 14. Definition of $\Lambda$ \endsubhead We will define 
$$\Lambda:T(E\downarrow (B,A)) \ \sim \ \overline T(E\downarrow (B,A)) \tag 14.1$$ 
Again, to simplify notation we will consider the absolute case. Let 
$\alpha:E\rightarrow B\times \Bbb R^s$ be a fiberwise embedding with normal bundle $N_{\alpha}$, let 
$\alpha:N_{\alpha}\rightarrow B\times \Bbb R^s$ also denote a regular extension of $\alpha$.
Let $Y=D(N_{\alpha})$, \ $p_Y:Y\rightarrow E$ \ the projection, and let $\ddot Y$ denote the topological boundary of Y, 
$$\ddot Y=(D(N_{\alpha})\downarrow \dot E)\cup S(N_{\alpha})$$ 
Let $E\subset Y$ by means of the zero section, and let
$$W=Y\times_B E @>p_W>> Y,\qquad p_W(y,e)=y \tag 14.2$$
We will apply the construction of the previous section to $(W,W|\ddot Y)$. Let \newline
$\theta_E:\dot E\times I\rightarrow E$  be a fiberwise collar with 
$\theta_E(\dot e,1)=\dot e$, and let 
$$\rho:E\rightarrow E(3/4)$$ 
denote the canonical retraction. Let 
$\theta:\dot W\times I\rightarrow W$  denote the induced fiberwise collar. We have
$$W(s)=Y\times_B E(s),\qquad W(s,t)=Y\times_B E(s,t)$$
Define a decomposition $\{Y_1,Y_2\}$ of $Y$ by
$$Y_1=p_Y^{-1}(E(1/4)),\qquad Y_2=p_Y^{-1}(E(1/4,1))$$
 
Let $\bar\alpha$ denote \, $N_{\alpha}@>\alpha>> B\times \Bbb R^s @>\text{proj.}>>\Bbb R^s$, \, and let $\epsilon^*>0$ be 
chosen such that:
\bigskip
\noindent
(14.3) If $(y,e)\in (Y|\dot E)\times_B E(3/4)$, \  then $|\bar\alpha(y)-\bar\alpha(e)|>\epsilon^*$. 
\medskip
\noindent
(14.4) If $(y,e)\in S(N_{\alpha})\times_B E$, then $|\bar\alpha(y)-\bar\alpha(e)|>\epsilon^*.$
\bigskip
\noindent

Let $\zeta_{\epsilon}:W\rightarrow S^s$ be given by
$$\zeta_{\epsilon}(y,e)=
\cases \dsize\frac{\bar\alpha(y)-\bar\alpha(e)}{\epsilon-|\bar\alpha(y)-\bar\alpha(e)|},&
\quad |\bar\alpha(y)- \bar\alpha(e)|<\epsilon \\
       \infty,&\quad |\bar\alpha(y)-\bar\alpha(e)|\ge\epsilon.
\endcases$$
and define for $)<\epsilon<\epsilon^*$
$$\gamma_{\epsilon}:(W,W\downarrow \ddot Y)\longrightarrow (E^+\wedge S^s,*),
\qquad \gamma_{\epsilon}(y,e)=e\wedge\zeta_{\epsilon}(y,\rho(e)) \tag 14.5$$ 
Note that by (8.3) and (8.4), \ $\gamma_{\epsilon}$ sends $W\downarrow \ddot Y$ to the basepoint 
as indicated.

There is the section
$$\Delta:Y\rightarrow W, \qquad \Delta(y)=(y,\rho p_Y(y)) \tag 14.6$$

\proclaim{(14.7) Lemma} Relative to the decomposition  \, $\{Y_1,Y_2\}$, \, the family $\gamma_{\epsilon}$ 
satisfies the conditions {\bf I, II} of section 10.
\endproclaim
\demo{Proof} If $C$ is a closed subset of $Y$, let
$$\Psi(C)=\{ \ (y,e) \in C\times_B E\ \ | \ \ y=\rho(e) \ \}$$
By compactness, if $N$ is a neighborhood of $\Psi (C)$, there is $\epsilon>0$
such that
$$\Psi(C)\subset\text{supp}(\gamma_{\epsilon}|C)\subset N$$
To establish {\bf I}, it is then sufficient to show that 
$$\Psi(Y_1)\subset \Delta(Y_1)$$
Let $(y,e)\in \Psi(Y_1)$. Then $y=\rho(e)$, so $\rho(e)\in E(1/4)$ and therefore
$\rho(e)=e$. Then $y=e$, so $\rho p_Y(y)=e$, and we have $(y,e)=\Delta(y)$.

\medskip
For {\bf II}, it is sufficient to show that
$$\Psi(Y_2)\subset \theta (\delta(Y_2))\times (0,1])$$
Let $(y,e)\in \Psi(Y_2)$. Then $y=\rho(e)$ implies that
$\rho(e)\in E(1/4,3/4)$, so $e\in E(1/4,1)$. Write $e=\theta_E(\dot e,s)$, $\dot e\in \dot E$,
$s\in [1/4,1]$. Then $\rho(e)=\theta_E(\dot e,s')$ for some $s'$,
so $y=\rho(e)=\theta_E(\dot e,s')$, which implies that $\delta(y)=(y,\dot e)$.
Then 
$$(y,e)=(y,\theta_E(\dot e,s))=\theta((y,\dot e),s)=\theta(\delta(y),s)
\in \theta(\delta(Y_1)\times (0,1])$$
\enddemo

Define
$$\Gamma(W):F(\gamma)T(W\downarrow (Y,\ddot Y))\ \sim \ F(\gamma)\overline T(W\downarrow (Y,\ddot Y))\qquad\qquad 
\Gamma(W)=\frak h(W)\cdot\overline{\frak h}(W)^{-1} \tag 14.8$$
Let
$$\widehat W=Y/_B\ddot Y\times_B E @>\widehat w>> Y/_B\ddot Y,\qquad 
\widehat w(\hat y,e)=\hat y$$
and consider
\medskip 
$$
\CD
(W, W|\ddot Y) @>\widetilde k>> (\widehat W,\widehat W|B) @>\widehat\gamma >> (E^+\wedge S^s,*) \\
@VV w V @VV\widehat w V @.\\
(Y, \ddot Y)@>k>> (Y/_B\ddot Y,B)@. \\ 
\endCD 
$$
\bigskip
\noindent
where $\widehat\gamma$ is the quotient map of $\gamma$. By naturality we obtain from $\Gamma$, 
$$\Gamma^{(1)}:F(\widehat\gamma)T(\widehat W,\widehat W\downarrow B)L(k) \ \sim \ 
F(\widehat\gamma)\overline T(\widehat W,\widehat W\downarrow B)L(k)
\tag 14.9$$
and since $L(k)$ is a homotopy equivalence, the latter defines a path 
$$\Gamma^{(2)}:F(\widehat\gamma)T(\widehat W,\widehat W\downarrow B) \, \sim \,
F(\widehat\gamma)\overline T(\widehat W,\widehat W\downarrow B) \tag 14.10$$

We have
$$
\CD
E\times (S^s,\infty) @>\widetilde \alpha!>> (\widehat W,\widehat W|B) @>\widehat\gamma >> (E^+\wedge S^s,*)\\
@VV{p\times 1}V @VV\widehat w V @.\\
B\times (S^s,\infty) @>\alpha! >> (Y/_B\dot Y,B) @.\\ 
\endCD 
$$
where $\alpha!$ is the umkehr map associated to the embedding $\alpha$. By naturality, we obtain   
$$\Gamma^{(3)}:F(\widehat\gamma\widetilde\alpha!)T(E\times (S^s,\infty)) \, \sim \,  
F(\widehat\gamma\widetilde\alpha!)\overline T(E\times (S^s,\infty))   \tag 14.11$$

Let $q:E\times (S^s,\infty)\rightarrow (E^+\wedge S^s,*)$ denote the quotient map. A specific homotopy $H$ from q 
to $\widehat\gamma\widetilde\alpha!$ is given by
$$H_t(e,x)=
\cases
e\wedge \dsize\frac{x-t\pi\alpha\rho_1(e)}{\epsilon(t)-t|x-t\pi\alpha\rho_1(e)|}, \qquad 
t|x-t\pi\alpha\rho_1(e)|<\epsilon(t),\\
*, \qquad t|x-t\pi\alpha\rho_1(e)|\ge \epsilon(t),\\
\endcases
$$ 
where $\pi:B\times \Bbb R^s\rightarrow\Bbb R^s$ is projection and $\epsilon(t)=1-t+t\epsilon$. Then 
$\widehat F(H)$, together with $\Gamma^{(3)}$ define
$$\Gamma^{(4)}:F(q)T(E\times (S^s,\infty)) \, \sim \,  F(q)\overline T(E\times (S^s,\infty)) \tag 14.12$$         
Since $F(q)$ is a homeomorphism the latter defines 
$$\Gamma^{(5)}:T(E\times (S^s,\infty)) \, \sim \,  \overline T(E\times (S^s,\infty))\tag 14.13$$
The multiplicativity paths
$$T(E,(S^s,\infty)):\omega_F\left(T(E)\wedge 1)\right)\ \sim T(E\times (S^s,\infty))\omega_L$$
and 
$$\overline T(E,(S^s,\infty)):\omega_F\left(\overline T(E)\wedge 1)\right) \, \sim \,
\overline T(E\times (S^s,\infty))\omega_L$$
together with $\Gamma^{(5)}$, define
$$\Gamma^{(6)}:\omega_F\left(T(E)\wedge 1)\right) \, \sim \, \omega_F\left(\overline T(E)\wedge 1)\right)\tag 14.14$$
Finally, since $F$ is a homology theory, $\Gamma^{(6)}$ defines a path class 
$$\Lambda:T(E\downarrow (B,A)) \, \phom \, \overline T(E\downarrow (B,A))$$

\bigskip
Given $(E_0<E_1)\downarrow (B,A)$, let $\bold E_{01}=(\Phi<E_0<E_1)\downarrow (B,A)$. Define
\smallskip
$$\Lambda((E_0<E_1)\downarrow (B,A)):T((E_0<E_1)\downarrow (B,A)) \ \phom \ \overline T((E_0<E_1)\downarrow (B,A)) \tag 14.15$$
so that the following commutes.
\medskip
$$
\CD
T(E_1) @>T(\bold E_{01})> \phom > F(i_{01})T((E_0)+T((E_0<E_1)\\
@V \phom V \Lambda(E_1) V         @V \phom V F(i_{01})\Lambda(E_0) + \Lambda(E_0<E_1) V\\
\overline T(E_1) @>\overline T(\bold E_{01})> \phom > F(i_{01})\overline T((E_0)+\overline T((E_0<E_1)\\
\endCD \tag 14.16
$$

\bigskip
\subhead 15. Proof of theorem (1.5) \endsubhead The uniqueness of $\Lambda$ follows from the fact that it is constructed from paths of type {\bf I--IV}.   

\medskip
{\bf I}. Suppose that $E$ is an orthogonal disk bundle and \, $\Delta_E:B@>>>E-\dot E$ \, is a section. We wish to show commutativity of the diagram
$$
\CD
T(E\downarrow (B,A)@>\Lambda(E)> \phom >\overline T(E\downarrow (B,A)\\
@V \phom V \bold d(E,\Delta_E) V   @V \phom V \overline{\bold d}(E,\Delta_E) V\\
F(\Delta_E)\varepsilon(B,A) @>1>\phom> F(\Delta_E)\varepsilon(B,A)\\
\endCD \tag 15.1
$$
By theorem (12.1), $\Gamma(W)$ in (14.8) is given by
$$\Gamma(W)=F(\gamma)\bold d(W,\Delta)\cdot F(\gamma)\overline{\bold d}(W,\Delta)^{-1}  \tag 15.2$$
Let $\Delta_W$ denote the section of $W$ corresponding to $\Delta_E$. By (7.6),
$$\Gamma(W)=F(\gamma)\bold d(W,\Delta_W)\cdot F(\gamma)\overline{\bold d}(W,\Delta_W)^{-1}  \tag 15.3$$
Following the same procedure by which $\Lambda(E)$ is obtained from $\Gamma(W)$, we obtain from (15.3) 
$$\Lambda(E)=\bold d(E,\Delta_E)\cdot \overline{\bold d}(E,\Delta_E)^{-1} \tag 15.4$$

\medskip
{\bf II}. A general excision can be written as a composition of a fiberwise diffeomorphism and a strict excision, that is, one of the form

$$\Cal E= \matrix
E_1 & < & E\\
\vee & & \vee\\
E_0 & < & E_2\\
\endmatrix 
$$
The commutativity of $\Lambda$ with fiberwise diffeomorphisms is easily verified. A proof that $\Lambda$ commutes with a strict excision is given 
in the following section.

\medskip
{\bf III}. Given $E_{012}=E_0<E_1<E_2$, we wish to show that 
$$
\CD
T(E_0<E_2) @>T(E_{012})> \phom > F(i_{12})T(E_0<E_1)+T(E_1<E_2)\\
@V \phom V \Lambda(E_0<E_2) V         @V \phom V F(i_{12})\Lambda(E_0<E_1) + \Lambda(E_1<E_2) V\\
\overline T(E_0<E_2) @>\overline T(E_{012})> \phom > F(i_{12})\overline T(E_0<E_1)+\overline T(E_1<E_2)\\
\endCD \tag 15.5
$$
commutes. The associativity map $T(\Phi<E_0<E_1<E_2)$ determines a commutative diagram
$$
\CD
T(E_2)@>T(\bold E_{02})>\phom> F(i_{02})T(E_0)+T(E_0<E_2)\\
@V \phom V T(\bold E_{12}) V     @V \phom V 1+T(E_{012}) V\\
F(i_{12})T(E_1)+T(E_1<E_2)@>F(i_{12})T(\bold E_{01})+1>\phom>
F(i_{02})T(E_0)+F(i_{12})T(E_0<E_1)+T(E_1<E_2)\\
\endCD
$$
where $\bold E_{ij}=\Phi<E_i<E_j$. Since by definition each $T(\bold E_{ij})$ commutes with $\Lambda$, so also does $T(E_{012})$.

\medskip
{\bf IV}. Verification that $\Lambda$ commutes with the multiplicativity path is straightforward.

\bigskip
\subhead 16. Excision property \endsubhead  Let
$$\Cal E=
\matrix
E_1 &< & E\\
\vee & & \vee\\
E_0 &< & E_2\\
\endmatrix 
$$
be an excisive diagram. Considering the absolute case for notational convenience, we wish to verify homotopy
commutativity of the diagram
$$\CD 
F(i_2) T(E_0<E_2)@> T(\Cal E)>\sim> T(E_1<E)\\
@V\Lambda V\sim V    @V\Lambda V \sim V\\
F(i_2) \overline T(E_0<E_2)@> \overline T(\Cal E)>\sim>  \overline T(E_1<E)\\
\endCD \tag 16.1
$$

Let $F_j=\text{fr}(E_j\subset E)$. Then $F_0=F_1\sqcup F_2$. Let \ $F_0\times [-2,2]\longrightarrow E$ \
be a neat bicollar such that $(e,0)\rightarrow e$, and  $F_0\times [0,2]\subset E_0$. We will identify $F_0\times [-2,2]$ with its image in $E$. Let $\eta:\dot F_0\times I\rightarrow F_0$ be a collar. We have $\dot F_0\times [-2,2]\subset \dot E$, and we extend
$\eta$ to 
$$\eta':(\dot F_0\times [-2,2])\times I\longrightarrow E$$
by $\eta'((e,s),t)=(\eta(e,t),s)$. We may choose collars
$$\bar\theta:\dot E\times I\longrightarrow E,\qquad \bar\theta_j:\dot E_j\times I\longrightarrow E_j, \quad j=0,1,2 $$
such that
\medskip
\itemitem{(i)} $\bar\theta$ extends $\eta'$, and \ $\bar\theta(e,t)\not\in F_0\times [-2,2]$, \ $e\not\in \dot F_0\times [-2,2]$.
\medskip
\itemitem{(ii)} $\bar\theta_j(e,t)=\bar\theta(e,t)$, \ $e\in\dot E_j-(F_j\times [0,1))$, and 
$\bar\theta_j(e,t)=\bar\theta_0(e,t)$, \ $e\in \dot F_j\times [0,1]$.
\medskip
\itemitem{(iii)} $\bar\theta_j(e,t)=\bar\theta_0(e,t)$, \ $e\in \dot F_j\times [0,1]$.
\medskip
\itemitem{(iv)} $\bar\theta_0(e,t)\in F_0\times [0,1]$, \ $e\in \dot F_0\times [0,1]$.
\medskip
\noindent
Let $\bar\rho:E\rightarrow E(3/4)$ and $\bar\rho_j:E_j\rightarrow E_j(3/4)$, denote the associated retractions.
\bigskip
Let $\alpha:E\rightarrow B\times \Bbb R^s$ be a fiberwise embedding. With the notation of section 14, we have
$$Y=D(N_{\alpha}), \qquad \ddot Y=S(N_{\alpha})\cup D(N_{\alpha})\downarrow \dot E$$
Let 
$$Y_j=Y\downarrow E_j, \qquad \ddot Y_j=S(N_{\alpha})\downarrow E_j\cup D(N_{\alpha})\downarrow \dot E_j$$
$$W_j=Y_j\times_B E_j, \qquad W^j=Y\times_B E_j$$
There are the families 
$$\gamma:(W\downarrow(Y,\ddot Y))\longrightarrow (E^+\wedge S^s,*),\qquad \gamma(y,e)=
e\wedge\zeta_{\epsilon}(y,\bar\rho(e)) $$
$$\gamma_j:(W_j\downarrow(Y_j,\ddot Y_j))\longrightarrow (E^+\wedge S^s,*),\qquad
\gamma_j(y,e)=e\wedge\zeta_{\epsilon}(y,\bar\rho_j(e))$$
\smallskip
$$\gamma^j:(W^j\downarrow(Y,\ddot Y))\longrightarrow (E^+\wedge S^s,*),\qquad
\gamma^j(y,e)=e\wedge\zeta_{\epsilon}(y,\bar\rho_j(e))$$

\bigskip
A path 
$$\mu:\gamma \ \sim \ \gamma^j:W^j\downarrow (Y,\ddot Y)\longrightarrow (E^+\wedge S^s,*)$$ 
is defined as follows. Let $\kappa:E_j\times I\rightarrow E_j$, denote
the standard deformation retraction of $E_j$ onto \newline $E_j-(F_j\times [0,1))$. The conditions on the collars imply that 
$\bar\rho\kappa(e,1)=\bar\rho_j\kappa(e,1)$, \ $e\in E_j$. So we have $\nu:\bar\rho \sim \bar\rho_j$ by
$$
\nu(e,t)=
\cases
\bar\rho\kappa(e,2t),\qquad &0\le t\le 1/2\\
\bar\rho_j\kappa(e,2-2t),\qquad &1/2\le t\le 1\\
\endcases
$$
Let $\mu(y,e,t)=e\wedge\zeta(y,\nu(e,t))$.

\bigskip
Let $A\subset Y$. We have 
$$\Gamma(W):F(\gamma)T(W\downarrow (A,A\cap\ddot Y)) \ \sim \ F(\gamma)\overline T(W\downarrow (A,A\cap\ddot Y))\tag 16.2$$
$\Gamma(W)=\frak h(W)\cdot\overline{\frak h}(W)^{-1}$, \, and
$$\Gamma(W_j):F(\gamma_j)T(W_j\downarrow (A\cap Y_j,A\cap \ddot Y_j)) \ \sim \ 
F(\gamma_j)\overline T(W_j\downarrow (A\cap Y_j,A\cap \ddot Y_j))\tag 16.3$$
$\Gamma(W_j)=\frak h(W_j)\cdot\overline{\frak h}(W_j)^{-1}$. We derive a path
$$\Gamma(W^j):F(\gamma)T(W^j\downarrow (A,A\cap\ddot Y)) \ \sim \ F(\gamma)\overline T(W^j\downarrow (A,A\cap\ddot Y))
\tag 16.4$$
as follows. Let $C_j=\ddot Y_j\cup(Y-Y_j)$. We have a commutative diagram

$$\CD
(E^+\wedge S^s,*) @> 1 >> (E^+\wedge S^s,*) @> 1 >> (E^+\wedge S^s,*)\\
@A \gamma^j AA  @A \gamma^j AA   @A \gamma_j AA \\
W^j\downarrow (A,A\cap\ddot Y) @> \tilde a>> W^j\downarrow (A,A\cap C_j) @< \tilde b <<
W_j\downarrow (A\cap Y_j,A\cap  \ddot Y_j)\\
@VVV        @VVV       @VVV\\
(A,A\cap\ddot Y) @> a >> (A,A\cap C_j) @< b << (A\cap Y_j,A\cap \ddot Y_j)\\
\endCD$$
By naturality and the fact that $b$ is an excision, $\Gamma(W_j)$ in (10.3) defines 
$$\alpha:F(\gamma^j)T(W^j\downarrow (A,A\cap\ddot Y)) \ \sim \ F(\gamma^j)\overline T(W^j\downarrow (A,A\cap\ddot Y))$$
Then \ $\Gamma(W^j)=\beta\cdot\alpha\cdot\overline\beta^{-1}$, \ where \
$\beta=\widehat F(\mu_j)T(W^j\downarrow (A,A\cap\ddot Y))$. 

\bigskip
Together with the additivity path, $\Gamma(W)$ and $\Gamma(W^j)$ define
\bigskip
$$\Gamma(W^j<W):F(\gamma)T((W^j<W)\downarrow (A,A\cap\ddot Y))\ \sim \
F(\gamma)\overline T((W^j<W)\downarrow (A,A\cap\ddot Y)) \tag 16.5$$

$$\Gamma(W^i<W^j):F(\gamma) T((W^i<W^j)\downarrow (A,A\cap\ddot Y)) \ \sim \
F(\gamma)\overline T((W^i<W^j)\downarrow (A,A\cap\ddot Y)) \tag 16.6$$

Let
$$\Cal W=
\matrix
W^1 &< & W\\
\vee & & \vee\\
W^0 &< & W^2\\
\endmatrix 
$$
It is straightforward to show that the homotopy commutativity of (16.1) follows from that of
\bigskip

$$\CD 
F(\gamma)T((W^0<W^2)\downarrow (Y,\ddot Y))@> F(\gamma)T(\Cal W)>\sim> 
F(\gamma)T((W^1<W)\downarrow(Y,\ddot Y))\\
@V\sim V\Gamma(W^0<W^2) V  @V \sim V\Gamma(W^1<W) V\\
F(\gamma)\overline T((W^0<W^2)\downarrow (Y,\ddot Y))@>F(\gamma)\overline T(\Cal W)>\sim> 
F(\gamma)\overline T((W^1<W)\downarrow(Y,\ddot Y))\\
\endCD \tag 16.7
$$
Suppressing $(Y,\ddot Y)$ from the notation, the homotopy commutativity of the latter follows from the homotopy commutativity of
$$\CD 
F(\gamma)T(W^0<W)@> F(\gamma)\bold a(\Cal W)>\sim> 
F(\gamma)T(W^0<W^1)+ F(\gamma)T(W^0<W^2)\\
@V\sim V\Gamma(W^0<W) V   @V\sim V\Gamma(W^0<W^1) + \Gamma(W^0<W^2) V\\
F(\gamma)\overline T(W^0<W)@>F(\gamma)\overline{\bold a}(\Cal W)>\sim> 
F(\gamma)\overline T(W^0<W^1)+F(\gamma)\overline T(W^0<W^2)\\
\endCD \tag 16.8
$$

\bigskip
Let
$$X_j=Y\downarrow (E_j-(F_j\times[0,2))),\qquad j=0,1,2$$
We will define for $A\subset X_j$
$$H^j(\Cal W\downarrow (A,A\cap\ddot Y)):L(A,A\cap\ddot Y)\wedge (I^2)^+\longrightarrow F(E^+\wedge S^s,*), \qquad j=0,1,2 
\tag 16.9$$
having edge-path diagram (16.8). 

\bigskip
{\it Definition of  $H^1(\Cal W)$.}  Bundles are over $(A,A\cap\ddot Y)$, where $A\subset X_1$.
Let
$$H^1(\Cal W^1):L(A,A\cap\ddot A)\wedge (I^2)^+\longrightarrow F(E^+\wedge S^s,*) \tag 16.10$$
be the standard map with edge path 
$$\CD 
F(\gamma)T(W^0<W^1)@> \frak s >\sim> 
F(\gamma)T(W^0<W^1)+ *\\
@V\sim V\Gamma(W^0<W^1) V   @V\sim V\Gamma(W^0<W^1) + 1 V\\
F(\gamma)\overline T(W^0<W^1)@> \frak s>\sim> 
F(\gamma)\overline T(W^0<W^1)+*\\
\endCD 
$$
 
We have  
$$\text{supp}(\gamma|W\downarrow A)\subset W^1-\text{fr}(W^1\subset W)$$ 
Consequently, $\frak h(W^1)$ is defined and moreover, \, $\Gamma(W^1)=\frak h(W^1)\cdot \overline{\frak h}(W^1)^{-1}$. 
From (13.1), we have $\frak h(W/W^1)$. Define $\gamma(W/W^1)$ by
$$\Gamma(W/W^1)(\_,t_1,t_2)=
\cases
\frak h(W/W^1)(\_,t_2,2t_1)\qquad &0\le t_1\le 1/2\\
\frak h(W/W^1)(\_,t_2,2-2t_1)\qquad &1/2\le t_1\le 1\\
\endcases
$$
It has edge path

$$
\CD
F(\gamma)T(W) @>\Gamma(W)> \sim > F(\gamma)\overline T(W)\\
@V \sim V \frak r(W/W^1) V   @V\sim V \overline{\frak r}(W/W^1) V \\
F(\gamma)T(W^1) @>\Gamma(W^1) > \sim > F(\gamma)\overline T(W^1)\\
\endCD
$$

Let
$$\bold W=W_0<W\qquad\qquad \bold W^j=W^0<W^j,\quad j=0,1,2$$
In the diagram
$$
\CD
F(\gamma)T(W^0<W) @>\frak r(\bold W/\bold W^1)> \sim > F(\gamma)T(W^0<W^1)\\
@V \sim VV               @V \sim VV\\
-F(\gamma)T(W^0)+F(\gamma)T(W) @>1+\frak r(W/W^1)> \sim > -F(\gamma)T(W^0)+F(\gamma)T(W^1)\\
@V \sim V -\Gamma(W^0)+\Gamma(W)  V   @V\sim V-\Gamma(W^0)+\Gamma(W^1)V \\
-F(\gamma)\overline T(W^0)+F(\gamma)\overline T(W) @>1+\overline {\frak r}(W/W^1)> \sim > 
-F(\gamma)\overline T(W^0)+F(\gamma)\overline T(W^1)\\
@V \sim VV               @V \sim VV\\
F(\gamma)\overline T(W^0<W) @>\overline{\frak r}(\bold W/\bold W^1) > \sim > F(\gamma)\overline T(W^0<W^1)\\
\endCD
$$
the map (8.2) is used to fill in the top and bottom square. This defines $\Gamma(\bold W/\bold W^1)$ \ with edge path
$$
\CD
F(\gamma)T(W^0<W) @>\Gamma(W^0<W)> \sim > F(\gamma)\overline T(W^0<W)\\
@V \sim V \frak r(\bold W/\bold W^1) V   @V\sim V \overline{\frak r}(\bold W/\bold W^1) V \\
F(\gamma)T(W^0<W^1) @>\Gamma(W^0<W^1) > \sim > F(\gamma)\overline T(W^0<W^1)\\
\endCD
$$
Similarly, since
$$\text{supp}(\gamma|W^2\downarrow A)\subset W^0-\text{fr}(W^0 < W^2)$$
we have \ $\Gamma(\bold W^2/\bold W^0)$, \ with edge path
$$
\CD
F(\gamma)T(W^0<W^2) @> \Gamma(W^0<W^2) > \sim > F(\gamma)\overline T(W^0<W^2)\\
@V \sim V\frak r(\bold W^2/\bold W^0) V   @V\sim V \overline{\frak r}(\bold W^2/\bold W^0) V \\
F(\gamma)T(W^0<W^0) @>\Gamma(W^0<W^0) > \sim > F(\gamma)\overline T(W^0<W^0)\\
\endCD
$$
Let  
$$ K^1(\Cal W):L(A,A\cap\ddot Y)\wedge (I^3)^+\longrightarrow F(E^+\wedge S^s,*)\tag 16.11$$
be such that
$$
K^1(\_,t_1,t_2,t_3)=
\cases
\Gamma(\bold W/\bold W^1)(\_,t_2,t_3)\qquad &t_1=0\\
\Gamma(\bold W^1/\bold W^1)(\_,t_2,t_3)+\Gamma(\bold W^2/\bold W^0)(\_,t_2,t_3)\qquad &t_1=1\\
\bold a(\Cal W/\Cal W^1)(\_,t_1,t_3)\qquad &t_2=0\\
\overline{\bold a}(\Cal W/\Cal W^1)(\_,t_1,t_3)\qquad &t_2=1\\
H^1(\Cal W^1)(\_,t_1,t_2)\qquad &t_3=1\\
\endcases
$$
where $\bold a(\Cal W/\Cal W^1)$ is from (9.2). 
Now define \ $H^1(\Cal W)(\_,t_1,t_2)=K^1(\_,t_1,t_2,0)$.

\bigskip
{\it Definition of  $H^2(\Cal W)$.} Here bundles are over $(A,A\cap\ddot Y)$ where $A\subset X_2$. Let
$$H^2(\Cal W^2):L(A,A\cap\ddot Y)\wedge (I^2)^+\longrightarrow F(E^+\wedge S^s,*) \tag 16.12$$
be the standard map with edge path 
$$\CD 
F(\gamma)T(W^0<W^2)@> \frak s >\sim> 
*+F(\gamma)T(W^0<W^2)\\
@V\sim V\Gamma(W^0<W^2) V   @V\sim V 1+\Gamma(W^0<W^2) V\\
F(\gamma)\overline T(W^0<W^2)@> \frak s>\sim> 
*+F(\gamma)\overline T(W^0<W^2)\\
\endCD 
$$
 
Define 
$$ K^2(\Cal W):L(A,A\cap\ddot Y)\wedge (I^3)^+\longrightarrow F(E^+\wedge S^s,*)\tag 16.13$$
such that
$$
K^2(\_,t_1,t_2,t_3)=
\cases
\Gamma(\bold W/\bold W^2)(\_,t_2,t_3)\qquad &t_1=0\\
\Gamma(\bold W^1/\bold W^0)(\_,t_2,t_3)+\Gamma(\bold W^2/\bold W^2)(\_,t_2,t_3)\qquad &t_1=1\\
\bold a(\Cal W/\Cal W^2)(\_,t_1,t_3)\qquad &t_2=0\\
\overline{\bold a}(\Cal W/\Cal W^2)(\_,t_1,t_3)\qquad &t_2=1\\
H^2(\Cal W^2)(\_,t_1,t_2)\qquad &t_3=1\\
\endcases
$$
where $\bold a(\Cal W/\Cal W^2)$ is from (9.3). Again, define $H^2(\Cal W)(\_,t_1,t_2)=K^2(\_,t_1,t_2,0)$.
\bigskip
{\it Definition of $H^0(\Cal W)$.} Bundles are over $(A,A\cap\ddot Y)$ where $A\subset X_0$. Let
$$ K^0(\Cal W):L(A,A\cap\ddot Y)\wedge (I^3)^+\longrightarrow F(E^+\wedge S^s,*)  \tag 16.14$$
be such that
$$
K^0(\Cal W)(\_,t_1,t_2,t_3)=
\cases
\Gamma(\bold W/\bold W^0)(\_,t_2,t_3)\qquad &t_1=0\\
\Gamma(\bold W^1/\bold W^0)(\_,t_2,t_3)+\Gamma(\bold W^2/\bold W^0)(\_,t_2,t_3)\qquad &t_1=1\\
\bold a(\Cal W/\Cal W^0)(\_,t_1,t_3)\qquad &t_2=0\\
\overline{\bold a}(\Cal W/\Cal W^0)(\_,t_1,t_3)\qquad &t_2=1\\
*\qquad &t_3=1\\
\endcases
$$
where $\bold a(\Cal W/\Cal W^0)$ is from (9.5). Define $H^0(\Cal W)(\_,t_1,t_2)=K^0(\_,t_1,t_2,0)$.

\bigskip
\proclaim{(16.15) Theorem} Assume $A\subset X_0$. Then
$H^0(\Cal W) \sim H^1(\Cal W)$ rel. $\partial(I^2)$.
\endproclaim
\demo{Proof}  Define
$$K^0(\Cal W^1):L(A,A\cap\ddot Y)\wedge (I^3)^+\longrightarrow F(E^+\wedge S^s,*)$$
by
$$K^0(\Cal W^1)(\_,t_1,t_2,t_3)(s)=
\cases
\Gamma(\bold W^1/\bold W^0)(\_,t_2,t_3)((1+t_1)s),\qquad &0\le s\le 1/(1+t_1)\\
*,\qquad &1/(1+t_1)\le s \le 1\\
\endcases
$$
\medskip\noindent
Then $K^1(\Cal W)\cdot K^0(\Cal W^1)$ defines
$$H^1(\Cal W)\ \sim * \qquad \text{rel.} \qquad P:L(A,A\cap\ddot Y)\wedge (\partial(I^2)\times I)^+
\longrightarrow F(E^+\wedge S^s,*) $$
where
$$
P(\_,t_1,t_2,t_3)=
\cases
\Gamma(\bold W/\bold W^1)\cdot\Gamma(\bold W^1/\bold W^0)(\_,t_2,t_3)\qquad &t_1=0\\
\Gamma(\bold W^1/\bold W^1)\cdot\Gamma(\bold W^1/\bold W^0)(\_,t_2,t_3)+
\Gamma(\bold W^2/\bold W^0)\cdot\Gamma(\bold W^0/\bold W^0)(\_,t_2,t_3)\qquad &t_1=1\\
\bold a(\Cal W/\Cal W^1)\cdot\bold a(\Cal W^1/\Cal W^0)(\_,t_1,t_3)\qquad &t_2=0\\
\overline{\bold a}(\Cal W/\Cal W^1)\cdot\overline{\bold a}(\Cal W^1/\Cal W^0)(\_,t_1,t_3)\qquad &t_2=1\\
\endcases
$$

By definition, $K^0(\Cal W)$ defines
$$H^1(\Cal W)\ \sim * \qquad \text{rel.} \qquad Q:L(A,A\cap\ddot Y)\wedge (\partial(I^2)\times I)^+
\longrightarrow F(E^+\wedge S^s,*) $$
where
$$
Q(\_,t_1,t_2,t_3)=
\cases
\Gamma(\bold W/\bold W^0)(\_,t_2,t_3)\qquad &t_1=0\\
\Gamma(\bold W^1/\bold W^0)(\_,t_2,t_3)+\Gamma(\bold W^2/\bold W^0)(\_,t_2,t_3)\qquad &t_1=1\\
\bold a(\Cal W/\Cal W^0)(\_,t_1,t_3)\qquad &t_2=0\\
\overline{\bold a}(\Cal W/\Cal W^0)(\_,t_1,t_3)\qquad &t_2=1\\
\endcases
$$
\medskip\noindent
So it is sufficient to exhibit a homotopy 
$$R:P \, \sim \, Q \qquad \text{rel.}\qquad  L(A,A\cap\ddot Y)\wedge(\partial(I^2)\times\{0,1\})^+ $$

Let $\phi:I^2@>>>I^2$ be given by
$$
\align
\phi(t_1,0)&=
\cases
(2t_1,0),\qquad &0\le t_1\le 1/2\\
(1,2t_1-1),\qquad &1/2\le t_1\le 1\\
\endcases\\
\qquad \phi(t_2,1)&=(0,t_2)\\
\qquad \phi(t_1,t_2)&=(1-t_2)\phi(t_1,0)+t_2\phi(t_1,1)\\
\endalign
$$
Let $\bold a^*(\Cal W/\Cal W^1/\Cal W^0)$ denote the reparametrization of $\bold a(\Cal W/\Cal W^1/\Cal W^0)$
in (10.6) given by
$$\bold a^*(\Cal W/\Cal W^1/\Cal W^0)(\_,t_1,t_2,t_3)=\bold a(\Cal W/\Cal W^1/\Cal W^0)(\_,t_1,\phi(t_2,t_3))$$

We have $\Gamma(\bold W/\bold W^1/\bold W^0)$ such that
$$
\Gamma(\bold W/\bold W^1/\bold W^0)(\_,t_1,t_2,t_3)=
\cases
\frak r(\bold W/\bold W^1/\bold W^0)(\_,t_1,t_2) \qquad & t_3=0\\
\overline{\frak r}(\bold W/\bold W^1/\bold W^0)(\_,t_1,t_2)\qquad & t_3=1\\
\Gamma(\bold W/\bold W^0)(\_,t_2,t_3)  \qquad & t_1=0\\
\Gamma(\bold W^1/\bold W^0)(\_,t_2,t_3) \qquad & t_1=1\\
\Gamma(\bold W/\bold W^1)(\_,t_1,t_3) \qquad & t_2=0\\
\Gamma(\bold W^0)(\_,t_3) \qquad & t_2=1\\
\endcases
$$
It is derived from (13.2) in a similar manner as $\Gamma(\bold W/\bold W^1)$ from (13.1).
Let 
$$\Gamma^*(\bold W/\bold W^1/\bold W^0)(\_,t_1,t_2,t_3)=\Gamma(\bold W/\bold W^1/\bold W^0)(\_,t_1,\phi(t_2,t_3))$$
The required homotopy $R$ is 
$$R(\_,t_1,t_2,t_3,t_4)=
\cases
\Gamma^*(\bold W/\bold W^1/\bold W^0)(\_,t_2,t_3,t_4)\qquad &t_1=0\\
\Gamma^*(\bold W^1/\bold W^1/\bold W^0)(\_,t_2,t_3,t_4)+\Gamma^*(\bold W^2/\bold W^0/\bold W^0)(\_,t_2,t_3,t_4)\qquad &t_1=1\\
\bold a^*(\Cal W/\Cal W^1/\Cal W^0)(\_,t_1,t_3,t_4)\qquad &t_2=0\\
\overline{\bold a}^*(\Cal W/\Cal W^1/\Cal W^0)(\_,t_1,t_3,t_4)\qquad &t_2=1\\
\endcases
$$
\enddemo

\bigskip
By a similar argument we have
\bigskip
\proclaim{(16.16) Theorem} Assume $A\subset X_0$. Then
$H^0(\Cal W) \sim H^2(\Cal W)$ rel. $\partial(I^2)$.
\endproclaim

\bigskip
Let $\Cal C$ denote the category of spaces and inclusions
$$
\CD
(X_1,X_1\cap\ddot Y) @>>> (Y,\ddot Y)\\
@AAA                                     @AAA \\
(X_0,X_0\cap\ddot Y) @>>> (X_2,X_2\cap\ddot Y)\\
\endCD
$$

\bigskip\noindent
Define
$$
\align
S_0:F(\gamma)T(W^0<W)@> F(\gamma)\bold a(\Cal W)>\sim> 
&F(\gamma)T(W^0<W^1)+ F(\gamma)T(W^0<W^2) \tag 16.17\\
@>\Gamma(W^0<W^1) + \Gamma(W^0<W^2)> \sim >
&F(\gamma)\overline T(W^0<W^1)+F(\gamma)\overline T(W^0<W^2)\\ 
\endalign 
$$

$$
\align
S_1:F(\gamma)T(W^0<W)@> \Gamma(W^0<W) >\sim> 
&F(\gamma)\overline T(W^0<W) \tag 16.18\\
@> F(\gamma)\overline{\bold a}(\Cal W)> \sim >
&F(\gamma)\overline T(W^0<W^1)+F(\gamma)\overline T(W^0<W^2)\\ 
\endalign 
$$
\bigskip\noindent
The homotopy commutativity of (10.8) corresponds to $S_0(Y,\ddot Y)\ \sim \ S_1(Y,\ddot Y)$ \, rel $\dot I$.
The paths $S_0$, $S_1$ are defined on the objects of $\Cal C$, and it is straightforward to extend $S_0$, $S_1$ to simplicial
transformations on $\Cal C$. We have defined $H^j(\Cal W):S_0 \ \sim S_1$ \, on the objects of the subcategory
$\{(X_0,X_0\cap\ddot Y)\longrightarrow (X_j,X_j\cap\ddot Y)\}$, \, $j=1,2$, \, and again, it is straightforward to
extend  $H^j(\Cal W)$ to a simplicial tranformation on the subcategory. Moreover, from (16.15) and (16.16),
$$H^1(\Cal W)(X_0,X_0\cap\ddot Y) \ \sim \ H^2(\Cal W)(X_0,X_0\cap\ddot Y)\qquad \text{rel. }\partial(I^2)$$
It follows now from the gluing property discussed in the Appendix, that $S_0(Y,\ddot Y)\ \sim \ S_1(Y,\ddot Y)$ \, rel $\dot I$.
This completes the verification of the excision property.

\bigskip
\subhead 17. $A$-theory transfer \endsubhead In this section we will define $A$-theory transfer \cite{4}  
as a refined transfer. We briefly recall the definition of Waldhausen's space $A(X)$ \cite{8}. Let $R(B,X)$ denote the category of retractive spaces over $X$, parametrized by $B$. The objects are 
$$ (Y,r,s)=
\CD
B\times X @>s>> Y @>r>> B\times X\\
@VV \pi_1 V @VV p V  @VV \pi_1 V\\
B @>1>> B @>1>> B\\
\endCD
$$
where $Y@> p >>B$ is a fibration and $r$, $s$, are fiber preserving maps such that $rs=1$. The morphisms are commutative 
diagrams. It is assumed that the fiber of $Y$ has the homotopy type of a CW-complex rel. $s(X)$. It is a category with 
cofibrations and weak equivalences. As a functor from pairs of spaces to categories, it is contravariant in $B$, and 
covariant in $X$. For $f:B@>>>B'$,
\ $R(f):R(B',X)@>>>R(B,X)$ is given by pullback, and for $g:X@>>>X'$, \ $R(g):R(B,X)@>>>R(B,X')$ is given by pushout.

\medskip
Now $\Cal S_q(R(B,X))$ is defined to be the category with objects functors
$$A:\text{Ar}([q])\longrightarrow R(B,X)$$
such that \ $A_{ii}=*$, for all $i$, and \ $A_{ij}@>>>A_{ik}@>>>A_{jk}$ \ is a cofibration sequence, for all $i\le j\le k$.
The morphisms are natural transformations.

\medskip
Let $R_p(X)=R(\Delta_p,X)$. We have a simplicial map
$$\mu:\Delta(B)\times \Cal N.(w\Cal S.(R(B,X)))\longrightarrow \Cal N.(w\Cal S.(R.(X))), \qquad 
\mu(\sigma,Y)=R(\sigma)(Y) \tag 17.1$$
Let
$$W(X)=|\Cal N.(w\Cal S.(R.(X))| \tag 17.2$$
For $X\in \Cal T_0$, define
$$\widehat W(X)=\frac{\text{Fiber}(W(X)\longrightarrow W(x_0))}
{\text{Fiber}(W(x_0)\longrightarrow W(x_0))}  \tag 17.3$$
For $(X,A)\in \Cal T^2$, define
$$\widetilde A(X,A)=\widehat W(X/A),\qquad A(X,A)=\Omega \widetilde A(X,A) \tag 17.4$$
There is the natural map \ $R(B,X)@>>>S_1(R(B,X))$, \ by \ $(Y,r,s)\longrightarrow(*@>>>(Y,r,s))$. In particular, it
defines a map
$$|wR.(X)|\longrightarrow \Omega W(X) \tag 17.5$$

\bigskip
By definition, $A(X)=A(X,\Phi)=\Omega\widehat W(X^+)$. Let us recall the correspondence with the usual definition of $A(X)$ as 
$\Omega W(X)$. If $Y$ is pointed, \ $\{y_0\}@>i>> Y@>j>>\{y_0\}$, \ we have
$$\lambda:\Omega W(Y) \longrightarrow \text{Fiber}(\Omega W(Y) @>\Omega W(j) >>\Omega W(y_0) ),\qquad\qquad 
\lambda(\rho)=(\rho- W(ij) \rho,\mu)$$
where \ $\mu:*  \sim   W(j) \rho- W(j) \rho$ \ is the canonical path. Then for $X$ unpointed, we have a natural homotopy equivalence
$$\zeta:\Omega W(X) \longrightarrow A(X) \tag 17.6$$
by
$$
\align
\Omega W(X) @>>>\Omega W(X^+) &@>\lambda>> \text{Fiber}(\Omega W(X^+ @>>>\Omega W(+) )\\
     &=\Omega\left(\text{Fiber}( W(X^+) @>>> W(+) )\right)
     @>>>\Omega\left(\frac{\text{Fiber}( W(X^+) \longrightarrow W(+) )}{\text{Fiber}( W(+) \longrightarrow  W(+))}\right)
     =A(X)\\     
\endalign
$$

\bigskip
In $R(X,X)$ there is the retractive space 
$$
\widetilde X=
\CD
X\times X@>s>>X\sqcup (X\times X)@>r>>X\times X\\
@VV \pi_1 V  @VV p V  @VV \pi_1 V\\
X@> 1>> X @>1>> X\\
\endCD
$$
where \ $p=1_X\sqcup \pi_1$, \ $r(x)=(x,x)$, \ $r(x,y)=(x,y)$, \ $s(x,y)=(x,y)$.  Define
$$\widetilde{\varepsilon}:\Delta_p(X)\longrightarrow R_p(X), \qquad\qquad \widetilde{\varepsilon}(\sigma)=
R(\sigma)(\widetilde X)\tag 17.7$$
Then
$$|X|@>|\widetilde{\varepsilon}|>>|wR.(X)|@>>>\Omega W(X)$$ 
defines an augmentation
$$\varepsilon:L(X,A)\longrightarrow A(X,A) \tag 17.8$$

\medskip
Given retractive spaces
$$Y_j=
\CD
X_j\times B @>s_j>> Y_j @>r_j>> X_j\times B\\
@VV \pi_1 V @VV p_j V  @VV \pi_1 V\\
B @>1>> B @>1>> B\\
\endCD \qquad \qquad j=0,1 
$$
the fiber product $Y_0\times_B Y_1$ has the form 
$$Y_0\times_B Y_1=
\CD
(X_0\times X_1)\times B  @>s_{01}>> Y_0\times_B Y_1 @>r_{01}>> (X_0\times X_1)\times B \\
@VV \pi_1 V @VV p_{01} V  @VV \pi_1 V\\
B @>1>> B @>1>> B\\
\endCD  
$$
a retractive space over $X_0\times X_1$ parametrized by $B$. Thus we have a pairing 
$$R(B,X_0)\times R(B,X_1)\longrightarrow R(B,X_0\times X_1) \tag 17.9$$
This in turn defines a pairing
$$A(X_0,A_0)\wedge A(X_1,A_1)\longrightarrow A((X_0,A_0)\times (X_1,A_1)) \tag 17.10$$
A multiplicative structure for $A$, in the sense of section 2,
$$\omega_A:A(X_0,A_0)\wedge L(X_1,A_1)\longrightarrow A((X_0,A_0)\times (X_1,A_1)) \tag 17.11$$
is defined by
$$A(X_0,A_0)\wedge L(X_1,A_1)@>1\wedge\varepsilon>> A(X_0,A_0)\wedge A(X_1,A_1)@>>> A((X_0,A_0)\times (X_1,A_1))$$

\medskip
We turn now to the definition of transfer $T_A$. Suppose given a $qr$-diagram
$$
E=
\CD
(E_{0r}< \dots < E_{qr})\downarrow (B_r,A_r)\\
@AA(\widetilde f_{r-1)r},f_{r-1)r}) A\\
\vdots\\
@AA(\widetilde f_{01},f_{01}) A\\
(E_{00}< \dots < E_{q0})\downarrow (B_0,A_0)\\
\endCD
$$
Define  
$$Y_k(E)\in S_q(R(B_0,E_{qr})),\qquad\quad
Y_k(E)(i,j)=(B_0\times_{B_k} E_{jk})\cup_{B_0\times_{B_k} E_{ik}} (B_0\times E_{qr}) \tag 17.12$$
If $E$ is an excisive diagram, then
$$Y(E)=Y_0(E)@>>> \dots @>>> Y_r(E)\tag 17.13$$
is a sequence of weak equivalences. Thus, we have
$$Y(E)\in \Cal N_r(wS_q(R(B_0,E_{qr})), \qquad\qquad Y(E):[q]\times [r]@>>>\Cal N.(wS.(R(B_0,E_{qr}))$$
Define
$$\Delta(B_0)\times [q]\times [r]\longrightarrow \Cal N.(wS.(R.(E_{qr})) \tag 17.14$$
by
$$\Delta(B_0)\times [q]\times [r]@>1\times Y(E)>>\Delta(B_0)\times \Cal N.(wS.(R(B_0,E_{qr}))@>\mu>>\Cal N.(wS.(R.(E_{qr}))$$
Let
$$t(Y(E)):|B_0|\times \Delta_q\times \Delta_r \longrightarrow W(E_{qr}) \tag 17.15$$
denote its realization. This defines a bisimplicial transformation on the subcomplex of $\Cal N^2(\Cal B)$ consisting of regular excisive diagrams.

\bigskip
Let $e_{j_0},\dots ,e_{j_q}$ be vertices of $\Delta_p$, and let $\phi(e_{j_0},\dots ,e_{j_q}):\Delta_q@>>>\Delta_p$ \
denote the linear map which sends $e_k$ to $e_{j_k}$. We define
$$h^q:I^q@>>> \Delta_q,\qquad\qquad 0\le k\le 4 \tag 17.16$$
as follows. Assuming $h^q$ has been defined, let 
$$h(e_{j_0},\dots ,e_{j_q})=\phi(e_{j_0},\dots ,e_{j_q})h^q:I^q@>>>\Delta_p$$

\item{q=0.} Let $h^0:I^0@>>>\Delta_0$, \ $h^0(0)=e_0$.
\medskip
\item{q=1.} Let $h^1:I@>>>\Delta_1$, \ $h^1(t)=(1-t)e_0+te_1$.
\medskip
\item{q=2.} Let \ $h^2:I^2@>>>\Delta_2$ \ be the linear homotopy such that \
$h^2_0=h(e_0,e_2)$, \ $h^2_1=h(e_0,e_1)\cdot h(e_1,e_2)$.
\medskip
\item{q=3.} let \ $h^3:I^3@>>>\Delta_3$ \ be such that 

\medskip\noindent
($i$) $h^3(t_1,t_2,0)=e_0$, \ $h_3(t_1,t_2,1)=e_3$.

\smallskip\noindent
($ii$) The adjoint of $h^3:\partial(I^2)\times I@>>>\Delta_3$ has edge path
$$
\CD
h(e_0,e_3) @>h(e_0,e_2,e_3)>\sim> h(e_0,e_2)\cdot h(e_2,e_3)\\
@V\sim V h(e_0,e_1,e_3) V   @V\sim V h(e_0,e_1,e_2)\cdot 1 V\\
h(e_0,e_1)\cdot h(e_1,e_3) @>1\cdot h(e_0,e_0,e_0)>\sim> h(e_0,e_1)\cdot h(e_1,e_2)\cdot h(e_2,e_3)\\
\endCD
$$
\medskip
\item{q=4.} Let \ $h^4:I^4@>>>\Delta_4$ \ be a map with boundary values 
$$
h^4(t_1,t_2,t_3,t_4)=
\cases
h(e_0,e_2,e_3,e_4)(t_2,t_3,t_4),\qquad &t_1=0\\
h(e_0,e_1,e_3,e_4)(t_1,t_3,t_4),\qquad &t_2=0\\
h(e_0,e_1,e_2,e_4)(t_1,t_2,t_4),\qquad &t_3=0\\
h(e_0,e_1)(t_4)+h(e_1,e_2,e_3,e_4)(t_2,t_3,t_4),\qquad &t_1=1\\
h(e_0,e_1,e_2)(t_1,t_4)+h(e_2,e_3,e_4)(t_3,t_4),\qquad &t_2=1\\
h(e_0,e_1,e_2,e_3)(t_1,t_2,t_4)+h(e_3,e_4)(t_4),\qquad &t_3=1\\
e_0,\qquad &t_4=0\\
e_4,\qquad &t_4=1\\
\endcases
$$

\bigskip\noindent
{\it Transfer.} We have 
$$|B_0| \times I^q\times \Delta_r@>1\times h^q\times 1>>|B_0| \times \Delta_q\times \Delta_r@>t(Y(E))>>W(E_{qr}),
\qquad 1\le q\le 4$$
Its adjoint with respect to the last coordinate of $I^q$, has the form  
$$|B_0|\times I^{q-1}\times\Delta_r^+@>>>\Omega W(E_{qr})$$ 
We now have
$$|B_0|\times I^{q-1}\times\Delta_r@>>>\Omega W(E_{qr})@>\zeta>>A(E_{qr})@>>>A(E_{qr},E_{qr}\downarrow A_r)$$
Transfer $T_A$ is the induced simplicial tranformation
$$T_A(\bold E):L(B_0,A_0)\wedge(I^{q-1})^+\wedge\Delta_r^+\longrightarrow A(E_{qr},E_{qr}\downarrow A_r),\qquad\qquad 
1\le q\le 4  \tag 17.17$$

\bigskip\noindent
{\it Disk bundles.} We shall define a simplicial homotopy
$$\bold d_A:L\wedge I^+\longrightarrow A \qquad \text{on} \qquad \Cal N(\Cal D) \tag 17.18$$
such that \ $(\bold d_A)_0=d_A$, \ $(\bold d_A)_1= T_A$. Let $E\in \Cal N_r(\Cal D)$,
$$E=(E_0,\Delta_0)@>(\widetilde f_{01},f_{01})>> \ \dots \dots \ 
@>(\widetilde f_{(r-1)r},f_{(r-1)r})>>(E_r,\Delta_r).\qquad E_j=E_j\downarrow (B_j,A_j)$$
and let
$$\overline {E}=
\CD
(\Phi<E_r)\\
@AAA\\
\vdots\\
@AAA\\
(\Phi<E_0)\\
\endCD
$$
By definition, $T(E)=T(\overline E)$. We have \ $Y_k(\overline{E})\in S_1(R(B_0,E_r)$,
$$Y_k(\overline E)(0,1)= B_0\times E_r@>s>>(B_0\times_{B_k}E_k)\sqcup B_0\times E_r@>r>>B_0\times E_r,\qquad r(b_0,e_k)=(b_0,\widetilde f_{kr}(e_k))$$
Let
$$X(E)=B_0\times E_r@>s>>B_0\sqcup B_0\times E_r@>r>>B_0\times E_r,\qquad r(b_0)=(b_0,\Delta_r f_{0r}(b_0))$$
and define \ $Z( E)\in \Cal N_{r+1}(wS_1(R(B_0,E_r))$ \ by
$$Z(E)=X(E)@>>>Y_0(\overline E) @>>> \ \dots \ @>>> Y_r(\overline E)$$
where \ $X(E)@>>>Y_0(\overline E)$ is given by \ $b_0@>>>(b_0,\Delta_0(b_0))$. From (17.15), we have
$$t(Z(E)):|B_0|\times \Delta_1\times \Delta_{r+1} \longrightarrow W(E_r) \tag 17.19$$
Define \ $\psi:\Delta_r\times I@>>>\Delta_{r+1}$ \ by
$$\psi(\lambda_0 e_0 + \dots + \lambda_r e_r,t)=(1-t)e_0+t(\lambda_0 e_1 + \dots + \lambda_r e_{r+1})$$
Let
$$\bold d_A^*=t(Z(E))(!\times 1\times \psi) \tag 17.20$$
Its $\Delta_1$-adjoint defines
$$\bold d_A(\bold E):L(B_0,A_0)\wedge\Delta_r^+\wedge I^+\longrightarrow A(E_r,E_r\downarrow A_r), \qquad (\bold d_A)_0=d_A,\qquad
(\bold d_A)_1=T_A \tag 17.21$$

\bigskip\noindent
{\it Expansions.} Note that $T_A$ can be defined more generally for $qr$-diagrams for which
$$(\widetilde f_{(k-1)k},f_{(k-1)k}):E_{j(k-1)}/_{B_{k-1}} E_{i(k-1)}\longrightarrow E_{jk}/_{B_k} E_{ik},
\qquad 0\le i\le j\le q, \qquad 0\le k\le r$$
is a fiberwise weak equivalence, since (17.13) is still a sequence of weak equivalences under this more general assumption.
We will call such diagrams {\it homotopy excisive}.

Suppose given $E_0<E_1$ and an expansion $(\phi,\lambda)$ as in section 6. Recall  
$$\bold c(\phi,\lambda):A(i_{01})T_A(E_0) \ \sim \ T_A(E_1),\qquad\qquad \overline{\bold c}(\phi,\lambda):T_A(E_0<E_1) \ \sim \ * \tag 17.22$$
Since $\lambda$ is a fiberwise weak equivalence, $T_A(\lambda)$ is defined and
determines a homotopy commutative diagram
$$
\CD
A(i_{01})T_A(E_0)@>T_A(i_{01})>\sim > T_A(E_1)\\
@V \sim V \widehat A(\lambda)(T_A(E_0)\wedge 1) V     @V \sim V 1 V\\
A(\phi)T_A(E_0)@>T_A(\phi)>\sim > T_A(E_1)\\
\endCD 
$$
Therefore 
$$\bold c(\phi,\lambda)=T_A(i_{01}):A(i_{01})T_A(E_0) \ \phom \ T_A(E_1) \tag 17.23$$

We have a homotopy excisive diagram
$$
\CD
\Phi<E<E\\
@AAA\\
\Phi<E_1<E\\
\endCD
$$
where the arrow is inclusion. It gives a homotopy commutative diagram
$$
\CD
T_A(E_1)@>>\sim > A(i_{01})T_A(E_0)+T_A(E_0<E_1)\\
@V \sim V 1 V     @V \sim V T_A(i)+T_A(j_{01}) V\\
T_A(E_1)@>>\sim > T_A(E_1)+T_A(E_1<E_1)\\
\endCD 
$$
where $j_{01}$ is inclusion. Since $T_A(E_1<E_1)=*$, \ and the bottom path is the standard path, it follows from(17.23) that
$$ \overline{\bold c}(\phi,\lambda)=T_A(j_{01}):T_A(E_0<E_1) \ \sim \ * \tag 17.24$$

given an excisive diagram
$$
\Cal E=
\matrix
E_0 & < & E_1\\
\vee & & \vee\\
\widetilde E_0 & < & \widetilde E_1\\
\endmatrix
$$
and expansion $(\phi,\lambda)$ for $E_0<E_1$, which restricts to an expansion $(\widetilde \phi,\widetilde \lambda)$ 
for $\widetilde E_0<\widetilde E_1$, we have a commutative diagram of homotopy excisions
$$ \CD
E_0<E_1 @> j_{01} >> E_1<E_1\\
@AA i A     @AA i' A\\
\widetilde E_0<\widetilde E_1 @>\widetilde j_{01} >> \widetilde E_1<\widetilde E_1\\
\endCD
$$
where $i,i'$ are inclusion. By definition, $T_A(\Cal E)=T_A(i)$. So by (17.24), we have a commutative diagram
$$
\CD
T_A(E_0<E_1) @> \overline{\bold c}(\phi,\lambda) > \phom > *\\
@A \phom A T_{A}(\Cal E) A    @A \phom A 1 A\\
A(i)T_A(\widetilde E_0<\widetilde E_1) @> \overline{\bold c}(\widetilde\phi,\widetilde\lambda) > \phom > *\\
\endCD \tag 17.25
$$

\bigskip\noindent
{\it Permutation paths.} Suppose given an excisive diagram
$$
\Cal E=
\matrix
E_1&< &E\\
\vee && \vee\\
E_0&< &E_2\\
\endmatrix,\qquad\qquad \Cal E=\Cal E\downarrow (B,A)
$$
as in section 6. Since $Y(E_{01})$, $Y(E_{02})$ are in general position, we have
$$H:L\wedge (I)^+\longrightarrow A_0((E/E\downarrow A)\vee (E/E\downarrow A))$$
by
$$H=\left(A(i_{01})(T_A(E_{011})^{-1}\cdot T_A(E_{001}),A(i_{02})(T_A(E_{002})^{-1} \cdot 
T_A(E_{022}))\right)$$
With the notation of section 5, \ $\psi H \ \sim \ \widetilde H$ rel. $\dot I$, and 
$$A_0(\lambda)H \ \sim \ \bold b(\Cal E) \text{ rel. } \dot I \tag 17.26$$

\bigskip
Let $X$ be pointed. There is the augmentation \ $\varepsilon:|X|@>>>A_0(X)$. A fundamental result of Waldhausen \newline [ \ \, \ \ ] asserts that the stabilization of $\varepsilon$,
$$\varepsilon^s:Q(|X|)\longrightarrow A_0^s(X) \tag 17.27$$
is a homotopy equivalence. In particular, $A_0^s$ is a homology theory on $\Cal T_0$, or equivalently $A^s$ is a homology theory on $\Cal T^2$.
Define
$$T_{A^s}:L@>>>A^s,\qquad\qquad L@>T_A>>A@>>>A^s \tag 17.28$$
The properties of $T_A$ discussed above then imply that $T_{A^s}$ is a refined transfer for $A^s$.

\bigskip
\subhead 18. Stable homotopy transfer \endsubhead  Let $F$ be a homology theory. We shall define a refined transfer for $F$. This is a matter of elaborating on the standard construction of transfer as in \cite{2}.

\bigskip
For $\sigma\in\Omega X$, and $0\le a,b\le 1$, define $\sigma^{[a,b]}\in \Omega(X)$  by
$$
\sigma^{[a,b]}(s)=
\cases
*, &\qquad 0\le s\le a\\
(s-a)/(b-a), &\qquad a\le s\le b\\
*, &\qquad b\le s\le 1\\
\endcases
$$

\bigskip
Given a vector bundle $V$, let $V^{\bullet}$ denote its Thom space. Let $\alpha:E@>>>B\times R^s$ be a fiberwise embedding. 
The tangent map
$$T_B(\alpha):T_B(E)\longrightarrow T_B(B\times R^s)=(B\times R^s)\times R^s$$
has the form \ $T_B(\alpha)=(\alpha,\alpha_*)$, \ where \ $\alpha_*:T_B(E)@>>>R^s$.

\medskip
The standard normal bundle of $\alpha$ may be written
$$N_{\alpha}=\{ \ (e,n)\in E\times R^s \ | \ n\perp \alpha_*(T_B(E)_e) \ \}$$
We have
$$\overline\alpha:N_{\alpha}\longrightarrow B\times R^s, \qquad\qquad \overline\alpha(e_b,n)=(b,\alpha_2(e_b)+n) \tag 18.1$$
where $\alpha(e_b)=(b,\alpha_2(e_b))$, and
$$\widetilde\alpha:T_B(E)\oplus N_{\alpha}\longrightarrow E\times R^s,\qquad\qquad \widetilde\alpha(v_e+(e,n))=(e,\alpha_*(v_e)+n) \tag 18.2$$
Let $\phi:N_{\alpha}@>>>B\times R^s$ be an embedding such that $\phi|E=\alpha$. We will
say that $\phi$ is  {\it standard} if $\phi=\overline\alpha$ on a neighborhod of $E$. Let $\Cal E_B^s(E)$ denote the space of such pairs $(\alpha,\phi)$. There is the umkehr map
$$\phi!:B^+\wedge S^s\longrightarrow N_{\alpha}^{\bullet}/(N_{\alpha}\downarrow \dot E)^{\bullet} \tag 18.3$$
Let $i_s:R^s@>>>R^{s+1}$, \ $i_s(x_1,\dots,x_s)=(x_1,\dots,x_s,0)$. Define \ $\Cal E_B^s(E)@>>>\Cal E_B^{s+1}(E)$ \ by
$(\alpha,\phi)@>>>(\widehat\alpha,\widehat\phi)$, \ where $\widehat\alpha$, $\widehat\phi$ are respectively
$$E@>\alpha>>B\times R^s@>1\times i_s>>B\times R^{s+1},\qquad\qquad N_{\widehat\alpha}=N_{\alpha}\oplus R@>\phi\times i_s>>B\times R^{s+1}$$
Then \ $\Cal E_B(E)=\lim_s \Cal E_B^s(E)$ \ is weakly contractible. Note that
$$\widehat\phi!=\phi!\wedge 1:B^+\wedge S^{s+1}\longrightarrow N_{\widehat\alpha}^{\bullet}/(N_{\widehat\alpha}\downarrow \dot E)^{\bullet}=
N_{\alpha}^{\bullet}/(N_{\alpha}\downarrow \dot E)^{\bullet}\wedge S^1 \tag 18.4$$

\medskip
Assume an orthogonal structure for $T_B(E)$. A {\it vertical local vector field} on $E$ is a pair $(U,u)$ where $U$ is an open subset of $E-\dot E$, and $u:U@>>>T_B(E)$ is a section,
such that $u^{-1}(D_r(T_B(E))$ is compact, for all $r>0$. The definition does not depend on the choice of orthogonal structure.

A local vector field corresponds to a section $\overline u$ of the fiberwise one-point compactification of $T_B(E)$, such that 
$\overline u(\dot e)=\infty_{\dot e}$, $\dot e\in\dot E$. The correspondence is \ $\overline u\longrightarrow (U,u)$, where $U=\overline u^{-1}(T_B(E))$, and $u=\overline u|U$. Let $\Cal L_B(E)$ denote the space of vertical local vector fields on $E$ topologized by the above correspondence. 
For a given vector bundle $V$, we will also write
$$u:V^{\bullet}/(V\downarrow \dot E)^{\bullet}\longrightarrow (T_B(E)\oplus V)^{\bullet} \tag 18.5$$
for the map
$$
v_e\longrightarrow 
\cases
u(e)+v_e, \qquad &e\in U\\
\infty,\qquad &e \notin U\\
\endcases
$$

Given  $(\alpha,\phi)\in \Cal E_B^s(E)$ and $(U,u)\in \Cal L_B(E)$, let
$$ \tau(u):(B/A)\wedge S^s\longrightarrow (E/E\downarrow A)\wedge S^s \qquad\qquad \tau(u)=\tau(\alpha,\phi,U,u) \tag 18.6$$
denote the quotient of
$$B^+\wedge S^s@>\phi!>> N_{\alpha}^{\bullet}/(N_{\alpha}\downarrow \dot E)^{\bullet}          
@>u>>(T_B(E)\oplus N_{\alpha})^{\bullet}@>\widetilde\alpha^{\bullet}>>E^+\wedge S^s$$

\bigskip
Let \ $\Cal L^p(E)$ \ denote the space of $p$-tuples of local vector fields on $E$ which are pairwise disjoint. Given \ 
$((U_1,u_1),\dots ,(U_p,u_p))$ \ in $\Cal L^p(E)$, 
$$(U_0,u_0)=(U_1\sqcup\dots \sqcup U_p,u_1\sqcup\dots \sqcup u_p)$$
is also a local vector field on $E$. Note that
$$(\tau(u_1),\dots,\tau(u_p)):B/A\wedge S^s\longrightarrow \prod_1^p (E/E\downarrow A)\wedge S^s$$
factors through \ $\bigvee_1^p (E/E\downarrow A)\wedge S^s$. Define
$$q(u_1,\dots,u_p):L(B,A)\longrightarrow F^s(\bigvee_1^p (E/E\downarrow A),*)  \tag 18.7$$
by
$$\align
L(B,A)@>\varepsilon>>F(B,A)@>\omega_F^*>>\Omega^s F((B/A) \wedge S^s,*)&@>\Omega^s F(\tau)>>
\Omega^sF(\bigvee_1^p (E/E\downarrow A)\wedge S^s,*)\\
&@> \ \ \ \ \ \ \ >>F^s(\bigvee_1^p (E/E\downarrow A),*)\\
\endalign
$$
where $\tau=(\tau(u_1),\dots, \tau(u_p))$.

\bigskip
Given $E_{0\dots p}=E_0<\dots<E_p$, and an embedding $\alpha:E_p@>>>B\times R^s$, let $T_B(E_p)$ have the induced metric.
\smallskip
\noindent 
Let $v$ be vector field on $\overline{E_p-E_0}$ such that 

\hskip1in(i) $v$ points out of $E_p$ on $\dot E_p$.

\hskip1in (ii) $v$ points out of $E_{j-1}$ on $\text{fr}(E_{j-1}<E_j)$, $1<j\le p$.

\noindent
Let $d>0$ be such that $|v|>d$ on 
$$\dot E\cup\text{fr}(E_0<E_1)\cup \dots\cup \text{fr}(E_{k-1}<E_k)$$
Let $\Cal V_{\alpha}(E_{0\dots k})$ denote the space of such pairs $(v,d)$.  Let
$$\Cal M^s(E_{0\dots p})=\{ (\alpha,\phi,v,d) \ | \ (\alpha,\phi)\in \Cal E_B^s(E_p),  \ (v,d)\in \Cal V_{\alpha}(E_{0\dots k}) \} \tag 18.8$$
We have
$$\Cal M^s(E_{0\dots p})@>>>\Cal M^{s+1}(E_{0\dots p}),\qquad\qquad (\alpha,\phi,v,d)@>>>(\widehat\alpha,\widehat\phi,v,d)$$
Let \ $\Cal M(E_{0\dots p})=\lim_s \Cal M^s(E_{0\dots p})$. We see that $\Cal M(E_{0\dots p})$ is weakly contractible.
Since an excision \ $E_{0\dots p}@>>>F_{0\dots p}$ restricts to a fiberwise diffeomorphism \ 
$$\overline{E_p-E_0}@>>>\overline{F_p-F_0}$$
$\Cal M $ is a contravariant functor on $\Cal B^p$. Define
$$\Cal V_{\alpha}(E_{0\dots p})\longrightarrow \Cal L_B^p(E_k)$$ 
by \ $(v,d)@>>>((U_1,u_1),\dots, (U_p,u_p))$, \ where
$$U_j=\{ \ e\in E_j-E_{j-1} \ | \ |v(e)|<d \ \},\qquad\qquad u_j(e)=\frac{1}{d-|v(e)|} v(e)$$
It defines
$$\zeta(E_{0\dots p}):\Cal M(E_{0\dots p})\longrightarrow \Cal E_B(E_p)\times\Cal L_B^p(E_k)\tag 18.9$$
Define
$$\eta (E_{0\dots p}):F^s(E_p,E_p\downarrow A)\longrightarrow \text{Map}(\Cal M(E_{0\dots p})^+,F^s(E_p,E_p\downarrow A)) \tag 18.10$$
by \ $\eta (E_{0\dots p})(x)(y)=x$.  Since $\Cal M(E_{0\dots p})$ is weakly contractible, $\eta$ is a weak homotopy equivalence.

\bigskip
{\it Case} $p=1$. Given $E_0<E_1$, let 
$$t(E_0<E_1):L(B,A)\longrightarrow \text{Map} \, (\Cal M(E_0<E_1)^+,F^s(E_1,E_1\downarrow A))$$
denote the natural transformation
$$L(B,A)@>q>> \text{Map} \, (\Cal N(E_1)^+,F^s(E_1,E_1\downarrow A))@>\zeta^*>>
\text{Map} \, (\Cal M(E_0<E_1)^+,F^s(E_1,E_1\downarrow A))$$
Consider

$$
\CD
   @. F^s(E_1,E_1\downarrow A)\\
@.    @VV \eta(E_0<E_1) V\\
L(B,A) @>t(E_0<E_1)>> \text{Map} \, (\Cal M(E_0<E_1)^+,F^s(E_1,E_1\downarrow A))\\
\endCD \tag 18.11
$$
\medskip\noindent
By the lifting property there is a simplicial transformation 

$$T(E_0<E_1):L(B,A)@>>>F^s(E_1,E_1\downarrow A) \tag 18.12$$
such that $\eta T$ is simplicially homotopic to $t$.

\bigskip
For $\sigma:I@>>> \Omega X$ and $0\le a<b\le 1$, let $\sigma^{[a,b]}:I@>>>\Omega X$ denote
$$
\sigma^{[a,b]}(s)=
\cases
\sigma(0),\qquad &0\le s\le a\\
\sigma((s-a)/(b-a)),\qquad &a\le s\le b\\
\sigma(1),\qquad &b\le s\le 1\\
\endcases
$$
Let $P(I)$ denote the space of paths $\rho:I@>>>I$ such that $\rho(0)=0$ and $\rho(1)=1$. Given $h:A@>>>P(I)$, and $\sigma:I@>>> X$,
we have $h(\sigma):A@>>>\text{Map}(I,X)$ by $h(\sigma)(a)(s)=\sigma(h(a,s))$.

\bigskip
{\it Case} $p=2$. Let $\gamma:I@>>>I$ denote the identity map, and let
$$h:I@>>>P(I)\times P(I),\qquad\qquad  h:(\gamma,\gamma) \ \sim \ (\gamma^{[0,1/2]},\gamma^{[1/2,1]}) \tag 18.13$$
Consider
$$(q(u_1),q(u_2)):L(B,A)\longrightarrow F^s(E_,E_\downarrow A)\times F^s(E_,E_\downarrow A)$$
For 
$$h(q(u_1),q(u_2)):L(B,A)\wedge I^+\longrightarrow F^s(E,E\downarrow A)\times F^s(E,E\downarrow A) \tag 18.14$$
we have
$$h_0(q(u_1),q(u_2))=(q(u_1),q(u_2))\qquad\qquad h_1(q(u_1),q(u_2))=(q(u_1)^{[0,1/2]},q(u_2)^{[1/2,0]})$$
Note that $h_1$ factors through \ $F^s((E/E\downarrow A)\vee (E/E\downarrow A),*)$. We have a natural transformation on $\Cal B$
$$h:L(B,A)\wedge I^+\longrightarrow \text{Map}(\Cal N^2(E)^+,F^s(E,E\downarrow A)\times F^s(E,E\downarrow A)) \tag 18.15$$
Since $F^s$ is a homology theory, we may lift $h$ to a simplicial tranformation
$$\widetilde h:L(B,A)\wedge I^+\longrightarrow \text{Map}(\Cal N^2(E)^+,F^s((E/E\downarrow A)\vee (E/E\downarrow A),*) \tag 18.16$$
such that
$$\widetilde h_0(q(u_1),q(u_2))=q(u_1,u_2)\qquad\qquad \widetilde h_1(q(u_1),q(u_2))=(q(u_1)^{[0,1/2]},q(u_2)^{[1/2,0]})$$
Let
$$t=F^s(\lambda)\widetilde h:L(B,A)\wedge I^+\longrightarrow \text{Map}(\Cal N^2(E)^+,F^s(E,E\downarrow A) \tag 18.17$$
where $\lambda$ is the folding map. We have
$$t(u_1,u_2)_0=q(u_1\sqcup u_2)\qquad\qquad t(u_1,u_2)_1=q(u_1)+q(u_2)$$

\bigskip
Given $E_0<E_1<E_2$, let
$$t(E_0<E_1<E_2):L(B,A)\wedge I^+\longrightarrow \text{Map}(\Cal M(E_0<E_1<E_2)^+,F^s(E_2,E_2\downarrow A) \tag 18.18$$
denote
$$L(B,A)\wedge I^+@>t>> \text{Map}(\Cal N^2(E)^+,F^s(E,E\downarrow A)@>\zeta^*>>\text{Map}(\Cal M(E_0<E_1<E_2)^+,F^s(E_2,E_2\downarrow A)$$
and consider
$$
\CD
   @. F^s(E_2,E_2\downarrow A)\\
@.    @VV \eta(E_0<E_1<E_2) V\\
L(B,A)\wedge I^+ @>t(E_0<E_1<E_2)>> \text{Map} \, (\Cal M(E_0<E_1<E_2)^+,F^s(E_2,E_2\downarrow A))\\
\endCD \tag 18.19
$$
We have that $T(E_0<E_2)$ is a lift of $t(E_0<E_1<E_2)_0$, and $F(i_{12})T(E_0<E_1)+T(E_1<E_2)$ is a lift of $t(E_0<E_1<E_2)_0$.
So we may choose a simplicial transformation
$$T(E_0<E_1<E_2):L(B,A)\wedge I^+\longrightarrow F^s(E_2,E_2\downarrow A) \tag 18.20$$
such that 
$$T(E_0<E_1<E_2)_0=T(E_0<E_2),\qquad T(E_0<E_1<E_2)_1=F(i_{12})T(E_0<E_1)+T(E_1<E_2)$$

\bigskip
Given $h:I@>>>\text{Map}(I,X)$, and $0\le a<b\le 1$, define \ $h_{[a,b]}:I@>>>\text{Map}(I,X)$ \ by \ $h_{[a,b]}(t)=h(t)^{[a,b]}$.
Let 
$$\alpha:(\gamma^{[0,1/4]},\gamma^{[1/4,1/2]},\gamma^{[1/2,1]}) \ \sim \ (\gamma^{[0,1/3]},\gamma^{[1/3,2/3]},\gamma^{[2/3,1]})$$
$$\beta:(\gamma^{[0,1/2]},\gamma^{[1/2,3/4]},\gamma^{[3/4,1]}) \ \sim \ (\gamma^{[0,1/3]},\gamma^{[1/3,2/3]},\gamma^{[2/3,1]})$$
denote the paths such that each coordinate path \ $\gamma^{[a,b]} \ \sim \gamma^{[c,d]}$ \ is given by
$$t\longrightarrow \gamma^{[(1-t)a+tc,(1-t)b+td]}$$
Note that for $(\sigma,\tau,\rho)\in\prod_1^3\Omega X$,  
$$\alpha(\sigma,\tau,\rho):I\longrightarrow \prod_1^3\Omega X$$
factors through \ $\Omega (\bigvee_1^3 X)$. So we have
$$\overline\alpha(\sigma,\tau,\rho)=\Omega(\lambda)\alpha(\sigma,\tau,\rho):I\longrightarrow\Omega X,\qquad\qquad 
\overline\alpha(\sigma,\tau,\rho):(\sigma+\tau)+\rho \ \sim \ \sigma+\tau+\rho$$
Similarly, we have
$$\overline\beta(\sigma,\tau,\rho)=\Omega(\lambda)\beta(\sigma,\tau,\rho):I\longrightarrow\Omega X,\qquad\qquad 
\overline\beta(\sigma,\tau,\rho):\sigma+(\tau+\rho) \ \sim \ \sigma+\tau+\rho$$

\bigskip
{\it Case} $p=3$. We begin with
$$h^2:I^2\longrightarrow P(I)\times P(I)\times P(I)$$
having edge path
$$
\CD
(\gamma,\gamma,\gamma) @>h(\gamma,(\gamma,\gamma))>\sim> (\gamma^{[0,1/2]},\gamma^{[1/2,1]},\gamma^{[1/2,1]})\\
@V\sim V h((\gamma,\gamma),\gamma) V                       @V\sim V (1,h_{[1/2,1]}(\gamma,\gamma))\cdot\beta V\\
(\gamma^{[0,1/2]},\gamma^{[0,1/2]},\gamma^{[1/2,1]})@> (h_{[0,1/2]}(\gamma,\gamma),1)\cdot\alpha >\sim>(\gamma^{[0,1/3]},\gamma^{[1/3,2/3]},\gamma^{[2/3,1]})\\
\endCD
$$

It defines
$$h^2(q(u_1),q(u_2),q(u_3)):L(B,A)\wedge(I^2)^+\longrightarrow \prod_1^3 F(E,E\downarrow A) \tag 18.21$$
The restriction of $h^2$ to $L(B,A)\wedge\partial(I^2)$ has a specified lift to \ $F^s(\bigvee_1^3 E/E\downarrow A,*)$ \ given by
$$
\CD
q(u_1,u_2,u_3) @>\widetilde h(q(u_1),q(u_2,u_3)>\sim> (q(u_1)^{[0,1/2]},q(u_2,u_3)^{[1/2,1]})\\
@V\sim V \widetilde h(q(u_1,u_2),q(u_3)) V                       @V\sim V (1,\widetilde h_{[1/2,1]}(q(u_2),q(u_3)))\cdot \beta V\\
(q(u_1,u_2)^{[0,1/2]},q(u_3)^{[1/2,1]})@> (\widetilde h_{[0,1/2]}(q(u_1),q(u_2)),1)\cdot \alpha >\sim>(q(u_1)^{[0,1/3]},q(u_2)^{[1/3,2/3]},q(u_3)^{[2/3,1]})\\
\endCD
$$
By the lifting property, we have
$$\widetilde h^2(q(u_1),q(u_2),q(u_3)):L(B,A)\wedge(I^2)^+\longrightarrow F^s(\bigvee_1^3 E/E\downarrow A,*) \tag 18.22$$
with edge path as above. Let 
$$t(u_1,u_2,u_3)=F^s(\lambda)\widetilde h^2(q(u_1,q(u_2),q(u_3))$$ 
We now have a simplicial transformation
$$t:L(B,A)\wedge (I^2)^+\longrightarrow \text{Map}(\Cal N^3(E),F^s(E,E\downarrow A) \tag 18.23$$
such that $t(u_1,u_2,u_3)$ has edge path
$$
\CD
q(u_{13}) @> t(u_1,u_{23}) >\sim> q(u_1)+q(u_{23})\\
@V\sim V t(u_{12},u_3) V                       @V\sim V (1+t(u_2,u_3))\cdot \overline\beta V\\
q(u_{12})+q(u_3)@> (t(u_1,u_2)+1)\cdot \overline\alpha >\sim> q(u_1)+q(u_2)+q(u_3)\\
\endCD
$$

\bigskip
Given $E_{0123}=E_0<E_1<E_2<E_3$, again consider
$$
\CD
   @. F^s(E_3,E_3\downarrow A)\\
@.    @VV \eta(E_{0123}) V\\
L(B,A)\wedge (I^2)^+ @>t(E_{0123})>> \text{Map} \, (\Cal M(E_{0123})^+,F^s(E_3,E_3\downarrow A))\\
\endCD \tag 18.24
$$
where $t(E_{0123})$ is
$$L(B,A)\wedge (I^2)^+@>t>> \text{Map}(\Cal N^3(E_3)^+,F^s(E_3,E_3\downarrow A)@>\zeta^*>>\text{Map}(\Cal M(E_{0123})^+,F^s(E_3,E_3\downarrow A))$$
A lift of $t(E_{0123})$ over $L(B,A)\wedge\partial(I^2)^+$ is given by 
$$
\CD
T(E_{03}) @> T(E_{013}) >\sim> F^s(i_{13})T(E_{01})+T(E_{13})\\
@V\sim V T(E_{023}) V                       @V\sim V (1+T(E_{123}))\cdot \overline\beta V\\
F^s(i_{23})T(E_{02})+T(E_{23})@> (T(E_{012})+1)\cdot \overline\alpha >\sim> F^s(i_{13})T(E_{01})+F^s(i_{23})T(E_{12})+T(E_{23})\\
\endCD \tag 18.25
$$
Let
$$T(E_{0123}):L(B,A)\wedge (I^2)^+\longrightarrow F^s(E_3,E_3\downarrow A) \tag 18.26$$
be a simplicial tranformation which is a lift of $t(E_{0123})$, and has edge-path (18.25).

\bigskip
{\it Case} $p=4$. The procedure for defining a simplicial transformation
$$T(E_{01234}):L(B,A)\wedge (I^3)^+\longrightarrow F^s(E_4,E_4\downarrow A) \tag 18.27$$
with boundary values (6.3.3) is similar. We will omit the details.

\bigskip
{\it Disk bundles.} As before, $\Cal D$ is the category of pairs $(E\downarrow B,\Delta)$ such that $E$ is fiberwise diffeomorphic to an orthogonal 
disk bundle, and $\Delta:B@>>>E-\dot E$ is a section. The maps are diffeomorphisms which preserve the section. 

Let $\overline{\Cal E}_B^s(E)$ denote the space of fiberwise embeddings \
$\alpha:E@>>>B\times R^s$ \ having the following property: There is a subbundle $V^*\subset B\times R^s$ such that $\alpha(E)\subset D(V^*)$ 
and $\alpha:E@>>>D(V^*)$ is a section preserving fiberwise diffeomorphism, where $D(V^*)$ has the zero section. We will refer to such 
embeddings as {\it standard}. Let 
$\overline{\Cal E}_B(E)=\lim_s \overline{\Cal E}_B^s(E)$. If $\alpha$, $\beta$ are standard, then for $0\le t\le 1$,
$$\gamma(t)=\cos(\pi t/2)\alpha \oplus \sin(\pi t/2)\beta:E\longrightarrow R^{2s}$$ is also standard. It follows that 
$\overline{\Cal E}_B(E)$ is weakly contractible. 

\medskip
Let $(\widetilde\phi,\phi):E_0\downarrow B_0@>>>E_1\downarrow B_1$ be a fiberwise diffeomorphism preserving the zero section. If 
$\alpha:E_1@>>>B_1\times R^s$ is standard, its pullback 
$$p_0\times_{B_1}\alpha\widetilde\phi:E_0@>>>B_0\times R^s$$ 
is also standard. It results that $\overline{\Cal E}_B$ is a contravariant functor on $\Cal D$.

\medskip
Let $\alpha:E@>>>B\times R^s$ be standard. Let $\rho_{\alpha}$ denote the pullback by $\alpha$ of the radial vector field on $V^*$.
With the metric induced by the tangent map $T_B(\alpha):T_B(E)@>>>T_B(B\times R^s)$, we have $|\rho_{\alpha}(v)|=|\alpha(v)|$. There is the
local vector field $(U,u)$ associated to $(\rho_{\alpha},1)$:
$$U=E-\dot E,\qquad\qquad u(v)=\frac{1}{1-|\alpha(v)|} \rho_{\alpha}(v)$$
Let
$$\overline\tau_{\alpha}(E)=\tau(\alpha,\overline\alpha,U,u):B^+\wedge S^s\longrightarrow E^+\wedge S^s \tag 18.27$$
where $\tau(\alpha,\overline\alpha,U,u)$ is from (18.6). For a vector $v\in B\times R^s$, let $v_2$ denote its projection to $R^s$.
Now $(b,y)\in B\times R^s$, \  may be written uniquely in the form \ $(b,y)=v^*+n$, \ $v^*\in V^*$, $n\in (V^*)^\perp$. Then
$$
\overline\tau_{\alpha}(E)(b,y)=
\cases
( \alpha^{-1}(v^*),\frac{1}{1-|v^*|} v^*_2+n_2),\quad &|v^*|<1\\
*,\qquad &\text{otherwise}\\
\endcases 
$$
For $(\widetilde\phi,\phi):E_0\downarrow B_0@>>>E_1\downarrow B_1$, 
$$(\widetilde\phi\wedge 1)\overline\tau_{[\alpha\widetilde\phi]}(E_0)=\overline\tau_{\alpha}(E_1)(\phi\wedge 1),\qquad
[\alpha\widetilde\phi]:=p_0\times_{B_1} \alpha\widetilde\phi \tag 18.28$$

Define
$$\delta_{\alpha}(E):B^+\wedge S^s\wedge I^+\longrightarrow E^+\wedge S^s \tag 18.29$$
by
$$
\delta_{\alpha}(E)(b,y)(t)=
\cases
( \alpha^{-1}(tv^*),\frac{1}{1-t|v^*|} v^*_2+n_2),\quad &t|v^*|<1\\
*,\qquad &\text{otherwise}\\
\endcases 
$$
We have \ $\delta_{\alpha}(E)_0=\Delta\wedge 1$, \  $\delta_{\alpha}(E)_1=\overline\tau_{\alpha}(E)$.
For $(\widetilde\phi,\phi):E_0\downarrow B_0@>>>E_1\downarrow B_1$, 
$$(\widetilde\phi\wedge 1)\delta_{[\alpha\widetilde\phi]}(E_0)=\delta_{\alpha}(E_1)(\phi\wedge 1 \wedge 1) \tag 18.30$$

\medskip
From $\overline\tau_{\alpha}(E)$ we have
$$\overline t_{\alpha}(E):L(B,A)\longrightarrow F^s(E,E\downarrow A) \tag 18.31$$
by
$$L(B,A)@>\varepsilon>>F(B,A)@>\omega_F^*>>\Omega^s F((B/A) \wedge S^s,*)@>\Omega^s F(\overline\tau)>>
\Omega^sF((E/E\downarrow A)\wedge S^s,*)@>>>F^s(E,E\downarrow A))
$$
where $\overline\tau=\overline\tau_{\alpha}(E)$. Similarly, from $\delta_{\alpha}(E)$ we have
$$d_{\alpha}(E):L(B,A)\wedge I^+\longrightarrow F^s(E,E\downarrow A) \tag 18.32$$
In view of (11.22) and (11.24), these maps define natural transformations on $\Cal D$
$$\overline t(E):L(B,A)\longrightarrow \text{Map}(\overline{\Cal E}_B(E),F^s(E,E\downarrow A)) \tag 18.33$$
$$d(E):L(B,A)\wedge I^+\longrightarrow \text{Map}(\overline{\Cal E}_B(E),F^s(E,E\downarrow A)) \tag 18.34$$

\medskip
Define \ $\mu:\overline{\Cal E}_B(E)@>>>\Cal M(E)$ \ by \ $\alpha@>>>(\alpha,\overline\alpha,\rho_{\alpha},1)$.
We have
$$
\CD
   @. F^s(E, E\downarrow A)) @.\\
@.    @VV \eta(E) V     @.\\
L(B,A) @>t(E)>> \text{Map} \, (\Cal M(E)^+,F^s(E,E\downarrow A)
@>\mu^*>> \text{Map} \, (\overline{\Cal E}_B(E)^+,F^s(E,E\downarrow A)\\
\endCD 
$$
\medskip\noindent
and $\overline t(E)=\mu^*t(E)$. Since by definition, transfer is a lift of $t(E)$
it is also a lift of $\overline t(E)$. It follows that $d(E)$ can be lifted to a simplicial tranformation
$$\bold d(E,\Delta):L(B,A)\wedge I^+\longrightarrow F^s(E,E\downarrow A) \tag 18.35$$
such that \ $\bold d(E,\Delta)_0=F^s(\Delta)\varepsilon(B,A)$, \ $\bold d(D(V))_1=T(E,E\downarrow A)$.

\bigskip
This completes the definition of the transformations that constitute a transfer $T_{F^s}$ for $F^s$. It is straightforward to verify the normality properties. Finally, since $F$ is a homology theory, \ $F@>>>F^s$ \ is a weak equivalence, so $T_{F^s}$ can be pulled back to a transfer
for $F$.

\bigskip
\subhead Appendix  \endsubhead Let $\Cal C$ be a category of the form
$$
\CD
C_1 @>c_1>> C\\
@A c_{01} AA  @A c_2 AA\\
C_0 @>c_{02}>> C_2\\
\endCD
$$
and let $c_0=c_1c_{01}=c_2c_{02}$. Let  $F$, $G$ functors from $\Cal C$ to pointed spaces. Assume that
$$
\CD
F(C_1) @>F(c_1)>> F(C)\\
@A F(c_{01}) AA  @A F(c_2) AA\\
F(C_0) @>F(c_{02})>> F(C_2)\\
\endCD \tag *
$$
is a pushout. Let $(B,A)$ be a pointed pair and $S,\widehat S:F\wedge B\longrightarrow G$
simplicial transformations such that $S|F\wedge A=\widehat S|F\wedge A$. 
\proclaim{ Lemma} Suppose that for $j=1,2$ we have on the subcategory
$\{C_0@>c_{0j}>>C_j\}$ a simplicial homotopy
$$H^j:F\wedge B\wedge I^+\longrightarrow G$$
rel. $F\wedge A$ such that $H_0^j=S$ and $H_1^j=\widehat S$. Suppose further that we have a homotopy
$$P:H^1(C_0) \sim  H^2(C_0)\quad  rel.  \quad F(C_0)\wedge (B\wedge \dot I^+\cup A\wedge I^+)$$
Then $S(C) \sim \widehat S(C)$ rel. $F(C)\wedge A$.
\endproclaim
\demo{Proof} Define
$$K^j:F(C_j)\wedge B\wedge I^+\longrightarrow G(C)$$
with
$$K_0^j=S(C)(F(c_j)\wedge 1),\qquad K_1^j=\widehat S(C)(F(c_j)\wedge 1)$$
by
$$
K^j(x\wedge b\wedge s)=
\cases
S(C_j@>>>C)(x\wedge b\wedge 1-3s), &0\le s\le 1/3\\
G(c_j)H^j(C_j)(x\wedge b\wedge 3s-1), &1/3\le s\le 2/3\\
\widehat S(C_j@>>>C)(x\wedge b\wedge 3s-2), &2/3\le s\le 1\\
\endcases
$$
A homotopy
$$L:F(C_0)\wedge B\wedge I^+\wedge I^+\longrightarrow G(C)\qquad\qquad \text{rel. } F(C_0)\wedge B\wedge \dot I^+$$
such that
$$L_0=K^1(F(c_{01})\wedge 1),\qquad L_1=K^2(F(c_{02})\wedge 1)$$
is constructed as follows. Consider
\bigskip
{\eightpoint
$$\CD
S(C)F(c_0)@>\alpha_1^{-1}>\sim>G(c_1)S(C_1)F(c_{01})@>\beta_1>\sim>
G(c_1)\widehat S(C_1)F(c_{01})@>\widehat\alpha_1 >\sim>\widehat S(C)F(c_0)\\
@V 1 V\sim V  @V \gamma_1^{-1} V\sim V  @V \widehat\gamma_1^{-1} V\sim V @V 1 V\sim V\\
S(C)F(c_0)@>\alpha_0^{-1}>\sim>G(c_0)S(C_0)@>\beta_{01}>\sim>
G(c_0)\widehat S(C_0)@>\widehat\alpha_0 >\sim>\widehat S(C)F(c_0)\\
@V 1 V\sim V  @V 1 V\sim V  @V 1 V\sim V @V 1 V\sim V\\
S(C)F(c_0)@>\alpha_0^{-1}>\sim>G(c_0)S(C_0)@>\beta_{02}>\sim>
G(c_0)\widehat S(C_0)@>\widehat\alpha_0 >\sim>\widehat S(C)F(c_0)\\
@V 1 V\sim V  @V \gamma_2 V\sim V  @V \widehat\gamma_2 V\sim V @V 1 V\sim V\\
S(C)F(c_0)@>\alpha_2^{-1}>\sim>G(c_2)S(C_2)F(c_{02})@>\beta_2>\sim>
G(c_2)\widehat S(C_2)F(c_{02})@>\widehat\alpha_2 >\sim>\widehat S(C)F(c_0)\\
\endCD
$$}
where
$$
\alignat 3
\alpha_j&=S(C_j\rightarrow C)F(c_{01}),\qquad &\widehat\alpha_j&=\widehat S(C_j\rightarrow C)F(c_{01}),\qquad
&\beta_j&=G(c_j)H^j(C_j)F(c_{01})\\
\beta_{0j}&=G(c_0)H^j(C_0),\qquad  &\gamma_j&=G(c_j)S(C_0\rightarrow C_j),\qquad 
&\widehat\gamma_j&=G(c_j)\widehat S(C_0\rightarrow C_j)\\
\endalignat
$$

Now $L$ is constructed by filling in the squares with the following:

$$
\matrix
S(C_0\rightarrow C_1\rightarrow C)\qquad &G(c_1)H^1(C_1)F(c_{01})\qquad &\widehat S(C_0\rightarrow C_1\rightarrow C)\\
S(C_0\rightarrow C)\qquad\hfill &G(c_0)P\qquad\hfill &\widehat S(C_0\rightarrow C)\hfill\\
S(C_0\rightarrow C_2\rightarrow C)\qquad &G(c_2)H^1(C_2)F(c_{02})\qquad &\widehat S(C_0\rightarrow C_2\rightarrow C)\\
\endmatrix
$$

The restriction of $K^j$ to $F(C_1)\wedge A\wedge I^+$ has a canonical deformation to the constant homotopy.
Similarly for the restriction of $L$ to $F(C_0)\wedge A\wedge I^+\wedge I^+$. These deformations fit together
allowing the construction of 
$$\widetilde K^j:F(C_j)\wedge B\wedge I^+\longrightarrow G(C)\qquad \text{rel. } F(C_j)\wedge A $$
such that
$$\widetilde K_0^j=S(C)(F(c_j)\wedge 1),\qquad \widetilde K_1^j=\widehat S(C)(F(c_j)\wedge 1)$$
and
$$\widetilde L:F(C_0)\wedge B\wedge I^+\wedge I^+\longrightarrow G(C)\qquad \text{rel. } 
F(C_0)\wedge (B\wedge \dot I^+\cup A\wedge I)$$
such that
$$\widetilde L_0=\widetilde K^1(F(c_{01})\wedge 1),\qquad \widetilde L_1=\widetilde K^2(F(c_{02})\wedge 1)$$
Since ($*$) is a pushout, $\widetilde K^1$, $\widetilde K^2$, and $\widetilde L$ determine a homotopy
$S(C)\sim \widehat S(C)$ rel. $F(C)\wedge A$.
\enddemo


\Refs

\ref
\no
\by
\pages
\paper
\yr
\vol
\jour
\endref

\ref
\no 1
\by B. Badzioch, W. Dorabiala, B. Williams
\pages 660--680
\paper Smooth parametrized torsion -- a manifold approach
\yr 2009
\vol 221 (2)
\jour Adv. Math.
\endref

\ref
\no 2
\by J. C. Becker and D. H. Gottlieb
\pages 1--12
\paper The transfer map and fiber bundles
\yr 1975
\vol 14
\jour Topology
\endref


\ref
\no 3
\by J. C. Becker and R. E. Schultz
\pages 583--605
\paper Axioms for bundle transfers and traces
\yr 1998
\vol 227
\jour Math. Z.
\endref


\ref
\no 4
\by W. Dwyer, M. Weiss, B. Williams
\pages 1--104
\paper A parametrized index theorem for the algebraic K-theory Euler class
\yr 2003
\vol 190(1)
\jour Acta Math.
\endref

\ref
\no 5
\by M. Hirsch
\pages
\paper Differential Topology
\yr 1976
\vol 
\jour Graduate Texts in Mathematics Vol. 33, Springer, Berlin-Heidelberg-New York
\endref

\ref
\no 6
\by V. G. Turaev
\pages 627--662
\paper Euler structures, nonsingular vector fields, and torsions of Reidemeister type
\yr 1990
\vol 34
\jour Math. USSR Izv.
\endref

\ref
\no 7
\by V. G. Turaev\pages
\paper Introduction to Combinatorial Torsions
\yr 2001
\vol
\jour Birkh\"auser Basel
\endref

\ref
\no 8
\by F. Waldhausen
\pages 318--419
\paper Algebraic K-theory of spaces
\yr 1985
\vol 1126
\jour Springer Lecture Notes in Mathematics
\endref

\ref
\no 9
\by F. Waldhausen
\pages 141--184
\paper Algebraic K-theory of spaces, a manifold approach
\yr 1982
\vol 2
\jour CMS Conference Proceedings 
\endref

\ref
\no 10
\by F. Waldhausen
\pages 392--417
\paper Algebraic K-theory of spaces, concordance, and stable homotopy theory
\yr 1983
\vol 113
\jour Algebraic topology and algebraic K-theory, Ann. of Math. Stud.
\endref






\endRefs
\end{document}